\documentclass[english]{article}
\usepackage[T1]{fontenc}
\usepackage[latin9]{inputenc}
\usepackage{geometry}
\usepackage{babel}
\usepackage{amsmath}
\usepackage{amsthm}
\usepackage{amssymb}
\usepackage{graphicx}
\usepackage{setspace}
\PassOptionsToPackage{normalem}{ulem}
\usepackage{ulem}
\usepackage[unicode=true,pdfusetitle,
 bookmarks=true,bookmarksnumbered=false,bookmarksopen=false,
 breaklinks=false,pdfborder={0 0 0},pdfborderstyle={},backref=false,colorlinks=false]
 {hyperref}

\makeatletter
\theoremstyle{plain}
\newtheorem{thm}{\protect\theoremname}[section]
\theoremstyle{definition}
\newtheorem{defn}[thm]{\protect\definitionname}
\theoremstyle{remark}
\newtheorem{claim}[thm]{\protect\claimname}
\ifx\proof\undefined
\newenvironment{proof}[1][\protect\proofname]{\par
	\normalfont\topsep6\p@\@plus6\p@\relax
	\trivlist
	\itemindent\parindent
	\item[\hskip\labelsep\scshape #1]\ignorespaces
}{%
	\endtrivlist\@endpefalse
}
\providecommand{\proofname}{Proof}
\fi
\theoremstyle{definition}
\newtheorem{example}[thm]{\protect\examplename}
\theoremstyle{plain}
\newtheorem{cor}[thm]{\protect\corollaryname}
\theoremstyle{remark}
\newtheorem{rem}[thm]{\protect\remarkname}
\theoremstyle{plain}
\newtheorem{lem}[thm]{\protect\lemmaname}
\theoremstyle{plain}
\newtheorem{prop}[thm]{\protect\propositionname}

\usepackage{tikz}
\usetikzlibrary{positioning}
\usetikzlibrary{arrows.meta}

\usepackage{algorithm}
\usepackage{algorithmic}

\usepackage{hyperref}
\hypersetup{hidelinks,
backref=true,
pagebackref=true,
hyperindex=true,
breaklinks=true,
colorlinks=true,
linkcolor=[RGB]{164,65,255},
urlcolor=blue,
citecolor=[RGB]{164,65,255},
bookmarks=true,
bookmarksopen=false,
pdftitle={Title},
pdfauthor={Author}}

\makeatother

\providecommand{\claimname}{Claim}
\providecommand{\corollaryname}{Corollary}
\providecommand{\definitionname}{Definition}
\providecommand{\examplename}{Example}
\providecommand{\lemmaname}{Lemma}
\providecommand{\propositionname}{Proposition}
\providecommand{\remarkname}{Remark}
\providecommand{\theoremname}{Theorem}

\begin{document}
\title{Exposure Orders in Free Group Algebras: Minimal Schreier Transversals,
Free Bases, and Gröbner Bases}
\author{Matan Seidel}
\maketitle
\begin{abstract}
Consider the free group algebra $K\left[F\right]$, where $F$ is
a free group and $K$ a field. A well-order $\prec$ on $F$ is called
an \emph{exposure order} if words are greater than their proper prefixes.
We show that every one-sided ideal $I$ in $K\left[F\right]$ admits
a Schreier transversal, a basis, and a Gröbner basis -- each minimal
in a natural sense with respect to $\prec$. When $I$ is finitely
generated and $\prec$ is computable, we provide an algorithm for
computing these minimal structures from a finite generating set. This
extends the foundational works of Lewin \cite{Lewin1969} and Rosenmann
\cite{Rosenmann1993}, which relied on the \emph{shortlex} order,
to both a broader class of orders on $F$ and to infinitely generated
ideals, while retaining algorithmic capabilities for such orders in
case $I$ is finitely generated.

General exposure orders lack a form of compatibility with products
which we call \emph{suffix-invariance}, that shortlex enjoys, and
which prior Gröbner basis constructions in $K\left[F\right]$ relied
on. In its absence, reductions may strictly increase the support of
elements, requiring nontrivial conceptual adaptations to definitions
and algorithms. These adaptations clarify the notion of minimality
underlying prior constructions and demonstrate that algorithmic Gröbner
theory in $K\left[F\right]$ does not fundamentally require suffix-invariance,
although its presence---as in shortlex---results in a simpler theory.

Our framework further illuminates the flexibility of exposure orders:
with a suitable choice of $\prec$, any Schreier transversal for $I$
can be realized as minimal, and any basis for $I$ arising from the
constructions of Lewin or Rosenmann can likewise be realized as its
minimal basis, unifying both approaches under a single framework.
\end{abstract}
\tableofcontents{}

\section{\label{sec: Introduction}Introduction}

Fix a field $K$ and a free group $F$ with basis $S=\left\{ x_{1},x_{2},...,x_{r}\right\} $.
The free group algebra $\mathcal{A}:=K\left[F\right]$ is the $K$-vector
space with basis $F$, endowed with multiplication obtained by bilinearly
extending the group operation: $\left(\sum_{u\in F}\alpha_{u}u\right)\cdot\left(\sum_{v\in F}\beta_{v}v\right)=\sum_{u,v\in F}\alpha_{u}\beta_{v}uv$.

A classical theorem of Cohn and Lewin \cite{Cohn1964,Lewin1969} asserts
that $\mathcal{A}$ is a \emph{free ideal ring}: every one-sided ideal
of $\mathcal{A}$ is a free $\mathcal{A}$-module, with a well-defined
rank equal to the cardinality of any of its bases.\footnote{Theorem \ref{thm: The exposure basis is a basis} provides another
proof of the freeness of one-sided ideals of $\mathcal{A}$. For a
short and elegant proof, see Hog-Angeloni \cite{HogAngeloni2006}.} This result parallels the Nielsen-Schreier theorem for free groups
\cite{Nielsen1921,Schreier1927}, drawing an analogy between the one-sided
ideals of $\mathcal{A}$ and subgroups of $F$. Moreover, a subgroup
$H\leq F$ can be encoded as the right ideal $I_{H}\leq\mathcal{A}$
generated by $\left\{ h-1:h\in H\right\} $. If $B\subseteq H$ is
a basis for $H$ then $\left\{ w-1:w\in B\right\} $ forms a basis
for $I_{H}$, and hence $\text{rk}H=\text{rk}I_{H}$ (see \cite[Section~4]{Cohen2006}).

This analogy underlies recent approaches that use $\mathcal{A}$ to
study $F$. For example, Jaikin-Zapirain \cite{JaikinZapirain2024}
used $\mathcal{A}$ to address the  inert--compressed conjecture
for free groups. When $K$ is finite, a paper on the theory of word
measures on $\text{GL}_{n}\left(K\right)$ \cite{ErnstWest2024} by
Ernst-West, Puder, and the current author, uncovered their natural
connection to certain ideals in $\mathcal{A}$. Motivated by this
relevance of $\mathcal{A}$ to the study of free groups, and in particular,
by the need to analyze concrete examples in the context of word measures,
the present paper joins prior work on primitivity testing \cite{Seidel2025}
in the goal of further developing theoretical and algorithmic aspects
of $\mathcal{A}$.

Inspired by the foundational works of Lewin \cite{Lewin1969} and
Rosenmann \cite{Rosenmann1993} on bases for right ideals in $\mathcal{A}$,
we extend their work beyond the standard \emph{shortlex} order (see
Example \ref{exa: Examples of Exposure orders}) to a broader class
of well-orders on $F$, which we call \emph{exposure orders} (see
Definition \ref{def: exposure orders}). In this context, it is somewhat
unnatural that prior constructions and algorithms rely on a specific
order, including the approaches of Lewin and Rosenmann, as well as
parallel developments in the Gröbner-basis theory of $\mathcal{A}$
by Madlener and Reinert \cite{Madlener1993,Reinert1995}. While shortlex
does enjoy a form of compatibility with products, which we call \emph{suffix-invariance}
(see Definition \ref{def: suffix preserving}), working with more
general orders allows us to isolate the essential principles behind
these constructions while gaining flexibility. 

Given an exposure order $\prec\text{ }$on $F$ and a one-sided ideal
$I$ in $\mathcal{A}$, we define a Schreier transversal, a basis,
and a Gröbner basis, all of which are minimal in a natural sense with
respect to $\prec$. Defining these structures for general exposure
orders requires nontrivial adaptations of prior constructions and
algorithms, since prefix-based reductions are no longer guaranteed
to decrease support. We address these challenges in two ways. First,
the structures are constructed directly via their minimality, the
bases by means of transfinite recursion and the minimal Schreier transversal
by introducing an auxiliary order that produces a genuine minimum.
This approach also allows definitions to extend to infinitely-generated
ideals. Second, we develop more delicate arguments to ensure algorithmic
termination, ensuring that these constructions are not only theoretically
well-defined but also effectively computable in the finitely-generated
case.

\subsection{Ordering Supports of $\mathcal{A}$ via Orders on $F$}

It is standard in the combinatorial study of $\mathcal{A}$ to impose
an order $\prec$ on the base group $F$, with the goal of constructing
minimal structures and performing algorithmic reductions. The order
$\prec$ is typically extended to an order $\prec_{\text{max}}$ on
finite subsets of $F$, by examining the maximal word in which two
subsets differ (see Definition \ref{def: order on finite subsets}).
For $f\in\mathcal{A}$, its \emph{support} is the set of words with
nonzero coefficient in $f$. When $\prec$ is a well-order, the induced
order $\prec_{\text{max}}$ is also a well-order (see Claim \ref{claim: The order on finite subsetes is a well-order}).
This allows one to select elements with $\prec_{\text{max}}$-minimal
support from any nonempty subset of $\mathcal{A}$.
\begin{defn}
\textbf{\label{def: exposure orders}(Exposure Orders)} Let $\prec$
be a well-order on $F$. The order is called an \emph{exposure order}
if $u\preceq v$ whenever $u$ is a prefix of $v$.
\end{defn}
A main motivation for using exposure orders is their flexibility.
This feature has already proven useful to us in our earlier work \cite[Section~3]{ErnstWest2024},
where an ideal $I$ is ``explored'' along a fixed exposure order
on $F$, producing a basis for $I$. In that process, elements of
$I$ with $\prec_{\text{max}}$-smaller support are ``discovered''
earlier, so varying the exposure order yields distinct bases. This
flexibility made it possible to transfer desired properties from one
basis of $I$ to another (see \cite[Thm.~3.8]{ErnstWest2024}).

\subsection{Minimal Schreier Transversals}

Lewin \cite{Lewin1969} studied bases of a right ideals $I\leq\mathcal{A}$
by adapting Schreier\textquoteright s method from free groups. To
construct such bases, he examined the \emph{Schreier transversals}
of $I$: prefix-closed subsets of $F$ whose $K$-span contains exactly
one representative from each $I$-coset in $\mathcal{A}$. Lewin showed
that, relative to a shortlex order $\prec$ imposed on $F$, every
right ideal $I\leq\mathcal{A}$ admits a ``minimal'' Schreier transversal,
consisting of those words in $F$ whose support is $\prec_{\text{max}}$-minimal
within their $I$-coset. 

Let $\prec$ now denote an exposure order on $F$. A naive adaptation
of Lewin's definition of a minimal Schreier transversal may fail to
produce a Schreier transversal, since the resulting subset of $F$
need not be prefix-closed (see Example \ref{exa: T_I' not closed under prefixs}).
To overcome this, we introduce an auxiliary order $\prec_{\text{min}}$
on subsets of $F$, which is well-suited for comparing Schreier transversals
(see Section \ref{sec: Ordering the Subsets of F}). Using $\prec_{\text{min}}$,
our first main result establishes that a well-defined minimal Schreier
transversal can be selected.
\begin{thm}
\label{thm: minimal Schreier transversal exists}Let $I$ be a right
ideal in $\mathcal{A}$. The set of Schreier transversals for $I$
admits a minimum $T_{I,\prec}$ with respect to the order $\prec_{\text{min}}$.
\end{thm}
When $\prec$ is the shortlex order, the minimal Schreier transversal
$T_{I,\prec}$ coincides with Lewin's ``minimal'' Schreier transversal
(see Theorem \ref{thm: suffix preserving order main theorem}), thereby
realizing it as a genuine minimum. Moreover, for every Schreier transversal
$T$ for $I$, there exists an exposure order $\prec_{T}$ on $F$
such that $T=T_{I,\prec_{T}}$ (see Proposition \ref{prop: Every Schreier tranversal is minimal with respect to some exposure order}).
This highlights the flexibility of exposure orders: any Schreier transversal
can be realized as minimal by choosing a suitable exposure order.

\subsection{Minimal Bases}

Rosenmann \cite{Rosenmann1993} extended Lewin's ideas to construct
a minimal basis for finitely generated right ideals of $\mathcal{A}$
(recall that such ideals are free as $\mathcal{A}$-modules, and thus
admit a basis). Fixing a similar shortlex order\footnote{For words of equal length, Rosenmann \cite{Rosenmann1993} employs
a lexicographical order and Lewin \cite{Lewin1969} a reverse lexicographical
one.} $\prec$ on $F$, he combined Lewin's minimal Schreier transversal
$T$ with a Gröbner-style reduction framework to develop an algorithm
that transforms a finite generating set for a right ideal $I\leq\mathcal{A}$
into a basis that is minimal with respect to $\prec$. Here, minimality
refers to the lexicographical order induced by $\prec_{\text{max}}$
on $n$-tuples of supports, where $n=\text{rk}I$. If the basis elements
are further required to be monic, the basis is unique, and we call
it the \emph{Rosenmann basis }of $I$ (with respect to the fixed order
$\prec$).\footnote{The shortlex order $\prec$ depends on an arbitrary ordering of the
generators of $F$ and their inverses; this is its only degree of
freedom.}

We extend Rosenmann's construction in two directions: to general exposure
orders and to possibly infinitely generated ideals. A naive adaptation
is not straightforward. Defined as the output of an algorithm, Rosenmann's
construction is ill-suited for infinitely generated ideals. Moreover,
under general exposure orders, the reduction steps may strictly $\prec_{\text{max}}$-increase
the support, preventing a usual reliance on well-ordering to guarantee
termination.\footnote{Specifically, the ``div-reductions'' of Rosenmann \cite{Rosenmann1993},
or ``prefix-reductions'' of Madlener and Reinert \cite{Madlener1993}
-- see Example \ref{exa: word smaller than its phi}} We resolve these issues by defining the basis directly via its minimality,
using the minimal Schreier transversal of Theorem \ref{thm: minimal Schreier transversal exists}. 

Explicitly, given an exposure order $\prec$ on $F$ and a right ideal
$I\leq\mathcal{A}$, we construct the \emph{exposure basis} $B_{I,\prec}$
by transfinite recursion: starting with the empty set, at each stage
we add a monic element of $I$ whose support is $\prec_{\text{max}}$-minimal
among those supported on the current minimal Schreier transversal
(See Definition \ref{def: inductive definition for ideal}).
\begin{thm}
\label{thm: The exposure basis is a basis}Let $\prec$ be an exposure
order on $F$ and let $I\leq\mathcal{A}$ be a right ideal. The exposure
basis $B_{I,\prec}$ is a basis for $I$.
\end{thm}
Enjoying the flexibility of exposure orders, exposure bases unify
under a single framework the distinct ways in which Rosenmann and
Lewin constructed bases for $I$.\footnote{The Rosenmann and Lewin bases of $I$ may differ, even when the Lewin
basis is induced by the shortlex-minimal Schreier transversal of $I$.
This is because the Lewin basis prioritizes positive labels rather
than the minimality of the resulting basis.} Lewin used a Schreier transversal $T$ and showed \cite[Theroem~1]{Lewin1969}
that it induces a basis for $I$ by considering the ``positive''
(i.e., $S$-labeled) outgoing edges from $T$. We show that every
such Lewin basis can be realized as the exposure basis of $I$ for
a suitably chosen exposure order (see Proposition \ref{prop: Lewin bases are exposure bases}).
Moreover, when $\prec$ is the shortlex order, the exposure basis
$B_{I,\prec}$ coincides with the Rosenmann basis (Theorem \ref{thm: suffix preserving order main theorem}).

\subsection{Gröbner Bases}

The theory of Gröbner bases systematically transforms generating sets
of ideals, through a series of reductions, into special ones called
\emph{Gröbner bases}, which facilitate computations such as ideal
membership and normal forms. Originating in the works of Buchberger
\cite{Buchberger1965,Buchberger1970,Buchberger1985} on polynomial
rings, the theory was adapted to free associative algebras by Mora
\cite{Mora1985}, and to one-sided ideals of free group algebras by
Rosenmann \cite{Rosenmann1993} and of general monoid rings by Madlener
and Reinert \cite{Madlener1993}, with a subsequent specialization
to free group algebras by Reinert \cite{Reinert1995}. 

Gröbner base theory traditionally imposes two requirements on the
order $\prec$ used. The first is well-ordering, to ensure algorithmic
termination by guaranteeing that only finitely many reductions which
strictly $\prec_{\text{max}}$-decrease support can be performed on
a given polynomial. The second is some form of compatibility with
products, known as admissibility. In the free group algebra, all works
above use the shortlex order, relying on a form of one-sided admissibility
which we call\emph{ suffix}-\emph{invariance} (see Definition \ref{def: suffix preserving}).
Recently, Ceria, Mora and Roggero \cite{Ceria2019} advocated for
the removal of the admissibility condition, providing both motivation
and a resulting theory for commutative polynomial rings. Our work
aligns with this perspective: the class of exposure orders, while
including all suffix-invariant orders (see Claim \ref{claim: Suffix-invariant orders are exposure orders}),
does not guarantee compatibility with products. For instance, it might
happen that $e\prec x$ while $y\succ xy$.

Let $\prec$ be an exposure order. For any nonzero $f\in\mathcal{A}$,
the\emph{ head term} of $f$, denoted $\text{HT}_{\prec}\left(f\right)$,
is the $\prec$-maximal word in $\text{supp}\left(f\right)$. We introduce
the notion of \emph{combinatorially reducing systems} \emph{(CRSs)}
for proper right ideals $I\lneqq\mathcal{A}$ (with respect to $\prec$,
see Definition \ref{def: combinatorially reducing systems}) -- a
special type of Gröbner basis for $I$, tailored to the free group
algebra, inspired by the combinatorial properties satisfied by the
Gröbner bases from Rosenmann's algorithm \cite{Rosenmann1993}.

A CRS for $I$ is defined by combinatorial conditions that simultaneously
enable algorithmic reduction to a normal form (see Theorem \ref{thm: division with remainder}),
while preventing redundancy among its elements (e.g., distinct head
terms). In particular, for a finitely generated $I$, a CRS contains
exactly $2\text{rk}I$ elements (see Remark \ref{rem: cardinality of combinatorially reducing system}).
The defining conditions for a CRS are more relaxed than Reinert's
\emph{monic reduced prefix standard Gröbner basis} \cite[~Thm.~4.4.12]{Reinert1995}
in the sense that the head term of one element may lie in the support
of another, provided that the head terms of both elements end in the
same letter. This flexibility allows multiple CRSs for a single ideal,
and importantly, ones which remain stable under the exposure process
of Definition \ref{def: inductive definition for ideal}. When $\prec$
is shortlex, the Rosenmann's algorithm \cite{Rosenmann1993}, as well
as Reinert's reduced Gröbner basis above, and her $\mathcal{A}$-tailored
specialization \cite[~Thm.~5.2.8]{Reinert1995} are all special cases
of a CRS (see Remark \ref{rem: why we define seconds this way} for
further discussion).

CRSs facilitate a canonical division with remainder in a Schreier
transversal (compare with Reinert's standard representation \cite{Reinert1995}).
\begin{thm}
\textbf{(Division with Remainder in a Schreier Transversal) }\label{thm: division with remainder}Let
$\mathcal{Q}$ be a combinatorially reducing system for a right ideal
$I\lneqq\mathcal{A}$ with respect to an exposure order $\prec$ on
$F$, and let $T$ be a Schreier transversal for $I$. Every $f\in\mathcal{A}$
admits a unique decomposition $f=\sum_{q\in\mathcal{Q}}qg_{q}+r$,
where $r$ is supported on $T$ and each $g_{q}$ is supported on
words such that there is no cancellation in the product $\text{HT}_{\prec}\left(q\right)\cdot g_{q}$.
\end{thm}
Moreover, recall that Theorem \ref{thm: minimal Schreier transversal exists}
provides a minimal Schreier transversal $T_{I,\prec}$ for $I$ with
respect to $\prec$. While various CRSs for $I$ may exists, they
all determine $T_{I,\prec}$ via their prefix structure.
\begin{thm}
\label{thm: combinatorially reducing system defines minimal schreier transversal of the ideal which it generates}Let
$\mathcal{Q}$ be a combinatorially reducing system for a right ideal
$I\lneqq\mathcal{A}$ with respect to an exposure order $\prec$ on
$F$. Then $T_{I,\prec}$ is precisely the subset of words in $F$
that do not have any of $\left\{ \text{HT}_{\prec}\left(q\right):q\in\mathcal{Q}\right\} $
as a prefix.
\end{thm}
To construct a CRS for $I$, we associate to each element $f$ in
the exposure basis $B_{I,\prec}$ an element $s\in I$, called its
\texttt{second}\emph{ }(with $f$ referred to as a \texttt{first,}
following Rosenmann \cite{Rosenmann1993}). Let $b\in S\cup S^{-1}$
be the last letter of $\text{HT}_{\prec}\left(f\right)$. The \texttt{second}
$s$ is defined via the exposure process of $I$, as the remainder
of $fb^{-1}$ modulo the minimal Schreier transversal of the intermediate
ideal generated by the previously discovered exposure basis elements
(see Definition \ref{def: (HLL of first, second, highest term of second)}). 

We define the \emph{Gröbner basis of $I$ with respect to $\prec$},
denoted $B_{I,\prec}^{\text{gr}}$, to consist of the exposure basis
elements $B_{I,\prec}$ together with their associated \texttt{seconds}
(see Definition \ref{def: The Gr=0000F6bner Basis}).
\begin{thm}
\label{thm: introduction Grobner basis is combinatorially reducing }Let
$\prec$ be an exposure order on $F$ and let $I\lneqq\mathcal{A}$.
The Gröbner Basis $B_{I,\prec}^{\text{gr}}$ is a combinatorially
reducing system for $I$ with respect to $\prec$.
\end{thm}
Distinct CRSs may exist for $I$. We work with $B_{I,\prec}^{\text{gr}}$
for two reasons. First, it contains the exposure basis $B_{I,\prec}$.
Second, and importantly for both theory and computation, it respects
the inductive nature of the exposure process: for any intermediate
sub-ideal $I'\subseteq I$ arising in this process, we have $B_{I',\prec}^{\text{gr}}\subseteq B_{I,\prec}^{\text{gr}}$
(see Remark \ref{rem: why we define seconds this way}).

\subsection{Algorithmic Results}

The generality of exposure orders comes at an algorithmic cost: as
noted earlier, some standard reductions from Gröbner theory -- Rosenmann's
``div-reductions'' and Madlener and Reinert's ``prefix-reductions''
-- may strictly $\prec_{\text{max}}$-increase the support. Nonetheless,
when the exposure order $\prec$ and field operations of $K$ are
computable (i.e., algorithms exist for performing field operations
and for $\prec$-comparing elements of $F$), we are able to provide
explicit algorithms for working with an exposure order in $\mathcal{A}=K\left[F\right]$.
\begin{thm}
\label{thm: algorithms}Let $\prec$ be a computable exposure order
on a free group $F$ and $K$ a field with computable operations.
Let $I\leq\mathcal{A}$ be a finitely generated right ideal. There
exist explicit algorithms for:
\begin{itemize}
\item Computing the remainder $\phi_{I,\prec}\left(f\right)$ modulo the
minimal Schreier transversal $T_{I,\prec}$ of $I$, given $f\in\mathcal{A}$
and a combinatorially reducing system $\mathcal{Q}$ for $I$ --
Algorithm \ref{alg:reduction_mod_TI}.
\item More generally, computing, given such $f$ and $\mathcal{Q}$, the
canonical division with remainder of Theorem \ref{thm: division with remainder}
of $f$ by $\mathcal{Q}$ with remainder in $T_{I,\prec}$ -- Algorithm
\ref{alg:division_with_remainder}.
\item Computing the exposure basis $B_{I,\prec}$ and Gröbner basis $B_{I,\prec}^{\text{gr}}$
given a finite generating set for $I$ -- Algorithm \ref{alg:exposure_grobner}
and its simplified version, Algorithm \ref{alg:single_generator_basis}
, for ideals generated by a single element.
\item Expressing the \texttt{seconds} canonically as an $\mathcal{A}$-linear
combination of the exposure basis $B_{I,\prec}$ -- Algorithm \ref{alg:compute_matrix_C}.
\item More generally, expressing an element $f\in\mathcal{A}$ canonically
as an $\mathcal{A}$-linear combination of $B_{I,\prec}$ -- Algorithm
\ref{alg:express_using_exposure}.
\end{itemize}
\end{thm}
The coefficients computed by Algorithm \textit{\ref{alg:express_using_exposure}}
carry both theoretical and algorithmic significance, as they encode
structural properties of the expressed element $h$ within the ideal
$I$. A motivating example of this is given in \cite{Seidel2025},
where these coefficients are used to test whether a given $h\in I$
is \emph{primitive}, that is, whether it belongs to some basis of
$I$. In that work, the expression of the corresponding \texttt{seconds}
in terms of the Rosenmann basis is described only informally, as something
that could be computed during Rosenmann's algorithm. Here, we close
that gap by providing Algorithm \textit{\ref{alg:compute_matrix_C}}.

To develop the algorithms mentioned in Theorem \ref{thm: algorithms},
we address the challenge that reductions may $\prec_{\text{max}}$-increase
the support by employing several different techniques, all of which
rely on the directly-defined minimal Schreier transversal. 

First, termination in Algorithm \ref{alg:reduction_mod_TI}, which
computes $\phi_{I,\prec}$, follows from the structure of the minimal
Schreier transversal $T_{I,\prec}$. Each word $u$ outside of $T_{I,\prec}$
admits a unique shortest prefix $u_{\text{exit}}$ lying outside $T_{I,\prec}$,
and a corresponding suffix $u_{\text{suf}}$ such that $u=u_{\text{exit}}u_{\text{suf}}$.
We guarantee termination by showing that words outside of $T_{I,\prec}$
which replace $u$ in a single reduction step are either closer to
$T_{I,\prec}$ (i.e., $v_{\text{suf}}$ is shorter than $u_{\text{suf}}$)
or have $\prec$-smaller $v_{\text{exit}}$ (see Theorem \ref{thm: Algorithm: reduction modulo T_I given combinatorially reducing system for I}
for precise details). The crux here is that while $\phi_{I,\prec}$$\left(f\right)$
may have $\prec_{\text{max}}$-greater support than $f$ does for
general $f\in\mathcal{A}$, this does not happen when $f$ is supported
on $T_{I,\prec}$ and its neighbors (see Proposition \ref{prop: Reducing a neighbour of minimal Schreier Transversal T_I decreases support}).

Second, for computing exposure bases and Gröbner bases, we introduce
a combinatorial criterion for determining whether a given element
$f\in\mathcal{A}$ is a valid ``next'' exposure basis element. An
element satisfying the combinatorial conditions of this criterion
is called \emph{exposure-extending} for $I$ (see Definition \ref{def: exposure extending}).
The following theorem establishes that adding such an element to $I$
indeed extends the exposure basis in the expected way.
\begin{thm}
\textbf{\label{thm: Extending an exposure basis}(Extending an Exposure
Process)} Let $\prec$ be an exposure order and let $I\lneqq\mathcal{A}$.
Suppose that $f$ is an exposure-extending element for $I$, and let
$I'$ be the ideal generated by $I\cup\left\{ f\right\} $. Then
\[
B_{I',\prec}=B_{I,\prec}\cup\left\{ f\right\} ,\text{ and }B_{I',\prec}^{\text{gr}}=B_{I,\prec}^{\text{gr}}\cup\left\{ f,s\right\} ,
\]
where $s$ is the $I$-\texttt{second} of $f$ from Definition \ref{def: I-second}.
\end{thm}
For a shortlex order $\prec$, Rosenmann's approach to computing minimal
bases consisted of three distinct types of reductions -- orbit reductions,
div reductions, and LCM reductions -- each strictly $\prec_{\text{max}}$-decreasing
the support of an element. For each generator, the algorithm keeps
a list of possible elements to reduce it by, and terminates when no
further reductions are needed.

Algorithm \ref{alg:exposure_grobner} takes a somewhat simpler conceptual
approach, relying instead on the criterion of Theorem \ref{thm: Extending an exposure basis}
to decide termination. Rather than managing the interplay between
three types of reductions on both \texttt{firsts} and \texttt{seconds},
at each step, a generator $f$ is tested for being exposure-extending
with respect to a ``current'' sub-ideal $J$. If $f$ passes the
test, it is added to $B_{J}$, and its $J$-\texttt{second} is added
to $B_{J}^{\text{gr}}$. Otherwise, the test produces an element (either
$\phi_{J}\left(f\right)$ or the $J$-\texttt{second} of $f$) which
replaces $f$ but has smaller $\prec_{\text{max}}$-support in an
extended sense, where elements not supported on the ``current''
minimal Schreier transversal $T_{J,\prec}$ are considered as having
infinite support. 

This extended dynamic ordering of supports, with respect to the ``current''
sub-ideal $J$, is the heart of our argument for ensuring termination.
By considering $\infty$ as the ``support'' for elements outside
of normal form with respect to $J$, every reduction (in particular,
applying $\phi_{J}$) strictly decreases support, while the extended
set of values, obtained by adding $\infty$, remains well-ordered.

\subsection{Notation and Scope}

Throughout the paper, $K$ is a fixed field and $F$ the free group
with finite basis $S$. We denote the identity of $F$ by $e$, and
write $\mathcal{A}:=K\left[F\right]$ for the associated free group
algebra. Whenever we write a sum in $\mathcal{A}$ (for example $\sum_{u\in F}\alpha_{u}u$
with $\alpha_{u}\in K$), it is understood that only finitely many
terms are nonzero. For a subset $T\subseteq F$, we denote by $\text{Sp}_{K}\left(T\right)$
the $K$-linear span of $T$. The \emph{support} of an element $f\in\mathcal{A}$,
denoted by $\text{supp}\left(f\right)$, is the set of words in $F$
that have nonzero coefficient in $f$. We say that $f$ \emph{is supported
on} a subset $T\subseteq F$ if $f\in\text{Sp}_{K}\left(T\right)$.
If $I$ is a right ideal of $\mathcal{A}$, we write $I\leq\mathcal{A}$.
All results are formulated for right ideals of $\mathcal{A}$, and
unless explicitly stated otherwise, the term \emph{ideal} will always
mean a right ideal of $\mathcal{A}$. 

The framework developed here extends, with minor modifications, to
several more general settings:
\begin{enumerate}
\item \textbf{Left ideals.} All results transfer to left ideals of $\mathcal{A}$
by applying the inversion map $\iota$ (see \cite[Definition~4.1]{Seidel2025}),
which maps right ideals to left ideals and vice versa, as an automorphism
between $\mathcal{A}$ and its opposite ring $\mathcal{A}^{\text{op}}$. 
\item \textbf{Exposure orders based at arbitrary roots.} Instead of the
trivial word $e$, one may base exposure orders at an arbitrary element
$\rho\in F$, requiring $\rho u\preceq\rho v$ whenever $u$ is a
prefix of $v$. The right $\mathcal{A}$-module automorphism $L_{\rho}:\mathcal{A}\rightarrow\mathcal{A}$,
defined by $L_{\rho}\left(f\right)=\rho f$, maps $e$ to $\rho$.
Conjugation by $L_{\rho}$ transfers all results to this generalized
setting.
\item \textbf{Free groups of infinite rank.} The arguments remain valid
for free groups of arbitrary rank. The only modification required
is an ordinal of larger cardinality to guarantee termination of the
exposure process in Proposition \ref{prop:The-process-stabilizes}.
\item \textbf{Submodules of free $\mathcal{A}$-modules.} Since $\mathcal{A}$
is a free ideal ring, every submodule of a free $\mathcal{A}$-module
is itself free, with a well-defined rank. For a free right $\mathcal{A}$-module
with basis $\mathcal{E}$ of arbitrary (possibly infinite) rank, exposure
orders are adapted by requiring $\varepsilon u\preceq\varepsilon v$
for every $\varepsilon\in\mathcal{E}$ and words $u,v\in F$ such
that $u$ is a prefix of $v$. 
\end{enumerate}

\subsection{Paper Organization}

The paper is structured as follows. Section \ref{sec: Preliminary Concepts}
provides preliminary material, defining and illustrating the different
orders on $F$ and its subsets, and recalling relevant aspects of
well-orders and ordinals. Section \ref{sec: Induced Minimal Structures}
forms the core theoretical part of the paper, defining the minimal
structures induced by an exposure order on an ideal and establishing
their properties. Section \ref{sec: Algorithmic Framework} complements
the previous section by providing algorithmic criteria and explicit
algorithms for computing and using the minimal structures induced
by exposure orders in $\mathcal{A}$.

We describe the content of each section in more detail. Section \ref{sec Exposure orders on F}
introduces exposure orders on $F$, together with motivating examples
and machinery for constructing such orders. Section \ref{sec: Well Orders and Ordinals}
recalls relevant notions from the theory of well-orders and ordinal
numbers, and motivates their use. Section \ref{sec: Ordering the Subsets of F}
introduces and establishes basic properties of two ways to compare
subsets of $F$: the usual $\prec_{\text{\ensuremath{\max}}}$ for
comparing supports, and the auxiliary $\prec_{\text{min}}$ for comparing
Schreier transversals. Section \ref{sec: The Minimal Schreier Transversal}
proves Theorem \ref{thm: minimal Schreier transversal exists}, establishing
the existence of the $\prec_{\text{min}}$-minimal Schreier transversal
$T_{I,\prec}$, along with properties of $T_{I,\prec}$ and its associated
transversal function $\phi_{I,\prec}$. Section \ref{sec: Highest Terms and Minimal Monic Elements}
defines head terms and monic elements, and describes how elements
with $\prec_{\text{max}}$-minimal support are selected. Section \ref{sec: The Exposure Basis}
introduces the exposure process, in which the exposure basis $B_{I,\prec}$
is defined via transfinite recursion. Termination of the inductive
process is proven, and a useful form for the exposure basis elements
is provided (Proposition \ref{prop: Form for first using T_I}). Several
examples are given in Example \ref{exa: exposure process}. Section
\ref{sec: The Seconds} defines the \texttt{seconds} corresponding
to exposure basis elements and illustrates this construction with
several examples (see Example \ref{examples of seconds}). Section
\ref{sec: Combinatorially Reducing Systems} defines combinatorially
reducing systems for an ideal, showing that they allow canonical representations
(Lemma \ref{lem: Canonical-Representation-using-1}) and division
with remainder (Theorem \ref{thm: division with remainder}). In the
same section, Theorem \ref{thm: combinatorially reducing system defines minimal schreier transversal of the ideal which it generates}
is proven, showing that such a system determines $T_{I,\prec}$ by
its prefix structure, and Algorithm \ref{alg:reduction_mod_TI} for
computing $\phi_{I,\prec}$ is provided, along a correctness proof.
Section \ref{sec: Combinatorial Properties of the grobner basis}
proves Theorem \ref{thm: introduction Grobner basis is combinatorially reducing },
and discusses alternative combinatorially reducing systems for $I$.
Section \ref{sec: The Exposure Basis is a Basis} proves Theorem \ref{thm: The exposure basis is a basis},
showing that the exposure basis $B_{I,\prec}$ is indeed a basis for
$I$. Additionally, every Lewin basis is realized as an exposure basis
(Proposition \ref{prop: Lewin bases are exposure bases}). Section
\ref{sec: Suffix Preserving Orders} defines suffix-invariant orders
on $F$, proves that every such order is an exposure order, provides
examples and shows that these orders yield a simplified theory (Theorem
\ref{claim: Suffix-invariant orders are exposure orders}). Section
\ref{sec: Extending an Existing Exposure Basis} proves the criterion
of Theorem \ref{thm: Extending an exposure basis} for extending an
exposure basis. Finally, Section \ref{sec: Algorithm for Computing Bases}
provides algorithmic tools for exposure orders, including Algorithm
\ref{alg:exposure_grobner} for computing $B_{I,\prec}$ and $B_{I,\prec}^{\text{gr}}$,
Algorithm \ref{alg:division_with_remainder} for division with remainder,
and Algorithm \textit{\ref{alg:express_using_exposure}}\textit{\emph{
for expressing an element $f\in I$ using the exposure basis $B_{I}$.}}

\subsection*{Acknowledgements}

The author would like to thank Doron Puder for his invaluable support
and insightful advice in improving the clarity and readability of
the manuscript. The author is also grateful to Shachar Herpe for her
help in overcoming his fear of ordinal numbers. This work was supported
by the European Research Council (ERC) under the European Union\textquoteright s
Horizon 2020 research and innovation programme (grant agreement No
850956).

\section{\label{sec: Preliminary Concepts}Preliminary Concepts}

\subsection{\label{sec Exposure orders on F}Exposure orders on $F$}

We work with the right Cayley tree $\Gamma$ of $F$ with respect
to its fixed basis $S$. 3 of $\Gamma$ is $F$, and for each $s\in S$
and $v\in V$, there is a directed edge from $v$ to $vs$ labeled
$s$. When traversed in the opposite orientation -- from $vs$ to
$v$ -- we think of this edge as labeled $s^{-1}$.

A \emph{s}ubset $T\subseteq F$ is called \emph{prefix-closed} if
for every $w\in T$ and every prefix $u$ of $w$, we have $u\in T.$
Since the unique path in $\Gamma$ from any $w\in F$ to $e$ is made
of all prefixes of $w$, a nonempty subset $T\subseteq F$ is prefix-closed
if and only if it is a subtree of $\Gamma$ containing $e$. From
this perspective, we define a\emph{ prefix-neighbor }of a prefix-closed
subset $T\subseteq F$ to be a word outside of $T$ whose proper prefixes
all lie in $T$. We denote the set of prefix-neighbors of $T$ by
$\partial T$. Note that the notion of prefix-neighbors extends, but
is not identical to, the usual notion of neighbors from graph theory:
in particular, the empty set is prefix-closed and its only prefix-neighbor
is $e$.

Recall that a totally ordered set $\left(X,\prec\right)$ is called
\emph{well-ordered} if every nonempty subset of $X$ admits a minimum,
or, equivalently, if there exists no infinite strictly decreasing
sequence of elements of $X$. For a review of well-orders and our
motivation for working with them, see Section \ref{sec: Well Orders and Ordinals}.

Recall from Definition \ref{def: exposure orders} that an \emph{exposure
order }is a well-ordering of $F$ such that $u\preceq v$ whenever
$u$ is a prefix of $v$. Intuitively, one may think of such an order
as arising from transfinitely recursive process in which words are
``exposed'' one after the other, beginning with $e$ and proceeding
by successively adjoining prefix-neighbors.
\begin{claim}
Let $\prec$ be an exposure order on $F$. Then $\min F=e$ and for
every $u\in F\backslash\left\{ e\right\} $ the set $T_{u}:=\left\{ v\in F:v\prec u\right\} $
is a subtree of $\Gamma$ of which $u$ is a neighbor.
\end{claim}
\begin{proof}
The identity of $F$ is a prefix of any other word, so $\min F=e$.
Let $u\in F\backslash\left\{ e\right\} $. Since the order is an exposure
order, the set $T_{u}$ is prefix-closed and contains $e$, hence
a subtree of $\Gamma$. All proper prefixes of $u$ lie in $T_{u}$,
so $u$ is a neighbor of $T_{u}$.
\end{proof}
Given an exposure order on $F$, the following claim allows us to
construct more elaborate ones for use in subsequent examples. Recall
from the theory of well-ordered sets that if $X_{1}$ and $X_{2}$
are well-ordered sets, their \emph{sum} $X_{1}+X_{2}$ is defined
as their disjoint union, ordered such that every element of $X_{1}$
is smaller than every element of $X_{2}$, while preserving the original
ordering within each $X_{i}$. The resulting set $X_{1}+X_{2}$ is
again well-ordered.
\begin{claim}
\label{claim: Sum of two order on tree and complement is an exposure order too}Let
$\prec$ be an exposure order on $F$, and let $T\subseteq F$ be
prefix-closed. Then the order $\prec_{T}$ on $F$, obtained by viewing
$F$ as the sum of $\left(T,\prec\mid_{T}\right)$ and its complement
$\left(F\backslash T,\prec\mid_{F\backslash T}\right)$, is also an
exposure order.
\end{claim}
\begin{proof}
Let $w\in F$ and let $v$ be a proper prefix of $w$. Since $\prec$
is an exposure order, we have $v\prec w$. We must show that $v\prec_{T}w$
too. If $w,v\in T$ or $w,v\in F\backslash T$ then $v\prec_{T}w$.
Otherwise, $v\in T$ and $w\in F\backslash T$ and then, too, $v\prec_{T}w$.
\end{proof}
\begin{example}
\label{exa: Examples of Exposure orders}Let $F=\left\langle x,y\right\rangle $
be the free group on generators $x$ and $y$. The following examples
illustrate Definition \ref{def: exposure orders}.
\begin{enumerate}
\item \label{enu: shortlex}The \emph{shortlex} order (also known as \emph{length-lexicographic}
order) used by Rosenmann \cite{Rosenmann1993}, orders words in $F$
first by length and then lexicographically among words of the same
length, according to an arbitrary fixed order on $S\cup S^{-1}$.
For example, choosing $y^{-1}\prec x^{-1}\prec x\prec y$, the smallest
elements are:
\[
e\prec y^{-1}\prec x^{-1}\prec x\prec y\prec y^{-2}\prec y^{-1}x^{-1}\prec y^{-1}x\prec x^{-1}y^{-1}\prec x^{-2}\prec...
\]
Each word in $F$ has only finitely many smaller words, so the shortlex
order is a well-order: every strictly decreasing sequence must be
finite. It is also clearly an exposure order, since proper prefixes
are shorter, and thus smaller.
\item \label{enu: example needing transfinite induction}Let $T:=\left\{ x^{n}:n\geq0\right\} $
denote all nonnegative powers of $x$. Consider the exposure order
$\prec_{T}$, defined as in Claim \ref{claim: Sum of two order on tree and complement is an exposure order too}
using the shortlex order $\prec$ and the prefix-closed set $T$.
In this order, the smallest elements are: 
\[
e\prec_{T}x\prec_{T}x^{2}\prec_{T}x^{3}\prec_{T}...\prec_{T}y^{-1}\prec_{T}x^{-1}\prec_{T}y\prec_{T}y^{-2}\prec_{T}...
\]
Note that in this order, the element $y^{-1}$ has infinitely many
smaller words.
\item Define the order $\triangleleft$ on $F$ as follows. Let 
\[
T=\left\{ x^{n}:n\geq0\right\} ,\ Ty=\left\{ x^{n}y:n\geq0\right\} .
\]
Order $T$ using the shortlex order above. Order $Ty$ in reverse
order: $x^{n}y\triangleleft x^{m}y$ whenever $n>m$. Order the complement
$F\backslash\left(T\cup Ty\right)$ using shortlex, again. Now define
$\triangleleft$ by placing all elements of $Ty$ in the middle --
larger than all elements of $T$ and smaller than all remaining elements
of $F\backslash\left(T\cup Ty\right)$. The smallest elements are:
\[
\underbrace{e\triangleleft x\triangleleft x^{2}\triangleleft x^{3}\triangleleft...}_{T}\triangleleft\underbrace{...\triangleleft x^{3}y\triangleleft x^{2}y\triangleleft xy\triangleleft y}_{Ty}\triangleleft y^{-1}\triangleleft x^{-1}\triangleleft...
\]
This order satisfies the condition that whenever $u$ is a proper
prefix of $v$, then $u\triangleleft v$. However, it is not an exposure
order since it is not a well-order: it admits an infinite strictly
decreasing sequence $y\triangleright xy\triangleright x^{2}y\triangleright x^{3}y\triangleright...$
\end{enumerate}
\end{example}

\subsection{\label{sec: Well Orders and Ordinals}Well-Orders and Ordinals}

This section recalls some basic facts about well-orders and ordinals,
and motivates the use of such machinery for our cause. It may well
be skipped by readers familiar with these concepts. For a comprehensive
treatment of the facts recalled here, see \cite[Chapter~1]{Kunen2014}.

A first motivation for using well-orders, from an algorithmic point
of view, lies in their characterization via strictly decreasing sequences.
This characterization provides a method for ensuring that algorithms
terminate: some parameter taking values in a well-ordered set is shown
to strictly decrease in each iteration, and thus the algorithm must
terminate after finitely many such steps.

A second motivation for using well-orders is that they enable \emph{transfinite
induction} as a proof technique. Given a property $P_{x}$ defined
for each $x\in X$, where $X$ is well-ordered, one proves $P_{x}$
assuming that $P_{y}$ holds for all $y\prec x$; transfinite induction
then guarantees that $P_{x}$ holds for all $x\in X$. The order $\prec_{T}$
from Example \ref{exa: Examples of Exposure orders} (\ref{enu: example needing transfinite induction})
illustrates why this is sometimes necessary: a standard (finite) induction
would prove $P_{e}$, then $P_{x}$,$P_{x^{2}}$ and so on -- covering
just the subset $\left\{ x^{n}:n\geq0\right\} \subsetneqq F$.

A third motivation for using well-orders is to allow definitions by
\emph{transfinite recursion}. From the theory of well-ordered sets,
recall that \emph{ordinal numbers} (or \emph{ordinals}) represent
the order types of well-ordered sets. Every pair of ordinals is comparable
by inclusion. The set of countable ordinals, denoted by $\omega_{1}$,
is well-ordered with respect to this order, and moreover constitutes
the first uncountable ordinal. Each ordinal $\alpha$ is associated
with its \emph{successor ordinal} $\alpha+1$, which can be defined
as the minimal ordinal larger than $\alpha$.

\emph{In transfinite recursion}, one defines a family of objects $f_{\alpha}$,
indexed by ordinals $\alpha$, by specifying how to define $f_{\alpha}$
given all $f_{\beta}$ for $\beta<\alpha$. This guarantees that $f_{\alpha}$
is defined for all ordinals. A stopping condition can be imposed:
for example, to halt at the first ordinal $\alpha_{0}$ where a condition
on $\left\{ f_{\beta}\right\} _{\beta<\alpha_{0}}$ is met. In this
case, the process defines $f_{\beta}$ only for $\beta<\alpha_{0}$.

We aim to define a basis for a given ideal using transfinite recursion
(see Definition \ref{def: The Exposure Basis}). The following example
illustrates how transfinite recursion naturally arises in such settings.
\begin{example}
\label{exa: vector space example}Let $V$ be a vector space over
a field $K$. A basis for $V$ can be constructed by transfinite recursion
as follows. At each ordinal $\alpha$, suppose $\left\{ v_{\beta}:\beta<\alpha\right\} $
has been chosen to be linearly independent. If it spans $V$, halt.
Otherwise, choose $v_{\alpha}\in V\backslash\text{Sp}_{K}\left\{ v_{\beta}:\beta<\alpha\right\} $.
Suppose $V$ has a countable basis $\left(e_{i}\right)_{i=1}^{\infty}$.
We claim that the process must terminate before $\omega_{1}$. Otherwise,
we obtain an uncountable linearly independent set in a countable-dimensional
space. But since $V=\bigcup_{n=1}^{\infty}\text{Sp}_{K}\left\{ e_{1},e_{2},...,e_{n}\right\} $,
some finite-dimensional subspace must contain uncountably many of
the $v_{\beta}$, contradicting their independence. In contrast, constructing
a basis by ordinary (countable) induction may fail to span $V$ --
for example, if the sequence $v_{i}=e_{2i}$ is chosen.
\end{example}

\subsection{\label{sec: Ordering the Subsets of F}Ordering the Subsets of $F$}

In this section, we show that an exposure order $\prec$ on $F$ gives
rise to two distinct orders for comparing subsets of $F$. The first,
denoted $\prec_{\text{max}}$, is defined only on finite subsets of
$F$ and will be used to compare the supports of elements of $\mathcal{A}$.
It is a well-order, which lets us pick an element having minimal support
and ensures that algorithms terminate by tracking $\prec_{\text{max}}$-decreasing
supports. The second, denoted $\prec_{\text{min}}$, applies more
generally to all subsets of $F$, though it is not a well-order. It
will be used to compare linearly independent subsets modulo a one-sided
ideal $I$, and will allow us to associate to $I$ a \emph{minimal
Schreier transversal} (see Theorem \ref{thm: minimal Schreier transversal exists}).

We work in the generality of a totally-ordered set $X$, and prove
the properties mentioned above when $X$ is well-ordered. We use $\triangle$
to denote symmetric difference.
\begin{defn}
\textbf{\label{def: order on finite subsets}}Let $A$ and $B$ be
distinct subsets of a totally-ordered set $X$.
\begin{enumerate}
\item Define $A\prec_{\max}B$ if $\max\left(A\triangle B\right)\in B$.
\item Define $A\prec_{\min}B$ if $\min\left(A\triangle B\right)\in A$.
\end{enumerate}
\end{defn}
Since $A\neq B$, the symmetric difference $A\triangle B$ is nonempty.
However, it may lack a minimum or maximum, in which case $A$ and
$B$ are $\prec_{\text{min}}$-incomparable or $\prec_{\text{max}}$-incomparable,
respectively.

As an example, if finite subsets $A,B\subseteq X$ satisfy $A\subseteq B$,
then $A\preceq_{\text{max}}B$. In particular, the empty set is $\preceq_{\text{max}}$-smaller
than any nonempty finite subset of $X$. This aligns with our intended
use of $\preceq_{\text{max}}$ for comparing supports of elements
in $\mathcal{A}$: reductions that eliminate words from the support
should strictly decrease it. In contrast, if $A\subseteq B$ and the
sets are $\preceq_{\text{min}}$-comparable, then $A\succeq_{\text{min}}B$.
This fits the intended use of $\preceq_{\text{min}}$ for comparing
$K$-linearly independent subsets, where subsets containing extra
words are more desirable. Here too, the $\preceq_{\text{min}}$-smaller
subset corresponds to the preferred one.

We will use the following observation relating the minimum or maximum
of sets to that of their symmetric difference:
\begin{claim}
\label{claim: minimum of symmetric difference}Let $A$ and $B$ be
subsets of a totally ordered set.
\begin{enumerate}
\item If $A$ and $B$ have distinct minima $a$ and $b$, then $\min\left(A\triangle B\right)=\min\left\{ a,b\right\} $.
\item If $A$ and $B$ have distinct maxima $a$ and $b$, then $\max\left(A\triangle B\right)=\max\left\{ a,b\right\} $.
\end{enumerate}
We omit the straightforward proof.
\end{claim}
\begin{cor}
\label{cor: The two orders are transitive}The relations $\prec_{\text{min}}$
and $\prec_{\text{max}}$ of Definition \ref{def: order on finite subsets}
are partial orders on $X$.
\end{cor}
\begin{proof}
We prove for $\prec_{\text{min}}$, the proof for $\prec_{\text{max}}$
is similar. The relation $\prec_{\text{min}}$ is anti-reflexive and
anti-symmetric by definition. For transitivity, let $A,B,C\subseteq X$
satisfy $A\prec_{\text{min}}B\prec_{\text{min}}C$. Let $a=\min\left(A\triangle B\right)$
and $b=\min\left(B\triangle C\right)\in B$. Then $a\in A\backslash B$
and $b\in B\backslash C$. In particular, $a$ and $b$ are distinct
since $a\notin B$ and $b\in B$. By Claim \ref{claim: minimum of symmetric difference}
we have:
\[
\min\left(A\triangle C\right)=\min\left(\left(A\triangle B\right)\triangle\left(B\triangle C\right)\right)\in\left\{ a,b\right\} .
\]
Since both $a$ and $b$ are not contained in $C\backslash A$, we
deduce that $A\prec_{\text{min}}C$.
\end{proof}
\begin{claim}
\label{claim: The order on finite subsetes is a well-order}If $X$
is well-ordered, then:
\end{claim}
\begin{enumerate}
\item $\prec_{\text{min}}$ is a total order on all subsets of $X$.
\item $\prec_{\text{max}}$ is a well-ordering on the finite subsets of
$X$.
\end{enumerate}
\begin{proof}
Let $A$ and $B$ be distinct subsets of $X$. Then $A\triangle B\neq\emptyset$.
Since $X$ is well-ordered, $A\triangle B$ has a minimum, so $A$
and $B$ are $\prec_{\text{min}}$-comparable. If $A$ and $B$ are
also finite, then $A\triangle B$ is finite and thus admits a maximum,
making them $\prec_{\text{max}}$-comparable as well. By Corollary
\ref{cor: The two orders are transitive}, $\prec_{\text{min}}$ is
a total order on all subsets of $X$, and $\prec_{\text{max}}$ is
a total order on its finite subsets. 

To prove that $\prec_{\text{max}}$ is a well-order, we first observe
that if $\left(S_{i}\right)_{i=1}^{\infty}$ is a (weakly) decreasing
sequence of nonempty finite subsets of $F$ then the corresponding
sequence of maxima $\left(\max S_{i}\right)_{i=1}^{\infty}$ must
stabilize: Indeed, it is a decreasing sequence of elements of the
well-ordered set $F$. Denote this stable maximum by $\lim_{i\rightarrow\infty}\left(\max S_{i}\right)$.
Now suppose, for contradiction, that there exists a strictly decreasing
sequence $\left(S_{i}\right)_{i=1}^{\infty}$ of finite subsets of
$F$. Such a sequence is composed of nonempty sets alone, since the
empty set is smaller than any other finite subset of $F$. Let $\left(S_{i}\right)_{i=1}^{\infty}$
be such a strictly decreasing sequence for which $w=\lim_{i\rightarrow\infty}\left(\max S_{i}\right)$
is minimal. Without  loss of generality, by discarding finitely many
elements from the sequence, we may assume that $\max S_{i}=w$ for
every $i\in\mathbb{N}$. Since the sequence is strictly decreasing
with a common maximum, the sequence $\left(S_{i}\backslash\left\{ w\right\} \right)_{i=1}^{\infty}$
is also strictly decreasing. But for every $i\in\mathbb{N}$ we have
$\max\left\{ S_{i}\backslash\left\{ w\right\} \right\} \prec w$ and
so $\lim_{i\rightarrow\infty}\max\left\{ S_{i}\backslash\left\{ w\right\} \right\} \prec w$
as well, in contradiction.
\end{proof}
\begin{rem}
\label{rem: le_min is not a well-order}In contrast to $\prec_{\text{max}}$,
the order $\prec_{\text{min}}$ is not a well-order, even when restricted
to the finite sets of $X$, unless $X$ itself is a finite set. Indeed,
given any sequence of $\left(x_{i}\right)_{i=1}^{\infty}$ of distinct
elements of $X$, we obtain a strictly $\prec_{\text{min}}$-decreasing
sequence of finite subsets:
\[
\emptyset\succ_{\text{min}}\left\{ x_{1}\right\} \succ_{\text{min}}\left\{ x_{1},x_{2}\right\} \succ_{\min}\left\{ x_{1},x_{2},x_{3}\right\} \succ_{\text{min}}...
\]
\end{rem}

\section{\label{sec: Induced Minimal Structures}Induced Minimal Structures}

\subsection{\label{sec: The Minimal Schreier Transversal}The Minimal Schreier
Transversal of an Ideal}

From this point forward, we fix an exposure order $\prec$ on $F$. 

In this section, we establish Theorem \ref{thm: minimal Schreier transversal exists},
and examine properties of the minimal Schreier transversal whose existence
it guarantees. We begin by recalling the relevant notions.
\begin{defn}
\label{def: Schreier Transversal}Let $I\leq\mathcal{A}$ be a right
ideal. A subset $T\subseteq F$ is called a \emph{partial Schreier
transversal} for $I$ if it is closed under taking prefixes and $\text{Sp}_{K}T\cap I=\left\{ 0\right\} $.
A partial Schreier transversal $T\subseteq F$ is called a \emph{Schreier
transversal} for $I$ if, in addition, every $I$-coset in $\mathcal{A}$
intersects $\text{Sp}_{K}\left(T\right)$.
\end{defn}
Hence, if $T$ is a Schreier transversal for $I$, then every $I$-coset
has a unique representative in $\text{Sp}_{K}\left(T\right).$ 
\begin{claim}
\label{claim: A Schreir transversal is inclusion-maximal among partial Schreier transversals}Let
$T$ be a Schreier transversal for $I$. Then $T$ is maximal among
the partial Schreier transversals with respect to inclusion.
\end{claim}
\begin{proof}
Let $T'$ be a partial Schreier transversal such that $T\subseteq T'$.
Suppose towards contradiction that $T\neq T'$. Then there exists
$w\in T'\backslash T$. Since $T$ is a Schreier transversal, there
exists a representative $f\in\text{Sp}_{K}\left(T\right)$ in the
coset $w+I$. Then the difference  $f-w$ lies in $\text{Sp}_{K}\left(T'\right)\cap I$.
Since $T'$ is a partial Schreier transversal, $f-w=0$, so $w\in\text{Sp}_{K}\left(T\right)$,
in contradiction. 
\end{proof}
We recall the following classical result \cite[Lemma~1]{Lewin1969}.
\begin{lem}
\label{lem: partial ST can be extended}Let $T'$ be a partial Schreier
transversal for an ideal $I\leq\mathcal{A}$. Then there exists a
Schreier transversal $T$ for $I$ such that $T'\subseteq T$.
\end{lem}
The empty set is a partial Schreier transversal for any ideal $I\leq\mathcal{A}$.
Therefore, by Lemma \ref{lem: partial ST can be extended}, $I$ admits
a Schreier transversal. Among the possibly many Schreier transversals
of $I$, Theorem \ref{thm: minimal Schreier transversal exists} will
single out a canonical one: the \emph{minimal Schreier transversal}
$T_{I,\prec}$ (minimal with respect to $\prec_{\text{min}}$).
\begin{example}
\label{exa: The minimal Schreier transversal depends on the order}The
minimal Schreier transversal depends on the fixed exposure order $\prec$
on $F$. For example, let $F=\left\langle x\right\rangle $ be the
free group on a single generator $x$ and let $I$ be the right ideal
generated by $x^{2}-1$. In the shortlex order $\prec_{1}$ satisfying
$x\prec_{1}x^{-1}$ (see Example \ref{exa: Examples of Exposure orders}),
we have $T_{I,\prec_{1}}=\left\{ e,x\right\} $, whereas in the shortlex
order $\prec_{2}$ satisfying $x^{-1}\prec_{2}x$ we have $T_{I,\prec_{2}}=\left\{ e,x^{-1}\right\} $. 
\end{example}
We are now ready to prove Theorem \ref{thm: minimal Schreier transversal exists}.
Since we work with a fixed exposure order, we omit the dependence
on the order from the notation and simply write $T_{I}$.

\begin{proof}[Proof of Theorem \ref{thm: minimal Schreier transversal exists}]

Let $X$ be the set of Schreier transversals for $I$. Define $T_{I}$
to be the set of words of $F$ which lie in every small enough transversal.
Explicitly, for $u\in F$ we say that $u\in T_{I}$ if there exists
some $T_{u}\in X$ such that for every $T\in X$ with $T\preceq_{\text{min}}T_{u}$
we have $u\in T$. 

We first show that $T_{I}$ is a lower bound on $X$. Suppose otherwise,
and let $T\in X$ be a Schreier transversal for $I$ having, among
those smaller than $T_{I}$, a minimal $\min\left(T\triangle T_{I}\right)$.
Let $w=\min\left(T\triangle T_{I}\right)$. Since $T\prec_{\text{min}}T_{I}$
we have $w\in T\backslash T_{I}$, so by the definition of $T_{I}$
there exists some $T'\in X$ such that $T'\prec_{\text{min}}T$ and
$w\notin T'$. Let $u=\min\left(T\triangle T'\right)$. Since $w$
lies in this symmetric difference as well and $u$ is its minimum,
we have $u\preceq w$. But $u$ and $w$ are distinct since $u\notin T$.
By Claim \ref{claim: minimum of symmetric difference},
\[
\min\left(T'\triangle T_{I}\right)=\min\left(\left(T'\triangle T\right)\triangle\left(T\triangle T_{I}\right)\right)=\min\left\{ u,w\right\} =u.
\]
Since $u\in T'$, we have $T'\prec_{\text{min}}T_{I}$, in contradiction
to the minimality defining $T$.

We next show that $T_{I}$ is a partial Schreier transversal for $I$.
Let $u\in T_{I}$ and let $v$ be a prefix of $u$. Define $T_{v}:=T_{u}$
and let $T\in X$ satisfy $T\preceq_{\text{min}}T_{v}$. Since $T\preceq_{\text{min}}T_{u}$,
we have $u\in T$. As $T$ is closed under prefixes, $v\in T$ as
well. Thus, $v\in T_{I}$ so $T_{I}$ is prefix-closed. Let $f\in\text{Sp}{}_{K}\left(T_{I}\right)\cap I$.
Suppose toward contradiction that $f\neq0$. Let $T=\min\left\{ T_{u}:u\in\text{supp}f\right\} $.
Then $f\in\text{Sp}_{K}\left(T\right)\cap I$, in contradiction to
$T$ being a Schreier transversal for $I$.

Since $T_{I}$ is a partial Schreier transversal for $I$, by Lemma
\ref{lem: partial ST can be extended} there exists a Schreier transversal
$T_{I}'$ for $I$ such that $T_{I}\subseteq T_{I}'$. By the inclusion
relation we have $T_{I}'\preceq_{\text{min}}T_{I}$. On the other
hand, since $T_{I}'\in X$ and $T_{I}$ bounds $X$ from below, we
have $T_{I}\preceq_{\text{min}}T_{I}'$. Thus, $T_{I}=T_{I}'$ and
$T_{I}\in X$ as well.\end{proof}
\begin{defn}
\label{def: minimal Schreier transversal and transversal function}Let
$I\leq\mathcal{A}$. Let $T_{I}=T_{I,\prec}$ denote the \emph{minimal
Schreier transversal} of $I$, guaranteed by Theorem \ref{thm: minimal Schreier transversal exists}.
Let $\phi_{I}=\phi_{I,\prec}:\mathcal{A}\rightarrow\text{Sp}_{K}\left(T_{I}\right)$
denote the associated transversal function. Explicitly, $\phi_{I}\left(f\right)$
is the unique representative of the $I$-coset $f+I$ supported on
$\text{Sp}_{K}\left(T_{I}\right)$, and is called the \emph{remainder
of $f$ modulo $T_{I}$}.
\end{defn}
The transversal function $\phi_{I}$ is $K$-linear: for every $f_{1},f_{2}\in\mathcal{A}$
and $\alpha_{1},\alpha_{2}\in K$, the element $\alpha_{1}\phi_{I}\left(f\right)+\alpha_{2}\phi_{I}\left(g\right)$
lies in $\text{Sp}_{K}\left(T_{I}\right)$ and satisfies 
\[
\alpha_{1}\phi_{I}\left(f_{1}\right)+\alpha_{2}\phi_{I}\left(f_{2}\right)\in\alpha_{1}\left(f_{1}+I\right)+\alpha_{2}\left(f_{2}+I\right)=\alpha_{1}f_{1}+\alpha_{2}f_{2}+I,
\]
hence $\phi_{I}\left(\alpha_{1}f_{1}+\alpha_{2}f_{2}\right)=\alpha_{1}\phi_{I}\left(f_{1}\right)+\alpha_{2}\phi_{I}\left(f_{2}\right)$.

\begin{figure}
\centering{}\label{fig: vertex larger than its reminder and proof}\includegraphics[scale=0.5]{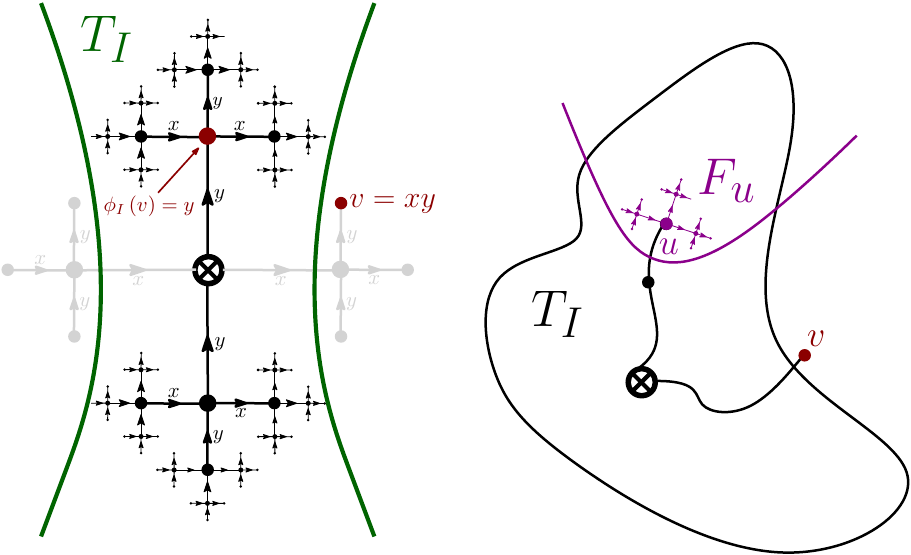}\caption{On the left, an illustration of Example \ref{exa: word smaller than its phi}:
the vertex $v=xy$ has $\prec_{\text{max}}$-smaller support than
its remainder $\phi_{I}\left(v\right)=y$ modulo $T_{I}$, with both
vertices depicted in red. On the right, the proof of Proposition \ref{prop: Reducing a neighbour of minimal Schreier Transversal T_I decreases support}
shows that this is impossible when $v$ is a prefix-neighbor of $T_{I}$:
Suppose for contradiction that the head term $u$ in its remainder
satisfies $u\succ v$. Then truncating $F_{u}$ (in purple) from $T_{I}$
and then adding $v$ (in red) results in a smaller partial Schreier
transversal, contradicting $\prec_{\text{min}}$-minimality.}
\end{figure}

We next observe that $T_{I}$ is in fact minimal among all partial
Schreier transversals for $I$. The fact that this set admits a minimum
may seem somewhat surprising when considering that $\prec_{\text{min}}$
is not necessarily a well-order (recall Remark \ref{rem: le_min is not a well-order}).
\begin{cor}
\label{cor: T_I  is minimal among partial Schreier transversals for I}Let
$T$ be a partial Schreier transversal for $I$. Then $T_{I}\preceq_{\text{min}}T$.
\end{cor}
\begin{proof}
By Lemma \ref{lem: partial ST can be extended}, there exists a Schreier
transversal $T'$ for $I$ such that $T\subseteq T'.$ By this inclusion
relation $T'\preceq_{\text{min}}T$ and by the definition of $T_{I}$
we have $T_{I}\preceq_{\text{min}}T'$, so by transitivity $T_{I}\preceq_{\text{min}}T$. 
\end{proof}
We next characterize $T_{I}$ in terms of its prefix-neighbors.
\begin{claim}
\label{claim: minimum of symmetric difference is neighbour of larger tree}Let
$T$ and $T'$ be prefix-closed sets such that $T\prec_{\text{min}}T'$.
Then $\min\left(T\triangle T'\right)$ is a prefix-neighbor of $T'$.
\end{claim}
\begin{proof}
Let $w=\min\left(T\triangle T'\right)$. Since $T\prec_{\text{min}}T'$
we have $w\in T\backslash T'$. Let $v$ be some proper prefix of
$w$. Then $v\in T$ because $T$ is prefix-closed and $v\prec w$.
By the minimality defining $w$, we have $v\notin$ $T\triangle T'$
so $v\in T'$. 
\end{proof}
The following example illustrates that although $T_{I}$ is minimal
in a certain sense, for general exposure orders, reducing an element
to its remainder modulo $T_{I}$ does not necessarily result in a
$\prec_{\text{max}}$-decrease of the support. This subtlety slightly
complicates the use of such reductions in algorithms, as a more delicate
reasoning is required to ensure termination. The example is illustrated
in the left part of Figure \ref{fig: vertex larger than its reminder and proof}. 
\begin{example}
\label{exa: word smaller than its phi}Let $I$ be the right ideal
generated by $x-1$. There is only one Schreier transversal for $I$:
the set of all words in $F$ that do not begin with $x$ or $x^{-1}$.
We omit a detailed proof, but note that it follows from Theorem \ref{Thm: properties of second + Grobner basis is combinatorially reducing}.
Specifically, since $x-1\in I$ and $x^{-1}-1\in I$, any Schreier
transversal for $I$ must exclude both $x$ and $x^{-1}$, while including
the identity element. The remaining subtree, which avoids these two
prefixes, is $K$-linearly independent modulo $I$. This unique transversal
is thus the minimal transversal $T_{I}$. Let $v=xy$. Then $\phi_{I}\left(v\right)=y$
since $y\in T_{I}$ and the difference $xy-y=\left(x-1\right)y$ lies
in $I$. Choosing some exposure order on $F$ for which $xy\prec y$
we have $v\prec\phi_{I}\left(v\right)$. The reduction $xy\rightarrow y$
using $x-1$, which under shortlex would constitute a valid div-reduction
in Rosenmann's algorithm \cite{Rosenmann1993} or a prefix-reduction
in Madlener-Reinert \cite{Madlener1993}, here produces a strict $\prec_{\text{max}}$-increase
in the support.
\end{example}

\begin{figure}
\centering{}\label{fig: no inclusion relation between T_I and T_J}\includegraphics[scale=0.5]{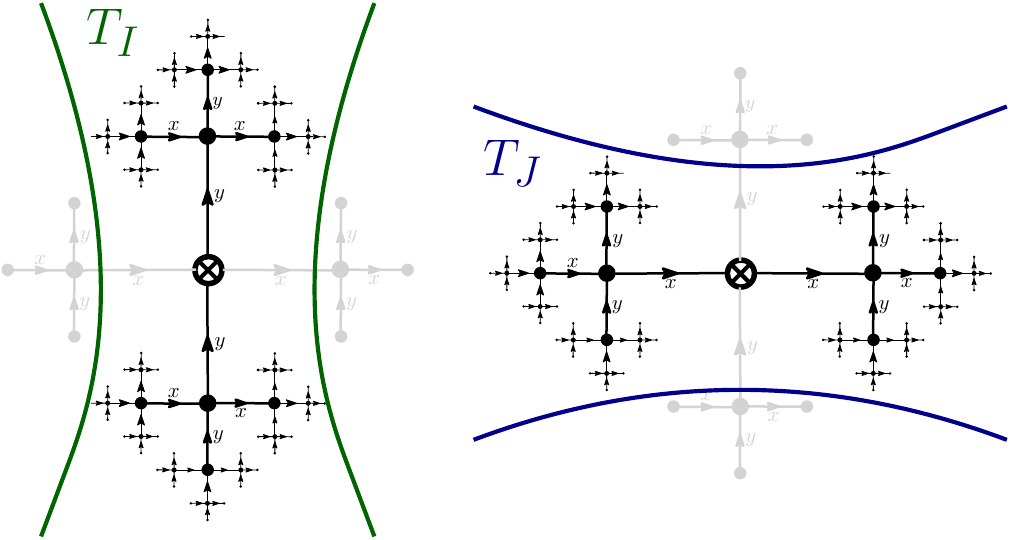}\caption{Illustration of Example \ref{exa: example with no inclusion relation when increasing ideal}:
The ideal $I=\left(x-1\right)\mathcal{A}$ contains $J=\left(x-1\right)\left(y-1\right)\mathcal{A}$,
but no inclusion relation exists between their respective minimal
Schreier transversals $T_{I}$ (on the left, in green) and $T_{J}$
(on the right, in blue). Edges and vertices outside of each Schreier
transversal are greyed out.}
\end{figure}

However, reducing a neighbor of $T_{I}$ does produce a $\prec_{\text{max}}$-decrease
in the support, and in fact, this property characterizes $T_{I}$,
as shown by the following proposition.
\begin{prop}
\label{prop: Reducing a neighbour of minimal Schreier Transversal T_I decreases support}Let
$T$ be a Schreier transversal for an ideal $I\leq\mathcal{A}$ with
associated transversal function $\phi:\mathcal{A}\rightarrow Sp_{K}T$.
Then $T=T_{I}$ if and only if for every neighbor $v$ of $T$ and
every $u\in\text{supp}\left(\phi\left(v\right)\right)$ we have $u\prec v$.
\end{prop}
\begin{proof}
\uline{(<=)} By minimality of $T_{I}$ we have $T_{I}\preceq_{\text{min}}T$.
Suppose for contradiction that $T\neq T_{I}$. Then $T_{I}\prec_{\text{min}}T$,
and let $v=\min\left(T_{I}\triangle T\right)\in T_{I}\backslash T$.
By Claim \ref{claim: minimum of symmetric difference is neighbour of larger tree},
$v$ is a neighbor of $T$. By the assumption, each $u\in\text{supp}\left(\phi\left(v\right)\right)$
satisfies $u\prec v$. Since $u\in T$, this implies by minimality
of $v$ that $u\in T_{I}$ as well. Hence, $v-\phi\left(v\right)$
is a nonzero element of $I$ supported on $T_{I}$, in contradiction.\\
\uline{(=>)} Assume that $T=T_{I}$ and suppose toward contradiction
that there exists a prefix-neighbor $v$ of $T$ and a word $u'\in\text{supp}\left(\phi\left(v\right)\right)$
such that $u'\succeq v$. Let $u$ be the maximal word in the support
of $\phi_{I}\left(v\right)$. Then $u\succeq u'\succeq v$, and in
fact $u\succ v$ since $u\in T$ but $v\notin T$. Denote by $F_{u}$
the set of words in $F$ having $u$ as a prefix (See the right part
of Figure \ref{fig: vertex larger than its reminder and proof}).
Note that $v\notin F_{u}$ since $u\succ v$. We claim that $T'=T_{I}\cup\left\{ v\right\} \backslash F_{u}$
is a partial Schreier transversal for $I$. It is still prefix-closed
since the added vertex $v$ is a prefix-neighbor of the original prefix-closed
set $T_{I}$. Let $f\in\text{Sp}_{K}\left(T'\right)\cap I$ and write
$f=\sum_{w\in T}\alpha_{w}w$ for coefficients $\alpha_{w}\in K$.
Since $f$ and $\alpha_{v}\left(v-\phi_{I}\left(v\right)\right)$
are elements of $I$ with the same coefficient behind $v$, their
difference satisfies $f-\alpha_{v}\left(v-\phi_{I}\left(v\right)\right)\in I\cap\text{Sp}_{K}\left(T_{I}\right)=\left\{ 0\right\} $.
Hence, $f=\alpha_{v}\left(v-\phi_{I}\left(v\right)\right).$ Comparing
the coefficient of $u$ on both sides -- nonzero in $v-\phi_{I}\left(v\right)$
and zero in $f$ -- we must have $\alpha_{v}=0$. But then $f\in I\cap\text{Sp}_{K}\left(T_{I}\right)$
so $f=0$. Thus, $T'$ is indeed a partial Schreier transversal for
$I$. But now 
\[
\min\left(T_{I}\triangle T'\right)=\min\left(\left\{ v\right\} \cup F_{u}\right)=v\in T',
\]
so $T'\prec_{\text{min}}T_{I}$ in contradiction to Corollary \ref{cor: T_I  is minimal among partial Schreier transversals for I}.
\end{proof}
\begin{claim}
\label{claim: ideal containment implies smaller schreier transversal}Let
$J\subseteq I$ be ideals. Then $T_{J}\preceq_{\text{min}}T_{I}$.
Moreover, if $J\subsetneqq I$, then $T_{J}\prec_{\text{min}}T_{I}$.
\end{claim}
\begin{proof}
$T_{I}$ is prefix-closed and satisfies $\text{Sp}_{K}\left(T_{I}\right)\cap J\subseteq\text{Sp}_{K}\left(T_{I}\right)\cap I=\left\{ 0\right\} $,
so $T_{I}$ is a partial Schreier transversal for $J$. By Corollary
\ref{cor: T_I  is minimal among partial Schreier transversals for I},
$T_{J}\preceq_{\text{min}}T_{I}$. Suppose now that $J\subsetneqq I$,
but $T_{J}=T_{I}$. Then there exists some $f\in I\backslash J$.
Consider $\phi_{J}\left(f\right)$, the remainder of $f$ modulo $T_{J}$.
By definition, $\phi_{J}\left(f\right)\in f+J\subseteq I$, and $\phi_{J}\left(f\right)$
is supported on $T_{J}=T_{I}$. Since $f\notin J$, it follows that
$\phi_{J}\left(f\right)\neq0$. But this contradicts the fact that
$T_{I}$ is a Schreier transversal for $I$.
\end{proof}
Although $J\subseteq I$ implies $T_{J}\preceq_{\text{min}}T_{I}$,
it does not follow that $T_{J}$ contains $T_{I}$, as shown in the
following example.
\begin{example}
\label{exa: example with no inclusion relation when increasing ideal}Let
$I$ be generated by $x-1$ as in Example \ref{exa: word smaller than its phi}.
Recall that $T_{I}$ is the set of words which do not have $x$ or
$x^{-1}$ as a prefix, regardless of the exposure order (Figure \ref{fig: no inclusion relation between T_I and T_J},
left). Let $J\subseteq I$ be the ideal generated by $\left(x-1\right)\left(y-1\right)=xy-y-x+1$.
Choose some exposure order on $F$ for which the smallest elements
are 
\[
e\prec x\prec xy\prec xy^{-1}\prec y\prec y^{-1}\prec...
\]
Then\footnote{We do not prove this for now as Algorithm \ref{alg:exposure_grobner}
allows computing $T_{J}$.} $T_{J}$ is the set of all words in $F$ which do not begin with
$y$ or $y^{-1}$ (Figure \ref{fig: no inclusion relation between T_I and T_J},
right). Thus, $y\in T_{I}\backslash T_{J}$ so $T_{I}\nsubseteq T_{J}$.
\end{example}

We finish this section with the following result, which illustrates
the flexibility of exposure orders.
\begin{prop}
\label{prop: Every Schreier tranversal is minimal with respect to some exposure order}Let
$T$ be a Schreier transversal for a right ideal $I\leq\mathcal{A}$,
and let $\prec_{T}$ be an exposure order such that $u\prec_{T}v$
for every $u\in T$ and $v\in F\backslash T$. Then $T_{I,\prec_{^{T}}}=T$.
In particular, every Schreier transversal for $I$ is its minimal
Schreier transversal with respect to some exposure order on $F$.
\end{prop}
\begin{proof}
For the first assertion, let $T'$ be another Schreier transversal
for $I$, distinct from $T$. By Claim \ref{claim: A Schreir transversal is inclusion-maximal among partial Schreier transversals},
two distinct Schreier transversals of the same ideal cannot contain
each other. Therefore, there exists some $u\in T\backslash T'$. By
the assumption on $\prec_{T}$, this word $u$ is smaller than any
word in $T'\backslash T$, so the $\prec_{T}$-minimum of $T\triangle T'$
must lie in $T\backslash T'$. Thus, $T\prec_{T,\text{min}}T'$ as
required.

For the second assertion, it suffices to show that there exists an
exposure order $\prec_{T}$ which satisfies the assumptions of the
first assertion. Let $\prec$ denote some exposure order (shortlex,
for example), and let $\prec_{T}$ denote the exposure order obtained
from $\prec$ by Claim \ref{claim: Sum of two order on tree and complement is an exposure order too},
applied to the prefix-closed subset $T$. By construction, any word
in $T$ is $\prec_{T}$-smaller than any word in its complement $F\backslash T$.
\end{proof}

\subsection{\label{sec: Highest Terms and Minimal Monic Elements}Head Terms
and Minimal Monic Elements of an Ideal}

Recall that we are working with a fixed exposure order $\prec$ on
$F$. In this section, we show how the derived well-order $\prec_{\text{max}}$
on finite subsets allows us to select elements with minimal support.
\begin{defn}
\textbf{(Head Term \& Monic Elements) }Let $0\neq f\in\mathcal{A}$.
The \emph{head term} \emph{of $f$}, denoted $\text{HT}_{\prec}\left(f\right)$,
is the $\prec$-maximal word in the support of $f$. The element $f$
is called \emph{monic} if the coefficient of $\text{HT}_{\prec}\left(f\right)$
in $f$ is $1$. Denote by $\text{MONIC}_{\prec}\left(f\right)$ the
monic $K^{\times}$-multiple of $f$.
\end{defn}
As we are working with a fixed exposure order, we omit it from the
notation and simply write $\text{HT}\left(f\right)$ and $\text{MONIC}\left(f\right)$. 
\begin{prop}
\label{prop: unique monic element supported on T with minimal support}Let
$I'\subsetneqq I$ be ideals in $\mathcal{A}$ and let $T'\subseteq F$
such that $I\backslash I'$ contains an element supported on $T'$.
Then there exists a unique monic element $f\in I\backslash I'$ supported
on $T'$ with $\prec_{\text{\ensuremath{\max}}}$-minimal support.
\end{prop}
\begin{rem}
Let $I'\subsetneqq I$. The following two choices for $T'$ satisfy
the assumptions of Proposition \ref{prop: unique monic element supported on T with minimal support}:
\begin{enumerate}
\item $T'=F$. In this case any element of $I\backslash I'$ is supported
on $T'$.
\item $T'$ is a Schreier transversal for $I'$. Then for any element of
$I\backslash I'$, its $I'$-coset is contained in $I\backslash I'$,
and contains an element supported on $T'$. In this case, a monic
element supported on $T'$ cannot belong to $I'$ since it is nonzero
(by being monic). The element $f$ produced by the proposition is
then more simply described as the unique monic element of $I$ of
minimal support among those supported on $T'$.
\end{enumerate}
\end{rem}
\begin{proof}
Let $X=\left\{ \text{supp}\left(f\right):f\in I\backslash I',\ \text{supp}\left(f\right)\subseteq T'\right\} $.
By the assumption on $T'$, the set $X$ is nonempty. By Claim \ref{claim: The order on finite subsetes is a well-order},
the order $\prec_{\text{max}}$ is a well-order on finite subsets
of $F$, so $X$ contains a $\prec_{\text{max}}$-minimal element
$S$. Choose an element $f\in I\backslash I'$ such that $\text{supp}\left(f\right)=S$.
Since $f\notin I'$, it is nonzero and can therefore be assumed to
be monic (if it is not already, replace $f$ with its monic $K^{\times}$-multiple).
To prove uniqueness, suppose $g\in I\backslash I'$ is another monic
element with $\text{supp}\left(g\right)=S$. Consider the difference
$h=f-g$. Since both $f$ and $g$ are elements of $I$ supported
on $S$, so is $h$. Moreover, the support of $h$ is strictly contained
in $S$ since $\max S$ cancels in the difference $f-g$. The minimality
of $S$ then implies that $h\in I'$. If $h\neq0$ then for a suitable
scalar $\lambda\in K$, the element $f-\lambda h$ lies in $I\backslash I'$
and has support strictly $\prec_{\text{max}}$-smaller than $S$,
in contradiction. Thus, $h=0$ and so $g=f$.
\end{proof}
\begin{prop}
\label{prop: HT(f) supported on neighbour of T_I}Let $f$ be a nonzero
element of $I$ which is supported on $T_{I}$ and its prefix-neighbors.
Then $\text{HT}\left(f\right)$ is a prefix-neighbor of $T_{I}$,
namely, $\text{HT}\left(f\right)\notin T_{I}$.
\end{prop}
\begin{proof}
Without loss of generality, assume $f$ is monic. Let $w=\text{HT}\left(f\right)$
and write $f'=f-w$. Since $f\in I,$ applying the transversal function
$\phi_{I}$ and using its $K$-linearity gives $\phi_{I}\left(f'\right)=-\phi_{I}\left(w\right).$
Assume for contradiction that $w\in T_{I}$. Then the equation becomes
$\phi_{I}\left(f'\right)=-w$. But $f'$ is supported on words smaller
than $w$, and by Proposition \ref{prop: Reducing a neighbour of minimal Schreier Transversal T_I decreases support},
so is $\phi_{I}\left(f'\right)$, in contradiction. It follows that
$w\notin T_{I}$.
\end{proof}

\subsection{\label{sec: The Exposure Basis}The Exposure Basis of an Ideal}

In this section we describe a process for constructing a basis for
an ideal $I\leq\mathcal{A}$, which we call the \emph{exposure basis}
of $I$, with respect to the fixed exposure order $\prec$. Similarly
to Example \ref{exa: vector space example}, the process adds basis
elements one after the other, by means of transfinite recursion. At
each step, after having added basis elements which generate a sub-ideal
$I'\subsetneqq I$, the new basis element added is a monic element
of $I$ with minimal support among those supported on $T_{I'}$. In
Section \ref{sec: The Exposure Basis is a Basis} we will prove Theorem
\ref{thm: The exposure basis is a basis}, showing that this procedure
indeed yields a basis for $I$.

\begin{figure}
\centering{}\label{fig: form for first}\includegraphics[scale=0.4]{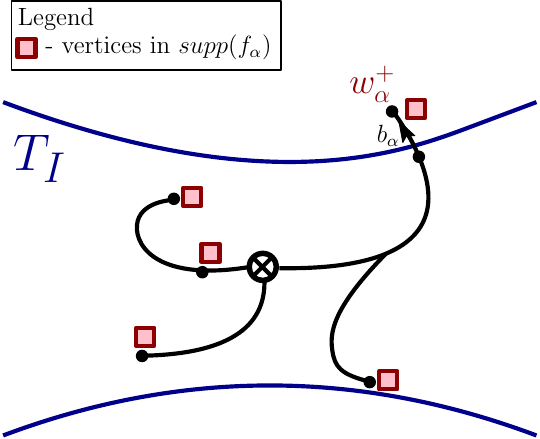}\caption{Schematic illustration of the support of an exposure element $f_{\alpha}$,
as in Proposition \ref{prop: Form for first using T_I}. The support
of $f_{\alpha}$ is indicated by red squares adjacent to the corresponding
vertices. The final letter $b_{\alpha}$ of the head term $w_{\alpha}^{+}$
exits $T_{I}$, making $w_{\alpha}^{+}$ its prefix-neighbor. The
rest of $\text{supp}\left(f_{\alpha}\right)$ lies within $T_{I}$.}
\end{figure}

The elements of the exposure basis are indexed by ordinal numbers.
For a given ordinal $\alpha$, the element $f_{\alpha}$ is only defined
after the elements $\left\{ f_{\beta}:\beta<\alpha\right\} $ have
already been defined. A stopping condition is provided for the process:
when the basis elements already generate $I$. Proposition \ref{prop:The-process-stabilizes}
shows that this stopping condition must occur before reaching the
first uncountable ordinal $\omega_{1}$. 
\begin{defn}
\label{def: inductive definition for ideal}Let $I\leq\mathcal{A}$
and let $\alpha$ be an ordinal, and suppose that we have already
defined $f_{\beta}=f_{\beta}\left(I,\prec\right)$ for every $\beta<\alpha$.
\begin{enumerate}
\item \textbf{(Current Ideal)} Let $I_{\alpha}=I_{\alpha}\left(I,\prec\right)$
be the ideal generated by $\left\{ f_{\beta}:\beta<\alpha\right\} $.
\item \textbf{(Stopping Condition \& Ordinal of $I$)} If $I_{\alpha}=I$
stop the recursive definition and call $\alpha$ the \emph{ordinal
of $I$} (with respect to the fixed order on $F$), which we denote
by $\omega_{I}=\omega_{I,\prec}$.
\item \textbf{(}\texttt{\textbf{First}}\textbf{ of $\alpha$)} If $I_{\alpha}\neq I$,
let $f_{\alpha}=f_{\alpha}\left(I,\prec\right)$ be the unique monic
element of $I$ having minimal support among those supported on $T_{I_{\alpha}}$,
as in Proposition \ref{prop: unique monic element supported on T with minimal support}.
We refer to $f_{\alpha}$ as the \texttt{first} of $\alpha$ and denote
its head term by $w_{\alpha}^{+}$.\footnote{A \texttt{second} for $\alpha$ will later be defined, with head term
$w_{\alpha}^{-}$. The terminology of \texttt{first} and \texttt{second}
originates from the work of Rosenmann in \cite{Rosenmann1993}.}
\end{enumerate}
To clarify, $f_{\alpha}$ is never defined for $\alpha\geq\omega_{I}$.
This process reflects a broader philosophy -- appearing also in \cite[Chapter~3]{ErnstWest2024}
-- of inductively \emph{exposing} elements of the ideal $I$ in accordance
with the fixed order on $F$. Here, at stage $\alpha$, having already
exposed elements generating a proper sub-ideal $I_{\alpha}\subsetneqq I$,
we seek a monic element $f_{\alpha}\in I$ with $\prec_{\text{max}}$-minimal
support among those supported on the current minimal Schreier transversal
$T_{I_{\alpha}}$. The minimal Schreier transversal of the \emph{updated}
ideal $I_{\alpha+1}$ generated by $I_{\alpha}$ and $f_{\alpha}$
will be shown to be obtained from $T_{I_{\alpha}}$ by truncating
precisely two branches from it: those lying behind the head terms
of $f_{\alpha}$ and its associated \texttt{second}. Thus, the process
can be viewed as a recursive refinement, repetitively seeking and
then pruning, until the ideal $I$ is fully exposed. 

We show that the process must stop at some countable ordinal by applying
the pigeonhole principle to the set of possible supports.
\end{defn}
\begin{prop}
\label{prop:The-process-stabilizes}The inductive process of Definition
\ref{def: inductive definition for ideal} terminates before reaching
the first uncountable ordinal $\omega_{1}$.
\end{prop}
\begin{proof}
Suppose otherwise, and consider the supports of $\left(f_{\alpha}\right)_{\alpha<\omega_{1}}$.
Since there are only countably many possible supports (finite subsets
of the countable set $F$), there exist two ordinals $\beta<\alpha<\omega_{1}$
such that $f_{\beta}$ and $f_{\alpha}$ have the same support. Consider
their difference $h=f_{\alpha}-f_{\beta}$ . Since $f_{\alpha}$ and
$f_{\beta}$ are both monic, the support of $h$ is strictly contained
in the support of $f_{\alpha}$. Hence, $h$ is supported on $T_{I_{\alpha}}$
as well but satisfies $\text{supp}h\prec\text{supp}f_{\alpha}$. Since
$f_{\beta}\in I_{\alpha}$ and $f_{\alpha}\in I\backslash I_{\alpha}$,
we have $h\in I\backslash I_{\alpha}$, in contradiction to the minimality
of $f_{\alpha}$.
\end{proof}
\begin{rem}
Although the results in this paper are stated for finitely generated
free groups for simplicity, they extend to infinitely generated ones
as well. The only Proposition requiring modification is Proposition
\ref{prop:The-process-stabilizes}, where replacing $\omega_{1}$
with an ordinal of cardinality strictly greater than that of $F$
yields an analogous result.
\end{rem}
\begin{example}
\label{exa: exposure process}We give several examples for the process
of Definition \ref{def: inductive definition for ideal}:
\begin{enumerate}
\item Let $I=\left\{ 0\right\} $ be the zero ideal. Then the process stops
already at the first ordinal, so $\omega_{I}=0$.
\item Let $I=\mathcal{A}$ be the improper ideal. Then $f_{0}$ is the unit
element $e\in K\left[F\right]$, since $e$ is monic and with $\prec_{\text{max}}$-minimal
support in every exposure order among the nonzero elements of $\mathcal{A}$.
Since $e$ already generates $\mathcal{A}$ as a right ideal, we have
$I_{1}=I$ and the process stops with ordinal $\omega_{I}=1$.
\item Let $F=\left\langle x,y\right\rangle $ be the free group on two generators
and let $I\leq K\left[F\right]$ be the augmentation ideal, i.e.,
the set of elements $\sum_{u\in F}\alpha_{u}u$ with $\sum_{u\in F}\alpha_{u}=0$.
Choose some exposure order $\prec$ on $F$ in which the smallest
elements are $e\prec x\prec xy\prec y\prec...$ . Then $f_{0}=x-1$
since the only monic element with $\prec_{\text{max}}$-smaller support
is $e$ and $e\notin I$. \\
The ideal $I_{1}$ is then generated by the single element $f_{0}=x-1$.
Recall from Example \ref{exa: word smaller than its phi} that the
minimal Schreier transversal $T_{I_{1}}$ for this ideal consists
of all words which do not begin with either $x$ or $x^{-1}$. \\
Since again $e\notin I$, the minimal monic element of $I$ supported
on $T_{I_{1}}$ is $f_{1}=y-1$. The process thus terminates at $\omega_{I}=2$,
as $I$ is generated by $x-1$ and $y-1$. Observe that $f_{1}$ has
$\prec_{\max}$-minimal support only among monic elements of $I$
supported on $T_{I_{1}}$. For example,  the element $f_{1}':=xy-1$
is another monic element in $I\backslash I_{1}$ with strictly smaller
support. Moreover, $f_{1}'$ lies in the same $I_{1}$-coset as $f_{1}$,
demonstrating that $f_{1}$ was chosen during the process without
minimizing support even within its $I_{1}$-coset.
\item Let $I$ be the ideal generated by the elements $\left\{ x^{i}y-x^{i}:i\geq0\right\} $.
Suppose that in the exposure order on $F$, the generator $y$ is
larger than any word of the form $x^{n}$ for $n\geq0$ or $x^{n}y$
for $n>0$, and that the smallest words in the order are:
\[
e\prec x\prec xy\prec x^{2}\prec x^{2}y\prec x^{3}\prec x^{3}y\prec x^{4}\prec...\prec y\prec...
\]
For every $n=0,1,2,...$ the $n$-th minimal transversal $T_{I_{n}}$
in the exposure process is the set of words which do not begin with
$x^{i}y$ or $x^{i}y^{-1}$ for any $0<i<n$, and the $n$-th basis
element is $f_{n}=x^{n+1}y-x^{n+1}$. After exhausting all natural
numbers, upon reaching the first infinite ordinal $\omega=\left\{ 0,1,2,3,...\right\} $,
the ideal $I_{\omega}$ is still not $I$ as it is only generated
by $\left\{ x^{i}y-x^{i}:i>0\right\} $. Its minimal transversal is
the set of words which do not begin with $x^{i}y$ or $x^{i}y^{-1}$
for any positive $i$. We then have $f_{\omega}=y-1$ and the process
stops at its successor ordinal $\omega_{I}=\omega+1$.
\end{enumerate}
\end{example}
\begin{defn}
\textbf{\label{def: The Exposure Basis}(The Exposure Basis) }Let
$I\leq\mathcal{A}$ be an ideal. The \emph{Exposure basis} $B_{I}=B_{I,\prec}$
of $I$ is the set $\left\{ f_{\beta}:\beta<\omega_{I}\right\} $,
where $\omega_{I}$ is the ordinal of $I$ and $f_{\beta}$ is the
\texttt{first} associated to the ordinal $\beta$ in Definition \ref{def: inductive definition for ideal}.
\end{defn}
Let $\alpha<\omega_{I}$ and let $I_{\alpha}$ be the intermediate
ideal constructed at stage $\alpha$ of the exposure process of $I$.
The following technical proposition shows that the exposure process
of $I_{\alpha}$ coincides with the initial segment of the exposure
process of $I$ up to ordinal $\alpha$.
\begin{prop}
\label{claim: Exposure basis for an ideal along the way}Let $I\leq\mathcal{A}$
be an ideal and let $\alpha<\omega_{I}$. Then the exposure basis
of $I_{\alpha}$ is given by $\left\{ f_{\beta}:\beta<\alpha\right\} $.
\end{prop}
\begin{proof}
For every $\beta<\omega_{I_{\alpha}}$, let $f_{\beta}'$ and $I_{\beta}'$
denote, respectively, the element and ideal associated with the ordinal
$\beta$ when performing the inductive process of Definition \ref{def: inductive definition for ideal}
on the ideal $I_{\alpha}$. We will show by transfinite induction
that $f_{\beta}'$ is defined and $f_{\beta}'=f_{\beta}$ for every
$\beta<\alpha$. Let $\gamma<\alpha$ be an ordinal, and suppose that
$f_{\beta}'$ is defined and $f_{\beta}'=f_{\beta}$ for every $\beta<\gamma$.
Then $I_{\gamma}'=I_{\gamma}$, since both are generated by $\left\{ f_{\beta}:\beta<\gamma\right\} $. 

Now, $f_{\gamma}$ is defined as the unique monic element of $I$
supported on $T_{I_{\gamma}}$ with minimal support. Then $f_{\gamma}\notin I_{\gamma}$,
and since $\gamma<\alpha$ we have $f_{\gamma}\in I_{\alpha}$. In
particular, $I_{\gamma}\neq I_{\alpha}$ so the inductive process
for $I_{\alpha}$ does not terminate at ordinal $\gamma$. Furthermore,
since $I_{\alpha}\subseteq I$, $f_{\gamma}$ is also the unique monic
element of $I_{\alpha}$ supported on $T_{I_{\gamma}}$ with minimal
support. Hence, by definition, $f_{\gamma}'=f_{\gamma}$, completing
the inductive step.

Finally, at ordinal $\alpha$, we have $\left\{ f_{\beta}':\beta<\alpha\right\} =\left\{ f_{\beta}:\beta<\alpha\right\} =I_{\alpha}$,
so the process terminates at ordinal $\alpha$, and the exposure basis
for $I_{\alpha}$ is exactly $\left\{ f_{\beta}:\beta<\alpha\right\} $.
\end{proof}
The following Proposition shows that the elements of the exposure
basis have a specific form with respect to $T_{I}$, as illustrated
in Figure \ref{fig: form for first}. Recall that $w_{\alpha}^{+}$
denotes the head term of $f_{\alpha}$.
\begin{prop}
\label{prop: Form for first using T_I}Let $\alpha<\omega_{I}$. Then
$f_{\alpha}=w_{\alpha}^{+}-\phi_{I}\left(w_{\alpha}^{+}\right)$ and
$w_{\alpha}^{+}=\min\left(T_{I_{\alpha}}\triangle T_{I}\right)$.
In particular, $w_{\alpha}^{+}$ is a neighbor of $T_{I}$.
\end{prop}
\begin{proof}
Since $I_{\alpha}\subsetneqq I$, we have $T_{I_{\alpha}}\prec_{\text{min}}T_{I}$
by Claim \ref{claim: ideal containment implies smaller schreier transversal}.
Let $w=\min\left(T_{I_{\alpha}}\triangle T_{I}\right)\in T_{I_{\alpha}}$.
By Claim \ref{claim: minimum of symmetric difference is neighbour of larger tree},
$w$ is a neighbor of $T_{I}$. Then, by Proposition \ref{prop: Reducing a neighbour of minimal Schreier Transversal T_I decreases support},
the support of $\phi_{I}\left(w\right)$ consists of words strictly
smaller than $w$. Since $\text{supp}\left(\phi_{I}\left(w\right)\right)\subseteq T_{I}$,
and $w$ is the minimal element of $T_{I_{\alpha}}\triangle T_{I}$,
it follows that $\text{supp}\left(\phi_{I}\left(w\right)\right)\subseteq T_{I_{\alpha}}$.
Now define $f:=w-\phi_{I}\left(w\right)$. Then $f\in I$ is monic
and supported on $T_{I_{\alpha}}.$ We will deduce that $f_{\alpha}=f$
by showing that the support of $f$ is $\prec_{\text{max}}$-minimal
among monic elements of $I$ supported on $T_{I_{\alpha}}$. Let $g\in I$
be another monic element supported on $T_{I_{\alpha}}.$ Since $g\in I\backslash\left\{ 0\right\} $,
there exists some $v\in\text{supp}\left(g\right)\backslash T_{I}$.
Since $v\in T_{I_{\alpha}}\backslash T_{I}$, by the minimality of
$w$ we have $v\succeq w$, so $\text{HT}\left(g\right)\succeq w=\text{HT}\left(f\right)$.
If $\text{HT}\left(g\right)\succ w$ then $\text{supp}\left(g\right)\succ_{\max}\text{supp}\left(f\right)$
as required. If instead $\text{HT}\left(g\right)=w$, then $g-w$
is supported on elements of $T_{I_{\alpha}}$ strictly smaller than
$w$, and hence its support is contained in $T_{I}$. Therefore, $g-f=\left(g-w\right)+\phi_{I}\left(w\right)$
is an element of $I$ supported on its Schreier transversal $T_{I}$,
so $g-f=0$.
\end{proof}

The inductive process thus selects at stage $\alpha$ the minimal
prefix-neighbor $w_{\alpha}^{+}$ of $T_{I}$ lying in $T_{I_{\alpha}}$,
and adds to the exposure basis the element $f_{\alpha}=w_{\alpha}^{+}-\phi_{I}\left(w_{\alpha}^{+}\right)$.
We next show that these head terms are selected in increasing order.
\begin{prop}
\label{prop: Firsts are increasing}Let $\beta<\alpha<\omega_{I}$.
Then $w_{\beta}^{+}\prec w_{\alpha}^{+}$.
\end{prop}
\begin{proof}
Since $f_{\beta}\in I_{\alpha}$ but $f_{\alpha}\notin I_{\alpha}$,
we conclude that $f_{\alpha}\neq f_{\beta}$. By Proposition \ref{prop: Form for first using T_I},
their respective head terms $w_{\alpha}^{+}$ and $w_{\beta}^{+}$
must be distinct. Moreover, that proposition gives: 
\[
w_{\alpha}^{+}=\min\left(T_{I_{\alpha}}\triangle T_{I}\right),\ w_{\beta}^{+}=\min\left(T_{I_{\beta}}\triangle T_{I}\right).
\]
 Applying Claim \ref{claim: minimum of symmetric difference}, we
obtain:
\[
\min\left(T_{I_{\alpha}}\triangle T_{I_{\beta}}\right)=\min\left(\left(T_{I_{\alpha}}\triangle T_{I}\right)\triangle\left(T_{I}\triangle T_{I_{\beta}}\right)\right)=\min\left\{ w_{\alpha}^{+},w_{b}^{+}\right\} .
\]
Now, since $I_{\beta}\subsetneqq I_{\alpha}$, Claim \ref{claim: ideal containment implies smaller schreier transversal}
implies that $\min\left(T_{I_{\alpha}}\triangle T_{I_{\beta}}\right)\in T_{I_{\beta}}\backslash T_{I_{\alpha}}$.
In particular, this minimum cannot be $w_{\alpha}^{+}$, which lies
in $T_{I_{\alpha}}.$ It follows that $\min\left\{ w_{\alpha}^{+},w_{\beta}^{+}\right\} =w_{\beta}^{+}$,
so $w_{\beta}^{+}\prec w_{\alpha}^{+}$.
\end{proof}

\subsection{\label{sec: The Seconds}The \texttt{Seconds}}

From this point and until the end of Section \ref{sec: Combinatorial Properties of the grobner basis},
in addition to the fixed exposure order $\prec$ on $F$, we fix a
proper ideal $I\lneqq\mathcal{A}$.

As a prefix-closed subset of $F$, the minimal Schreier transversal
$T_{I}$ can be characterized via its set of prefix-neighbors in the
Cayley graph of $F$. Specifically, if $\partial T_{I}$ denotes the
set of prefix-neighbors of $T_{I}$, then $T_{I}$ is precisely the
set of words in $F$ which do not have any $v\in\partial T_{I}$ as
a prefix. By Proposition \ref{prop: Form for first using T_I}, the
head terms of the elements in the exposure basis $B_{I}$ all lie
in $\partial T_{I}$. However, we will show (see Theorem \ref{Thm: properties of second + Grobner basis is combinatorially reducing})
that unless $I=\mathcal{A}$, these head terms do not account for
all the prefix-neighbors of $T_{I}$. The following example illustrates
this phenomenon.
\begin{example}
\label{exa: why we need seconds}Let the ordering $\prec$ of $F$
be such that the minimal words satisfy $e\prec x\prec...$, and let
$I=\left(x-1\right)\mathcal{A}$ be the ideal considered in Examples
\ref{exa: word smaller than its phi} and \ref{exa: example with no inclusion relation when increasing ideal}.
Then $T_{I}$ consists of all words in $F$ which do not have $x$
or $x^{-1}$ as a prefix. The exposure basis $B_{I}$ consists of
a single element $f_{0}=x-1$, whose head term is $w_{0}^{+}=x$,
which is indeed a neighbor of $T_{I}$. However, the other neighbor
of $T_{I}$, namely $x^{-1}$, does not appear as a head term of any
exposure basis element. To account for it, we associate to $f_{0}$
the element $s_{0}=x^{-1}-1$ which we call its \texttt{second}. The
set $\left\{ f_{0},s_{0}\right\} \subseteq I$ then characterizes
$T_{I}$ in the sense that $T_{I}$ is obtained from the Cayley graph
by removing the head terms of $f_{0}$ and $s_{0}$, along with the
subtrees behind them.
\end{example}
From this point onward we assume that $I\lneqq\mathcal{A}$$.$ Generalizing
the previous example, we associate to each exposure basis element
$f_{\alpha}$ an additional element $s_{\alpha}\in I$, called its
\texttt{second}. As in the example, we will show that the head terms
of the \texttt{firsts} $\left(f_{\alpha}\right)_{\alpha<\omega_{I}}$
and \texttt{seconds} $\left(s_{\alpha}\right)_{\alpha<\omega_{I}}$are
distinct words in $F$, which together form the set of all prefix-neighbors
of $T_{I}$ (see Theorem \ref{Thm: properties of second + Grobner basis is combinatorially reducing}). 

The general idea behind the construction of $s_{\alpha}$ is as follows:
The goal in mind is to compute the Schreier transversal $T_{I_{\alpha+1}}$
from the previous one, $T_{I_{\alpha}}$. By Proposition \ref{prop: Form for first using T_I},
the term $w_{\alpha}^{+}$ is not in $T_{I_{\alpha+1}}$. We therefore
remove it from $T_{I_{\alpha}}$, along with the entire subtree behind
it, forming a new subtree $T'\subseteq T_{I_{\alpha}}$. 

However, as demonstrated in Example \ref{exa: why we need seconds},
this truncation might not suffice: certain nonzero elements of $I_{\alpha+1}$
can still be supported on $T'$. To construct such an element explicitly,
observe that $f_{\alpha}$ is supported entirely on $T'$ except for
$w_{\alpha}^{+}$, which lies just outside $T'$ as a neighbor. Let
$b_{\alpha}$ denote the last letter of $w_{\alpha}^{+}$. We now
``push'' this head term back into $T'$ by multiplying $f_{\alpha}$
on the right by $b_{\alpha}^{-1}$, yielding the element $f_{\alpha}b_{\alpha}^{-1}\in I_{\alpha+1}$. 

A new problem then arises: $f_{\alpha}b_{\alpha}^{-1}$ might not
be supported on $T_{I_{\alpha}}$. Applying $\phi_{I_{\alpha}}$ to
this element yields another element of $I_{\alpha+1}$ which is not
only supported on $T_{I_{\alpha}}$, but, somewhat surprisingly, on
its truncated subtree $T'$. The monic $K^{\times}$-multiple of $\phi_{I_{\alpha}}\left(f_{\alpha}b_{\alpha}^{-1}\right)$
is the desired \texttt{second} $s_{\alpha}$. 

We hope the preceding discussion has provided some intuition. We now
formalize these ideas.
\begin{defn}
\textbf{(HTT - Head Term Tail) }Let $f\in\mathcal{A}\backslash K$.
The \emph{head term tail of $f$,} denoted $\text{HTT}\left(f\right)=\text{HTT}_{\prec}\left(f\right)$,
is the last letter of the (reduced) word $\text{HT}_{\prec}\left(f\right)$.
\end{defn}
Again, since we work with a fixed exposure order, we omit the dependence
of the head term tail on the fixed order and simply write $\text{HTT}\left(f\right)$. 

In order to define the \texttt{second}, observe that for any $f\in\mathcal{A}\backslash K$
supported on $T_{I}$, if we denote $b=\text{HTT}\left(f\right)$,
then $\phi_{I}\left(fb^{-1}\right)\neq0$. Indeed, since $f$ is nonzero
and supported on $T_{I}$, it cannot lie in $I$. Therefore, $fb^{-1}\notin I$
implying that its remainder modulo $T_{I}$ is nonzero.
\begin{defn}
\textbf{\label{def: I-second}($I$-}\texttt{\textbf{second}}\textbf{)}
Let $f\in\mathcal{A}\backslash K$ be supported on $T_{I}$, and let
$b=\text{HTT}\left(f\right)$. The \emph{$I$-}\texttt{second} \emph{of
$f$} is the monic $K^{\times}$-multiple of $\phi_{I}\left(fb^{-1}\right)$.
\end{defn}
Note that in the definition above, both $f$ and its $I$-\texttt{second}
are supported on $T_{I}$: $f$ by assumption, and the $I$-\texttt{second}
because it is a $K$-multiple of the output of $\phi_{I}$. 

Since $I\neq\mathcal{A}$, every nonzero $f\in I$ has a well-defined
head term tail. We can therefore name the head term tails of the exposure
basis elements and define their associated \texttt{seconds}.
\begin{defn}
\textbf{($b_{\alpha}$, $s_{\alpha}$ and $w_{\alpha}^{-}$)\label{def: (HLL of first, second, highest term of second)}}
Let $I\lneqq\mathcal{A}$. For every $\alpha<\omega_{I}$, let $b_{\alpha}$
be the head term tail of $f_{\alpha}$. Define $s_{\alpha}$ to be
the $I_{\alpha}$-\texttt{second} of $f_{\alpha}$, also called the
\texttt{second} \emph{of} $f_{\alpha}$. Finally, let $w_{\alpha}^{-}$
denote the head term of $s_{\alpha}$.
\end{defn}
\begin{example}
\label{examples of seconds}We provide several examples illustrating
Definition \ref{def: (HLL of first, second, highest term of second)}.
\end{example}
\begin{enumerate}
\item \label{enu: principal ideal example}Recall the ideal $I=\left(x-1\right)\mathcal{A}$
from example \ref{exa: why we need seconds}. We verify that Definition
\ref{def: (HLL of first, second, highest term of second)} indeed
yields $s_{0}=x^{-1}-1$. The ideal $I_{0}$ is the zero ideal, so
$T_{I_{0}}=F$ and $\phi_{I_{0}}$ is the identity map on $\mathcal{A}$.
By Definition \ref{def: (HLL of first, second, highest term of second)},
we indeed get:
\[
s_{0}=\text{MONIC}\left(f_{0}b_{0}^{-1}\right)=\text{MONIC}\left(\left(x-1\right)x^{-1}\right)=x^{-1}-1.
\]
\item We give a general description of $s_{0}$ for any nonzero ideal $I\lneqq\mathcal{A}$.
Recall that $f_{0}$ is the unique monic element of $I$ with minimal
support (since the demand of being supported on $T_{I_{0}}=F$ is
vacuous for for all elements of $\mathcal{A}$). As in item \ref{enu: principal ideal example},
$s_{0}$ is the monic $K^{\times}$-multiple of $f_{0}b_{0}^{-1}$,
so $\text{supp}\left(s_{0}\right)=\text{supp}\left(f_{0}\right)b_{0}^{-1}$.
By the minimality of $f_{0}$, $\text{supp}\left(s_{0}\right)\succeq_{\text{max}}\text{supp}\left(f_{0}\right)$,
so $w_{0}^{-}=\text{HT}\left(s_{0}\right)\succeq\text{HT}\left(f_{0}\right)=w_{0}^{+}$.
Consider any $v\in\text{supp}\left(s_{0}\right)$ such that $v\succeq w_{0}^{+}$
(for example, $v=w_{0}^{-})$. Then $vb_{0}\in\text{supp}\left(f_{0}\right)$,
so $vb_{0}\preceq\text{HT}\left(f_{0}\right)=w_{0}^{+}\preceq v$.
It follows that $v$ cannot be a proper prefix of $vb_{0}$, and so
$v$ must end with the letter $b_{0}^{-1}$. In particular, setting
$v=w_{0}^{-}$, we conclude that $\text{HTT}\left(s_{0}\right)=b_{0}^{-1}$.
This relation $\text{HTT}\left(s_{\alpha}\right)=\text{HTT}\left(f_{\alpha}\right)^{-1}$,
which we have shown here for $\alpha=0$, holds for all ordinals $\alpha<\omega_{I}$
-- see Theorem \ref{Thm: properties of second + Grobner basis is combinatorially reducing}.
\item \label{enu: involved example of seconds}We conclude with a more involved
example. We do not provide proofs for all the claims here. However,
all statements can be verified algorithmically, as will be explained
later. Let $K$ be the field of order $2$ and let $I=\left(y^{-2}+y+x\right)\mathcal{A}+\left(xy^{-1}+y\right)\mathcal{A}$.
Consider the shortlex order on $F$ such that $y^{-1}\prec x^{-1}\prec x\prec y$.
The minimal nonzero element of $I$ is $f_{0}=y^{-2}+y+x$, with head
term $w_{0}^{+}=y^{-2}$ and head term tail $b_{0}=y^{-1}$. Since
any nonzero element of $K\left[F\right]$ is monic, its associated
\texttt{second} is $s_{0}=f_{0}b_{0}^{-1}=y^{2}+xy+y^{-1}$, with
head term $w_{0}^{-}=y^{2}$. The minimal Schreier transversal $T_{I_{1}}$
of the ideal $I_{1}=f_{0}\mathcal{A}$ is the subtree with prefix-neighbors
$\left\{ y^{2},y^{-2}\right\} .$ The next \texttt{first} is $f_{1}=xy^{-1}+y$.
Multiplying by $y=b_{1}^{-1}$ yields $f_{1}y=y^{2}+x$, which is
not supported on $T_{I_{1}}$. However, subtracting $s_{0}$ from
$f_{1}y$ gives the representative supported on $T_{I_{1}}$ in its
$I_{1}$-coset: 
\[
s_{1}=\phi_{I_{1}}\left(f_{1}y\right)=xy+x+y^{-1}.
\]
The minimal Schreier transversal $T_{I}$ is then the subtree with
prefix-neighbors $\left\{ y^{2},yx,xy^{-1},xy\right\} .$ In correspondence
with Proposition \ref{prop: Firsts are increasing}, we have $w_{0}^{+}=y^{-2}\prec xy^{-1}=w_{1}^{+}$.
In contrast, we have $w_{0}^{-}=y^{2}\succ xy=w_{1}^{-}$, showing
that the \texttt{seconds} are not necessarily increasing. Furthermore,
as $w_{1}^{-}\in\text{supp}\left(s_{0}\right)$, the head term of
a ``later discovered'' \texttt{second} (or \texttt{first}) can appear
in the support of an ``earlier'' \texttt{second}. Note also that
the element $s_{0}-\text{HT}\left(s_{0}\right)$ is not supported
on $T_{I}$, unlike the form described in Proposition \ref{prop: Form for first using T_I}
which the \texttt{firsts} take. For a discussion of various ways to
define the seconds and the rationale behind our choice, see Remark
\ref{rem: why we define seconds this way}.
\end{enumerate}
\begin{figure}
\begin{centering}
\label{fig: involved example of seconds}\includegraphics[scale=0.4]{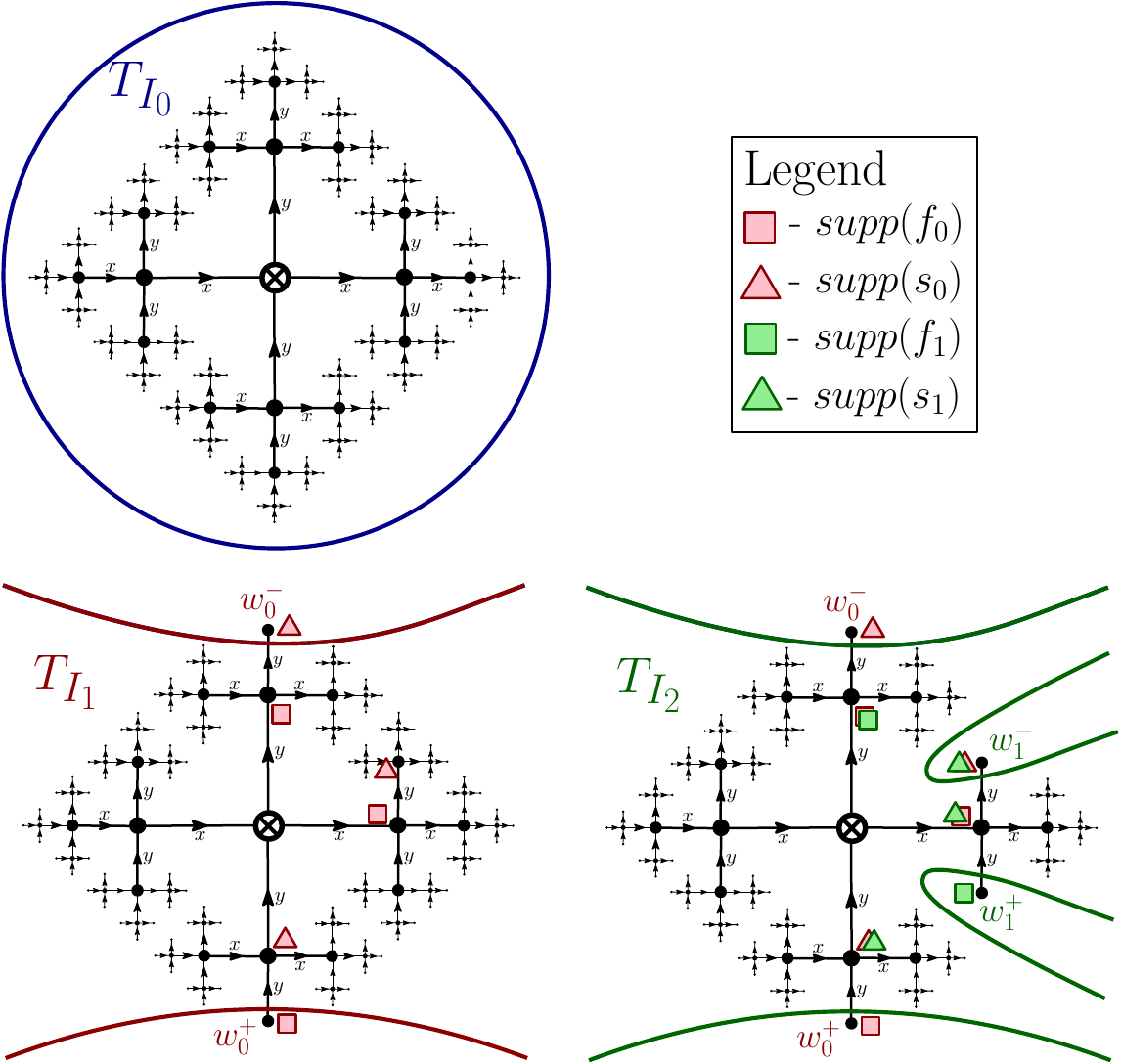}
\par\end{centering}
\caption{Illustration for Example \ref{examples of seconds} (\ref{enu: involved example of seconds}),
showing the evolution of the minimal Schreier transversal during the
two-stage exposure process of the ideal $I$:\protect \\
Top left: the initial Schreier transversal $T_{I_{0}}$ equals the
entire Cayley tree. \protect \\
Bottom left: following the exposure of $f_{0}=y^{2}+xy+y^{-1}$ and
its corresponding \texttt{second} $s_{0}=y^{-2}+y+x$, the Schreier
transversal $T_{I_{1}}$ excludes the prefixes $\left\{ y^{2},y^{-2}\right\} $.\protect \\
Bottom right: following the exposure of $f_{1}=xy^{-1}+y$ and its
corresponding \texttt{second} $s_{1}=xy+x+y^{-1}$, the Schreier transversal
$T_{I_{2}}=T_{I}$ excludes the prefixes $\left\{ y^{2},y^{-2},xy,xy^{-1}\right\} $.\protect \\
The supports of $f_{0},s_{0},f_{1},s_{1}$ are indicated by red squares,
red triangles, green squares, and \protect \\
green triangles, respectively.}

\end{figure}

\subsection{\label{sec: Combinatorially Reducing Systems}Combinatorially Reducing
Systems}

In this section, we establish sufficient combinatorial conditions
on a set $\mathcal{Q}\subseteq\mathcal{A}\backslash K$ which guarantee
that the subtree obtained by truncating their head terms is the minimal
Schreier transversal of the right ideal $I_{\mathcal{Q}}$ which they
generate. Moreover, we show that these conditions enable the reduction
of any element $f\in\mathcal{A}$ to $\phi_{I_{\mathcal{Q}}}\left(f\right)$
-- the unique representative of the coset $f+I_{\mathcal{Q}}$ supported
on this Schreier transversal -- in finitely many steps. In Section
\ref{sec: Combinatorial Properties of the grobner basis}, we will
prove that for any ideal $I\lneqq\mathcal{A}$, the \texttt{firsts}
$\left\{ f_{\alpha}:\alpha<\omega_{I}\right\} $ together with the
\texttt{seconds} $\left\{ s_{\alpha}:\alpha<\omega_{I}\right\} $
satisfy these combinatorial conditions.
\begin{defn}
\label{def: combinatorially reducing systems}Let $\mathcal{Q}$ be
a set of monic elements in $\mathcal{A}\backslash K$. For each $w\in F$,
denote by $F_{w}$ the set of words of $F$ having $w$ as a prefix.
Define the subtree 
\[
T_{\mathcal{Q}}^{\text{HT}}:=F\backslash\bigcup_{q\in\mathcal{Q}}F_{\text{HT}\left(q\right)}
\]
obtained by truncating all head terms of $\mathcal{Q}$ from the Cayley
tree. We say that $\mathcal{Q}$ is a \emph{combinatorially reducing}
\emph{system} or \emph{CRS}, with respect to the fixed order $\prec$,
if all the following conditions are satisfied:
\begin{enumerate}
\item The head terms $\left(\text{HT}\left(q\right)\right)_{q\in\mathcal{Q}}$
are distinct.
\item For each $q\in Q$, letting $b=\text{HTT}\left(q\right)$, the element
$qb^{-1}$ is a $K$-linear combination of elements $h\in Q$ with
$\text{HTT}\left(h\right)=b^{-1}$.
\item Each $q\in\mathcal{Q}$ is supported on $T_{\mathcal{Q}}^{\text{HT}}\cup\partial T_{\mathcal{Q}}^{\text{HT}}$,
where $\partial T_{\mathcal{Q}}^{\text{HT}}$ denotes the prefix-neighbors
of $T_{\mathcal{Q}}^{\text{HT}}$.
\end{enumerate}
\end{defn}
Note that if $\mathcal{Q}$ is a CRS then the head terms $\left\{ \text{HT}\left(q\right):q\in\mathcal{Q}\right\} $
are the (distinct) prefix-neighbors of $T_{\mathcal{Q}}^{\text{HT}}$.
In particular, $F_{\text{HT}\left(q_{1}\right)}\cap F_{\text{HT}\left(q_{2}\right)}=\emptyset$
for $q_{1}\neq q_{2}$ in $\mathcal{Q}$.

We describing a more precise form for the support of elements of $\mathcal{Q}.$
\begin{prop}
\label{prop: element of combinatorially reducing system is supported on T^HT =00005Ccup T^HT(b)}Let
$\mathcal{Q}\subseteq\mathcal{A}\backslash K$ be a CRS and let $q\in\mathcal{Q}$.
Then $q$ is supported on $T_{\mathcal{Q}}^{\text{HT}}\cup T_{\mathcal{Q}}^{\text{HT}}b$,
where $b=\text{HTT}\left(q\right)$.
\end{prop}
\begin{proof}
Let $u\in\text{supp}\left(q\right)\backslash T_{\mathcal{Q}}^{HT}$.
Since $\mathcal{Q}$ is a CRS, $u\in\partial T_{\mathcal{Q}}^{\text{HT}}$.
We claim that $u$ ends in $b$. Suppose otherwise. Then $u$ is a
proper prefix of $ub^{-1}$. Since $\mathcal{Q}$ is a CRS, $qb^{-1}$
is a $K$-linear combination of its elements, and each such element
is supported on $T_{\mathcal{Q}}^{\text{HT}}\cup\partial T_{\mathcal{Q}}^{\text{HT}}$.
Hence, $ub^{-1}\in T_{\mathcal{Q}}^{\text{HT}}\cup\partial T_{\mathcal{Q}}^{\text{HT}}$.
It follows that its proper prefix $u$ belongs to $T_{\mathcal{Q}}^{\text{HT}}$,
contradicting $u\in\partial T_{\mathcal{Q}}^{\text{HT}}$. Hence,
$u$ ends in $b$, so $ub^{-1}$ is a proper prefix of $u\in\partial T_{\mathcal{Q}}^{\text{HT}}$.
It follows that $ub^{-1}\in T_{\mathcal{Q}}^{\text{HT}}$, so $u\in T_{\mathcal{Q}}^{\text{HT}}b$,
as required.
\end{proof}
\begin{defn}
Let $I\leq\mathcal{A}$. We say that $\mathcal{Q}$ is a \emph{combinatorially
reducing system (CRS) for $I$}, with respect to the fixed order $\prec$,
if it is a combinatorially reducing system which generates $I$ as
a right ideal.
\end{defn}
\begin{lem}
\textbf{\label{lem: Canonical-Representation-using-1}(Canonical Representation
using a CRS)} Let $\mathcal{Q}$ be a combinatorially reducing system
for $I$, and let $f\in I$. Then:
\begin{enumerate}
\item \label{enu: canonical expression using grobner - existence-1}$f$
can be written as $f=\sum_{q\in\mathcal{Q}}qg_{q}$, for elements
$g_{q}\in\mathcal{A}$ such that for every $q\in\mathcal{Q}$, there
is no cancellation in $\text{\text{HT}}\left(q\right)\cdot g_{q}$,
i.e., $g_{q}$ is supported on words not beginning with $\text{\text{HTT}}\left(q\right)^{-1}$.
\item \label{enu: canonical expression using grobner - if some coefficient is nonzero then not supported on T_I'-1}If
in the representation above there exists some $q\in\mathcal{Q}$ with
$g_{q}\neq0$, then $f$ is not supported on $T_{\mathcal{Q}}^{\text{HT}}$.
\item \label{enu: canonical expression using grobner - coefficients are unique-1}The
coefficients $\left(g_{q}\right)_{q\in\mathcal{Q}}$ expressing $f$
as above are unique.
\end{enumerate}
\end{lem}
\begin{proof}
\begin{enumerate}
\item It suffices to prove the claim for $f=qp$ where $q\in\mathcal{Q}$
and $p\in\mathcal{A}$. Write $p=\sum_{u\in F}\lambda_{u}u\in\mathcal{A}$
for some scalars $\lambda_{u}\in K$. Let $b=\text{HTT}\left(q\right)$.
Decompose $p=g+h$, where $h=\sum_{u\in F_{b^{-1}}}\lambda_{u}u$
is the part of $p$ supported on words beginning with $b^{-1}$. Since
$\mathcal{Q}$ is a CRS, we can write $qb^{-1}=\sum_{i=1}^{n}\mu_{i}q_{i}$
for some scalars $\mu_{i}\in K^{\times}$ and elements $q_{i}\in\mathcal{Q}$
with $\text{HTT}\left(q_{i}\right)=b^{-1}$. Then
\[
qp=qg+qb^{-1}\cdot bh=qg+\sum_{i=1}^{n}q_{i}\cdot\mu_{i}bh.
\]
Here, $g=p-h=\sum_{u\notin F_{b^{-1}}}\lambda_{u}u$ is supported
on words not beginning with $b^{-1}$, and since words in $\text{supp}\left(h\right)$
start with $b^{-1}$, words in $\text{supp}\left(bh\right)$ do not
begin with $b$. Thus, the expression is in the desired form.
\item Suppose $f=\sum_{q\in\mathcal{Q}}qg_{q}$ as in part \ref{enu: canonical expression using grobner - existence-1},
with at least one $g_{q}\neq0$. Let $u\in F$ be a word of maximal
length appearing in the support of some $g_{q}$. Among those $q\in\mathcal{Q}$
with $u\in\text{supp}\left(g_{q}\right)$, pick $q_{\text{max}}$
having maximal head term. Let $w=\text{HT}\left(q_{\text{max}}\right)$.
Then $w$ is a neighbor of $T_{\mathcal{Q}}^{\text{HT}},$ and since
there is no cancellation in the product $wu$, we have $wu\notin T_{\mathcal{Q}}^{\text{HT}}$.
The claim will follow by showing that $wu$ is expressed uniquely
as a product in the sum forming $f$, so it cannot cancel. \\
Suppose there exists $q'\in\mathcal{Q}$, $w'\in\text{supp}\left(q'\right)$
and $u'\in\text{supp}\left(g_{q}\right)$ such that $wu=w'u'$. We
will show that $w'=w$, $u'=u$, and $q'=q_{\text{max}}$. Let $\ell=\left|u\right|$.
The vertex $wu$ lies at distance $\ell+1$ from $T_{\mathcal{Q}}^{\text{HT}}.$
Since $\text{supp}\left(q'\right)\subseteq T_{\mathcal{Q}}^{\text{HT}}\cup\partial T_{\mathcal{Q}}^{\text{HT}}$,
the distance of $w'$ from $T_{\mathcal{Q}}^{\text{HT}}$ is at most
$1$. Moreover, $\left|u'\right|\leq\ell$. Therefore, by length considerations,
the equality $wu=w'u'$ implies that $w'\in\partial T_{\mathcal{Q}}^{\text{HT}}$,
that $\left|u'\right|=\ell$, and that there is no cancellation in
$w'u'$. It follows that $u=u'$, as both are he common suffix of
length $\ell$ of $wu=w'u'$, and therefore $w=w'$. Since $u\in\text{supp}\left(g_{q'}\right)$,
we have $\text{HT}\left(q'\right)\preceq\text{HT}\left(q_{\text{max }}\right)=w$.
But $w\in\text{supp}\left(q'\right)$, so $w\preceq\text{HT}\left(q'\right)$.
Thus, $\text{HT}\left(q'\right)=w=\text{HT}\left(q_{\text{max}}\right)$.
By distinctness of head terms, $q'=q_{\text{max}}$.
\item Suppose $\sum_{q\in\mathcal{Q}}qg_{q}=\sum_{q\in\mathcal{Q}}qh_{q}$
for two sets of coefficients $\left(g_{q}\right)_{q\in\mathcal{Q}}$
and $\left(h_{q}\right)_{q\in\mathcal{Q}}$ as in part \ref{enu: canonical expression using grobner - existence-1}.
For each $q\in\mathcal{Q},$ words in $g_{q}$ and $h_{q}$ do not
begin with $\text{\text{HTT}}\left(q\right)^{-1}$, so the same holds
for their difference $g_{q}-h_{q}$. Then $\sum_{q\in\mathcal{Q}}q\left(g_{q}-h_{q}\right)$
is a sum in the form of part \ref{enu: canonical expression using grobner - existence-1}
that equals zero. Since $0$ is supported on $T_{\mathcal{Q}}^{\text{HT}}$,
by part \ref{enu: canonical expression using grobner - if some coefficient is nonzero then not supported on T_I'-1},
$g_{q}-h_{q}=0$, hence $g_{q}=h_{q}$.
\end{enumerate}
\end{proof}
We draw several corollaries from Lemma \ref{lem: Canonical-Representation-using-1}.
We begin with proving Theorem \ref{thm: division with remainder},
which provides a canonical division by a CRS with remainder in a Schreier
transversal $T$.

\begin{proof}[Proof of Theorem \ref{thm: division with remainder}]For
uniqueness, suppose that $f=\sum_{q\in\mathcal{Q}}qg_{q}+r=\sum_{q\in\mathcal{Q}}qg_{q}'+r'$
are two such expressions. Since $\mathcal{Q}\subseteq I$, we have
\[
r-r'=\sum_{q\in\mathcal{Q}}q\left(g_{q}'-g_{q}\right)\in I.
\]
But as both $r$ and $r'$ are supported on $T$ , so is $r-r'$.
Since $T$ is a Schreier transversal for $I$, we conclude that $r-r'=0$,
so $r'=r$. It follows that $\sum_{q\in\mathcal{Q}}q\left(g_{q}'-g_{q}\right)=0$.
This is a representation for the zero element in canonical form with
respect to $\mathcal{Q}$, so by Lemma \ref{lem: Canonical-Representation-using-1},
we must have $g_{q}'-g_{q}=0$ for every $q\in\mathcal{Q}$. Thus,
the two expression for $f$ coincide.

For existence, since $T$ is a Schreier transversal for $I$, the
$I$-coset $f+I$ contains a unique representative $r$ supported
on $T$. The difference $f-r$ lies in $I$. By Lemma \ref{lem: Canonical-Representation-using-1},
there exist coefficients $\left(g_{q}\right)_{q\in\mathcal{Q}}$ such
that $f-r=\sum_{q\in\mathcal{Q}}qg_{q}$, and each $g_{q}$ is supported
on words not beginning with $\text{\text{HTT}}\left(q\right)^{-1}$.
It follows that $f=\sum_{q\in\mathcal{Q}}qg_{q}+r$, as desired.\end{proof}
\begin{cor}
If $\mathcal{Q}$ is a combinatorially reducing system for $I$, then
$\text{Sp}_{K}\left(T_{\mathcal{Q}}^{\text{HT}}\right)\cap I=\left\{ 0\right\} .$
In particular, $I\neq\mathcal{A}$.
\end{cor}
\begin{proof}
Let $f\in\text{Sp}_{K}\left(T_{\mathcal{Q}}^{\text{HT}}\right)\cap I$.
By parts \ref{enu: canonical expression using grobner - existence-1}
and \ref{enu: canonical expression using grobner - if some coefficient is nonzero then not supported on T_I'-1}
of Lemma \ref{lem: Canonical-Representation-using-1}, we conclude
that $f=0$, so $\text{Sp}_{K}\left(T_{\mathcal{Q}}^{\text{HT}}\right)\cap I=\left\{ 0\right\} $.
To deduce the second statement from the first, observe that the unit
$e\in\mathcal{A}$ lies in $\text{Sp}_{K}\left(T_{\mathcal{Q}}^{\text{HT}}\right)$,
since each $q\in\mathcal{Q}$ lies outside of $K$, and thus $\text{HT}\left(q\right)\neq e$.
It follows that $e\notin I$, so $I\neq\mathcal{A}$. 
\end{proof}
\begin{lem}
\label{lem: Properties of the reduction algorithm for a neighbour of T_HT_Q}Let
$\mathcal{Q}$ be a combinatorially reducing system for $I$. If $f$
is supported on $T_{\mathcal{Q}}^{\text{HT}}\cup T_{\mathcal{Q}}^{\text{HT}}b$
for some letter $b\in S\cup S^{-1}$, then there exists an element
$f'\in f+I$ such that:
\begin{enumerate}
\item $f'$ is supported on $T_{\mathcal{Q}}^{\text{HT}}$.
\item $f'$ is obtained from $f$ by adding a $K$-linear combination of
elements from $\left\{ q\in\mathcal{Q}:\text{HTT}\left(q\right)=b\right\} $.
\item $\text{supp}\left(f'\right)\preceq_{\text{max}}\text{supp}\left(f\right)$,
with a strict inequality unless $f'=f$.
\end{enumerate}
\end{lem}
\begin{proof}
We construct a finite sequence $\left(f_{i}\right)$ of elements of
$\mathcal{A}$, each supported on $T_{\mathcal{Q}}^{\text{HT}}\cup T_{\mathcal{Q}}^{\text{HT}}b$,
beginning with $f_{0}:=f$. At the $i$-th step, suppose $f_{i}$
has been defined and is supported on $T_{\mathcal{Q}}^{\text{HT}}\cup T_{\mathcal{Q}}^{\text{HT}}b$.
If $\text{supp}\left(f\right)\subseteq T_{\mathcal{Q}}^{\text{HT}}$,
stop the process. Otherwise, there exists some $v_{i}\in\text{supp}\left(f_{i}\right)$
that is a neighbor of $T_{\mathcal{Q}}^{\text{HT}}$ ending in the
letter $b$. Since $\mathcal{Q}$ is a CRS, there exists $q_{i}\in\mathcal{Q}$
such that $\text{HT}\left(q_{i}\right)=v_{i}$. In particular, $\text{HTT}\left(q_{i}\right)=b$
so by Proposition \ref{prop: element of combinatorially reducing system is supported on T^HT =00005Ccup T^HT(b)},
$q_{i}$ is supported on $T_{\mathcal{Q}}^{\text{HT}}\cup T_{\mathcal{Q}}^{\text{HT}}b$.
Let $\lambda_{i}\in K^{\times}$ be the coefficient of $v_{i}$ in
$f_{i}$, and define $f_{i+1}:=f_{i}-\lambda_{i}q_{i}$. Then $f_{i+1}$
remains supported on $T_{\mathcal{Q}}^{\text{HT}}\cup T_{\mathcal{Q}}^{\text{HT}}b$
as well, since both $f_{i}$ and $q_{i}$ are.

We now claim that the supports of the $f_{i}$ form a strictly $\prec_{\text{max}}$-decreasing
sequence. Indeed, the step from $f_{i}$ to $f_{i+1}$ removes the
term $v_{i}=\text{HT}\left(q_{i}\right)$ from the support, and only
possibly adds strictly smaller terms from $\text{supp}\left(q_{i}\right)$.
Hence, as $\preceq_{\text{max}}$ is a well-order by Claim \ref{claim: The order on finite subsetes is a well-order},
the process terminates after $N$ steps for some integer $N\geq0$,
with $f_{N}$ supported on $T_{\mathcal{Q}}^{\text{HT}}$. Set $f':=f_{N}$.
Then 
\[
f'=f+\sum_{i=0}^{N-1}\lambda_{i}q_{i}\in f+I,
\]
since each $q_{i}\in I$ by assumption.
\end{proof}
We use the Lemma to prove the main theorem of this section -- that
$T_{\mathcal{Q}}^{\text{HT}}$ coincides with the minimal Schreier
transversal of the ideal generated by $\mathcal{Q}$.

\begin{proof}[Proof of Theorem \ref{thm: combinatorially reducing system defines minimal schreier transversal of the ideal which it generates}]Since
$T_{\mathfrak{\mathcal{Q}}}^{\text{HT}}$ is defined by forbidding
certain prefixes, it is clearly prefix-closed. Moreover, by Lemma
\ref{lem: Canonical-Representation-using-1} part \ref{enu: canonical expression using grobner - if some coefficient is nonzero then not supported on T_I'-1},
the $K$-space $\text{Sp}_{K}T_{\mathcal{Q}}^{\text{HT}}$ intersects
$I$ trivially. Thus, $T_{\mathcal{Q}}^{\text{HT}}$ is a partial
Schreier transversal for $I$. By Corollary \ref{cor: T_I  is minimal among partial Schreier transversals for I},
the minimal Schreier transversal $T_{I}$ satisfies $T_{I}\preceq_{\text{min}}T_{\mathcal{Q}}^{\text{HT}}$.

Suppose for contradiction that the inequality is strict: $T_{I}\prec_{\text{min}}T_{\mathcal{Q}}^{\text{HT}}$.
Let $v=\min\left(T_{I}\triangle T_{\mathcal{Q}}^{\text{HT}}\right)$.
By Claim \ref{claim: minimum of symmetric difference is neighbour of larger tree},
$v\in T_{I}$ and is a neighbor of $T_{\mathcal{Q}}^{\text{HT}}$.
Now, by Lemma \ref{lem: Properties of the reduction algorithm for a neighbour of T_HT_Q}
there exists some $f'\in v+I$ supported on $T_{\mathcal{Q}}^{\text{HT}}$
with $\text{supp}\left(f'\right)\prec\left\{ v\right\} $. Every $u\in\text{supp}\left(f'\right)$
then satisfies $u\prec v$ and $u\in T_{\mathcal{Q}}^{\text{HT}}$,
and by minimality of $v$ in the symmetric difference, we must also
have $u\in\text{supp}\left(T_{I}\right)$. Hence, the difference $v-f'$
is a nonzero element of $I$ supported on its Schreier transversal
$T_{I}$, in contradiction.\end{proof}

In light of Theorem \ref{thm: combinatorially reducing system defines minimal schreier transversal of the ideal which it generates},
the element $f'$ obtained in Lemma \ref{lem: Properties of the reduction algorithm for a neighbour of T_HT_Q}
is equal to $\phi_{I}\left(f\right).$ We can therefore reformulate
the Lemma in terms of $I$-\texttt{seconds}.
\begin{cor}
\label{cor: remainder of f modulo T_I is linear combination of Q and has smaller support}Let
$\mathcal{Q}$ be a combinatorially reducing system for $I$, and
let $f\in\mathcal{A}\backslash K$ be supported on $T_{I}$. Suppose
$\text{HTT}\left(f\right)=b$, and let $s$ be the $I$-\texttt{second}
of $f$. Then:
\begin{enumerate}
\item \label{enu: mu*s-f*b^-1 is a K-linear combination}There exists some
$\mu\in K^{\times}$ such that $\mu s-fb^{-1}$ is a $K$-linear combination
of $\left\{ q\in\mathcal{Q}:\text{HTT}\left(q\right)=b^{-1}\right\} $.
\item If $\text{HTT}\left(s\right)=b^{-1}$ then $f$ is also the $I$-\texttt{second}
of $s$, and there exists some $\nu\in K^{\times}$ such that $\nu f-sb$
is a $K$-linear combination of $\left\{ q\in\mathcal{Q}:\text{HTT}\left(q\right)=b\right\} $.
\end{enumerate}
\end{cor}
\begin{proof}
Recall that $s$ is defined as the monic $K^{\times}$-multiple of
$\phi_{I}\left(fb^{-1}\right)$. Let $\mu\in K^{\times}$ satisfy
$\mu s=\phi_{I}\left(fb^{-1}\right)$.
\begin{enumerate}
\item $fb^{-1}$ is supported on $T_{I}b^{-1}$ and $\mathcal{Q}$ is a
CRS for $I$. Lemma \ref{lem: Properties of the reduction algorithm for a neighbour of T_HT_Q}
then implies that $\phi_{I}\left(fb^{-1}\right)$ is obtained from
$fb^{-1}$ by adding a $K$-linear combination of $\left\{ q\in\mathcal{Q}:\text{HTT}\left(q\right)=b^{-1}\right\} $.
Since $\mu s=\phi_{I}\left(fb^{-1}\right),$ the claim follows.
\item Now assume $\text{HTT}\left(s\right)=b^{-1}$. Then $s\notin K$ and
is supported on $T_{I}$, so it admits an $I$-\texttt{second} --
the monic $K^{\times}$-multiple of $\phi_{I}\left(sb\right)$. From
the relation $\mu s=\phi_{I}\left(fb^{-1}\right)$ we deduce: 
\[
\mu s-fb^{-1}\in I\ \ \ \Longrightarrow\ \ \ \mu^{-1}f\in sb+I.
\]
Since $f$ is supported on $T_{I}$, this means that $\phi_{I}\left(sb\right)=\mu^{-1}f$.
The monic $K^{\times}$-multiple of $\phi_{I}\left(sb\right)$ is
$f$. The additional assertion then follows by applying part \ref{enu: mu*s-f*b^-1 is a K-linear combination}
again, this time to $s$ with respect to its $I$-\texttt{second}
$f$.
\end{enumerate}
\end{proof}
\begin{example}
The following simple example shows that the assumption $\text{HTT}\left(s\right)=\text{HTT}\left(f\right)^{-1}$
is necessary in Corollary \ref{cor: remainder of f modulo T_I is linear combination of Q and has smaller support}
for $f$ and $s$ to form an $I$-\texttt{second} pair. Let $I$ be
the zero ideal, so that $T_{I}=F$. Consider $f=x^{6}+x$. Then, regardless
of the order, $\text{HTT}\left(f\right)=x$, and the $I$-\texttt{second}
of $f$ is $s=x^{5}+1$, which satisfies $\text{HTT}\left(s\right)=x$
too. The $I$-\texttt{second} of $s$ is $x^{4}+x^{-1}$, which is
not equal to $f$. 
\end{example}
A generalization of the argument appearing in Lemma \ref{lem: Properties of the reduction algorithm for a neighbour of T_HT_Q}
gives an algorithm for computing $\phi_{I}$ given a CRS for $I$.

\begin{algorithm}
\caption{\textsc{ReduceModuloTI} - Computes $\phi_{I}$ Using a Combinatorially Reducing System for $I$}
\label{alg:reduction_mod_TI}
\begin{algorithmic}[1]
\REQUIRE An element $f \in \mathcal{A}$ and a combinatorially reducing system $\mathcal{Q}$ for a right ideal $I$
\ENSURE $\phi_I(f)$, the unique representative of $f + I$ supported on the minimal Schreier transversal $T_I$
\STATE $r \gets f$
\WHILE{there exists $u \in \mathrm{supp}(r)$ and $q_u \in \mathcal{Q}$ such that $\mathrm{HT}(q_u)$ is a prefix of $u$}
    \STATE Let $\gamma \in K^{\times}$ be the coefficient of $u$ in $r$
    \STATE $u_{\mathrm{suf}} \gets \mathrm{HT}(q_u)^{-1} u$  \COMMENT{so $u = \mathrm{HT}(q_u) \cdot u_{\mathrm{suf}}$ without cancellation}
    \STATE $r \gets r - \gamma q_u u_{\mathrm{suf}}$ \COMMENT{cancel the term $u$ from $r$}
\ENDWHILE
\RETURN $r$
\end{algorithmic}
\end{algorithm}
\begin{thm}
\textbf{\label{thm: Algorithm: reduction modulo T_I given combinatorially reducing system for I}(Reduction
Modulo $T_{I}$ Given a CRS for $I$)}: Let $\mathcal{Q}$ be a combinatorially
reducing system for the right ideal $I\leq\mathcal{A}$, and let $f\in\mathcal{A}$.
Then Algorithm \ref{alg:reduction_mod_TI} terminates in finite time
and computes $\phi_{I}\left(f\right)$, the unique representative
of the $I$-coset $f+I$ supported on the minimal Schreier transversal
$T_{I}$.
\end{thm}
\begin{proof}
\textit{\emph{Since $\mathcal{Q}$ generates $I$, each iteration
subtracts an element of the form $\gamma q_{u}u_{\text{suf}}\in I$.
Therefore, $r\in f+I$ throughout the execution. The algorithm halts
when no word in the support of $r$ has any $\text{HT}\left(q\right)$
for $q\in\mathcal{Q}$ as a prefix. This is precisely the condition
for $r$ to be supported on the subtree $T_{\mathcal{Q}}^{\text{HT}}$,
which equals $T_{I}$ by Theorem \ref{thm: combinatorially reducing system defines minimal schreier transversal of the ideal which it generates}.
Hence, upon halting, $r\in\left(f+I\right)\cap\text{Sp}_{K}T_{I}$,
so $r=\phi_{I}\left(f\right)$.}}

\textit{\emph{It remains to show that the algorithm terminates. Since
$T_{\mathcal{Q}}^{\text{HT}}$ is prefix-closed, every word $v\notin T_{\mathcal{Q}}^{\text{HT}}$
has exactly one prefix in $\partial T_{\mathcal{Q}}^{\text{HT}}$.
Denote this prefix by $v_{\text{exit}}$ and its corresponding suffix
by $v_{\text{suf}}$, so that $v=v_{\text{exit}}v_{\text{suf}}$. }}

\textit{\emph{Consider a single iteration involving some word $u\in\text{supp}\left(r\right)$.
Since $\text{HT}\left(q_{u}\right)$ is a prefix of $u$ and a prefix-neighbor
of $T_{\mathcal{Q}}^{\text{HT}}$, it follows that $u_{\text{exit}}=\text{HT}\left(q_{u}\right)$.}}
The reduction step removes $u$ from $\text{supp}\left(r\right)$
and possibly adds new words from $\text{supp}\left(q_{u}u_{\text{suf}}\right)$.
We claim that any such word $w\neq u$ which is not already inside
$T_{\mathcal{Q}}^{\text{HT}}$ either has a shorter suffix $w_{\text{suf}}$,
or the same suffix $u_{\text{suf}}$ but a strictly smaller prefix
$w_{\text{exit}}\prec u_{\text{exit}}$. Figure \ref{fig: reduction replaces u with closer or smalle with same suffix.}
illustrates this claim. To justify it, let $w\neq u$ be a word in
$\text{supp}\left(q_{u}u_{\text{suf}}\right)$ with $w\notin T_{\mathcal{Q}}^{\text{HT}}$.
Then $w=w'u_{\text{suf}}$ for some $w'\in\text{supp}\left(q_{u}\right)\backslash\left\{ u_{\text{exit}}\right\} $.
Let $b=\text{HTT}\left(q_{u}\right)$. By Proposition \ref{prop: element of combinatorially reducing system is supported on T^HT =00005Ccup T^HT(b)},
we have $w'\in T_{\mathcal{Q}}^{\text{HT}}\cup T_{\mathcal{Q}}^{\text{HT}}b$.
If $w'\in T_{\mathcal{Q}}^{\text{HT}}$ then $w_{\text{suf}}$ is
a proper suffix of $u_{\text{suf}}$, hence shorter. Otherwise, $w'$
is a prefix-neighbor of $T_{\mathcal{Q}}^{\text{HT}}$ ending in $b$,
so the product $w=w'u_{\text{suf}}$ is without cancellation since
$w'$ and $q_{u}$ end in the same letter. It follows that $w_{\text{exit}}=w'$
and $w_{\text{suf}}=u_{\text{suf}}$. Then $w_{\text{exit}}$ is a
word the support of $q_{u}$ which is distinct from its head term
$\text{HT}\left(q_{u}\right)=u_{\text{exit}}$. Hence, $w_{\text{exit}}\prec u_{\text{exit}}$.

Define an auxiliary order $<_{\text{len}}$ on $F$ by declaring shorter
words to be smaller, breaking ties arbitrarily. Now define a lexicographic
order $<_{\text{out}}$ on $F\backslash T_{\mathcal{Q}}^{\text{HT}}$
as follows:
\[
w\prec_{\text{out}}v\text{ if }w_{\text{suf}}<_{\text{len}}v_{\text{suf}}\text{ or }\left(w_{\text{suf}}=v_{\text{suf}}\text{ and }w_{\text{exit}}\prec v_{\text{exit}}\right).
\]
Since $\prec_{\text{out}}$ is a lexicographic order on two well-orders,
it is itself a well-order. The order induced by $\prec_{\text{out}}$
on finite subsets of $F\backslash T_{\mathcal{Q}}^{\text{HT}}$ is
then a well-order too (See Definition \ref{def: order on finite subsets}
and Claim \ref{claim: The order on finite subsetes is a well-order}).
By the previous argument, each reduction strictly decreases $\text{supp}\left(r\right)\backslash T_{\mathcal{Q}}^{\text{HT}}$
with respect to $\prec_{\text{out}}$. Hence, only finitely many steps
may occur\textit{\emph{, and the algorithm terminates.}}

\begin{figure}
\centering{}\label{fig: reduction replaces u with closer or smalle with same suffix.}\includegraphics[scale=0.5]{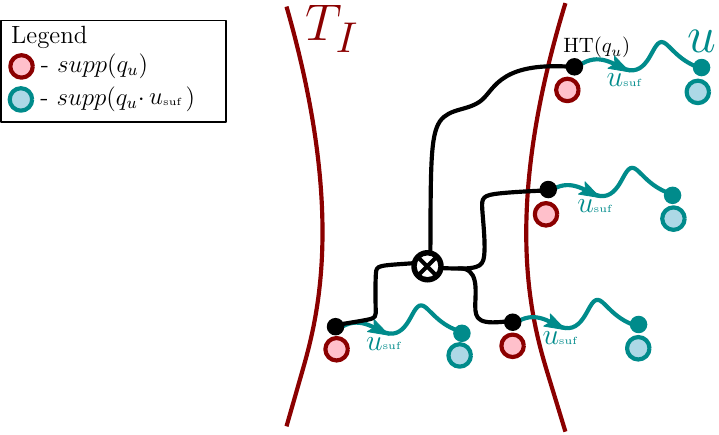}\caption{A single reduction of a word $u\protect\notin T_{\mathcal{Q}}^{\text{HT}}$
using $q_{u}\in\mathcal{Q}$ in Algorithm \ref{alg:reduction_mod_TI}.
The reduction replaces $u$ with words from $\text{supp}\left(f\right)$
with words which, if not closer to $T_{I}$ than $u$, must exit $T_{I}$
through a smaller vertex than $\text{HT}\left(q_{u}\right)$ and have
the same suffix $u_{\text{suf}}$. The respective supports of $q_{u}$
and $q_{u}u_{\text{suf}}$ are indicated by red and cyan circles.}
\end{figure}
\end{proof}
\begin{rem}
\label{rem: cardinality of combinatorially reducing system}Let $\mathcal{Q}$
be a CRS for $I$. To apply Algorithm \ref{alg:reduction_mod_TI},
we must be able to decide, given a word $u\in F$, whether the head
term of some $q_{u}\in\mathcal{Q}$ is a prefix of $u$, and if so,
to retrieve $q_{u}$. The most straightforward case in which this
is possible is when $\mathcal{Q}$ is finite. Note that the finiteness
of $\mathcal{Q}$ is a property depending only on the ideal $I$ which
it generates. Indeed, $\mathcal{Q}$ is in bijection with the set
of its distinct head terms $\left\{ \text{HT}\left(q\right):q\in\mathcal{Q}\right\} $,
which are precisely the prefix-neighbors of $T_{Q}^{\text{HT}}$.
By Theorem \ref{thm: combinatorially reducing system defines minimal schreier transversal of the ideal which it generates},
it follows that $\mathcal{Q}$ is in bijection with the prefix-neighbors
of $T_{I}$. Moreover, we will soon prove that $\left\{ f_{\alpha}\right\} _{\alpha<\omega_{I}}\cup\left\{ s_{\alpha}\right\} _{\alpha<\omega_{I}}$
forms a CRS for $I$ (see Theorem \ref{Thm: properties of second + Grobner basis is combinatorially reducing}),
and that $\left\{ f_{\alpha}\right\} _{\alpha<\omega_{I}}$ is a basis
for $I$ (see Theorem \ref{thm: The exposure basis is a basis}).
We conclude that if $I$ is finitely generated, then every CRS for
$I$ has cardinality equal to $2\text{rk}\left(I\right).$
\end{rem}

\subsection{\label{sec: Combinatorial Properties of the grobner basis}Combinatorial
Properties of the \texttt{Firsts} and \texttt{Seconds}}

A Gröbner basis for a right ideal $I$ is a generating set for $I$
which allows the reduction of elements of $\mathcal{A}$, using some
prespecified set of reduction rules, to a uniquely defined normal
form modulo $I$ (see \cite[Definition~4]{Madlener1993} for a precise
definition). Many Gröbner bases may exist for $I$ -- for example,
one may add more elements of $I$ to an existing one. In this section,
we use the fixed order to distinguish a particular Gröbner basis -
the set of \texttt{firsts} and \texttt{seconds}. This set will be
shown to be a combinatorially reducing system for $I$. The normal
form of an element $f\in\mathcal{A}$ is its remainder $\phi_{I}\left(f\right)$
modulo the minimal Schreier transversal $T_{I}$, and the reduction
of $f$ to $\phi_{I}\left(f\right)$ is done by Algorithm \ref{alg:reduction_mod_TI}.
\begin{defn}
\textbf{\label{def: The Gr=0000F6bner Basis}(The Gröbner Basis)}
The \emph{Gröbner Basis} associated to the proper right ideal $I\lneqq\mathcal{A}$
is the set $B_{I}^{\text{gr}}=B_{I,\prec}^{\text{gr}}=\left\{ f_{\alpha}\right\} _{\alpha<\omega_{I}}\cup\left\{ s_{\alpha}\right\} _{\alpha<\omega_{I}}$.
\end{defn}
We now prove Theorem \ref{thm: introduction Grobner basis is combinatorially reducing },
together with combinatorial properties of the \texttt{seconds}.
\begin{thm}
\label{Thm: properties of second + Grobner basis is combinatorially reducing}Let
$I\lneqq\mathcal{A}$. Then:
\begin{enumerate}
\item \textbf{(Properties of the }\texttt{\textbf{second}}\textbf{ $s_{\alpha}$)}
For every $\alpha<\omega_{I}$:
\begin{enumerate}
\item \label{enu: The second s_alpha is supported on T_I cup T_I * b_alpha^-1}The
\texttt{second} $s_{\alpha}$ is supported on $T_{I}\cup T_{I}b_{\alpha}^{-1}$.
\item \label{enu: HT and HLL of second}Its head term $w_{\alpha}^{-}=\text{HT}\left(s_{\alpha}\right)$
is a prefix-neighbor of $T_{I}$ ending in $b_{\alpha}^{-1}$, with
$w_{\alpha}^{-}\succ w_{\alpha}^{+}$.
\end{enumerate}
\item \textbf{(Theorem \ref{thm: introduction Grobner basis is combinatorially reducing },
restated)} \label{enu: Grobner basis is combinatorially reducing}The
\emph{Gröbner Basis} $B_{I}^{\text{gr}}$ is a combinatorially reducing
system for $I$.
\end{enumerate}
\end{thm}
\begin{proof}
Proceed by transfinite induction. Suppose that for every $\alpha<\omega_{I}$
the theorem holds for $I_{\alpha}$. Denote $\mathcal{Q}_{\alpha}:=B_{I_{\alpha}}^{\text{gr}}$
and $\mathcal{Q}:=B_{I}^{\text{gr}}$.
\begin{enumerate}
\item Let $\alpha<\omega_{I}$. By the induction hypothesis, $\mathcal{Q_{\alpha}}$
is a CRS for $I_{\alpha}$. We proceed by another transfinite induction.
This allows us to suppose further for every $\beta<\alpha$, that
the \texttt{second} $s_{\beta}$ is supported on $T_{I}\cup T_{I}b_{\beta}^{-1}$
and that $w_{\beta}^{-}$ ends in $b_{\beta}^{-1}$.
\begin{enumerate}
\item Recall that $s_{\alpha}$ is, by definition, the $I_{\alpha}$-\texttt{second}
of $f_{\alpha}$. Since $Q_{\alpha}$ is a CRS for $I_{\alpha}$,
it follows from Corollary \ref{cor: remainder of f modulo T_I is linear combination of Q and has smaller support}
that $s_{\alpha}$ is a $K$-linear combination of $f_{\alpha}b_{\alpha}^{-1}$
and other $q\in\mathcal{Q}_{\alpha}$ with $\text{HTT}\left(q\right)=b_{\alpha}^{-1}$.
It therefore suffices to show that both $f_{\alpha}b_{\alpha}^{-1}$
and any such $q$ are supported on $T_{I}\cup T_{I}b_{\alpha}^{-1}$.\\
By Proposition \ref{prop: Form for first using T_I}, we have $\text{supp}\left(f_{\alpha}\right)\subseteq T_{I}\cup T_{I}b_{\alpha}$.
Multiplying from the right by $b_{\alpha}^{-1}$ gives $\text{supp}\left(f_{\alpha}b_{\alpha}^{-1}\right)\subseteq T_{I}\cup T_{I}b_{\alpha}^{-1}$.
Now let $q\in\mathcal{Q}_{\alpha}$ with $\text{HTT}\left(q\right)=b_{\alpha}^{-1}$.
If $q=f_{\beta}$ for some $\beta<\alpha$, then by the same Proposition,
$\text{supp}\left(q\right)\subseteq T_{I}\cup T_{I}b_{\alpha}^{-1}$.
If, instead, $q=s_{\beta}$ for some $\beta<\alpha$, then $\text{supp}\left(s_{\beta}\right)\subseteq T_{I}\cup T_{I}b_{\beta}^{-1}$.
But since $w_{\beta}^{-}$ ends in $b_{\beta}^{-1}$, we have $b_{\alpha}^{-1}=\text{HTT}\left(s_{\beta}\right)=b_{\beta}^{-1}$,
so $\text{supp}\left(s_{\beta}\right)\subseteq T_{I}\cup T_{I}b_{\alpha}^{-1}$,
as required.
\item By the first part, $s_{\alpha}$ is supported on $T_{I}\cup T_{I}b_{\alpha}^{-1}$.
Proposition \ref{prop: HT(f) supported on neighbour of T_I} implies
that $\text{HT}\left(s_{\alpha}\right)=w_{\alpha}^{-}$ is a prefix-neighbor
of $T_{I}$. Moreover, as $w_{\alpha}^{-}\in\text{supp}\left(s_{\alpha}\right)\subseteq T_{I}\cup T_{I}b_{\alpha}^{-1}$,
its last letter must be $b_{\alpha}^{-1}$. Since $w_{\alpha}^{+}$
ends in $b_{\alpha}$, we have $w_{\alpha}^{-}\neq w_{\alpha}^{+}$
. Finally, because $s_{\alpha}\in I$ is monic and supported on $T_{I_{\alpha}}$
by definition, and $f_{\alpha}$ has $\prec_{\text{max}}$-minimal
support among such elements, it follows that $w_{\alpha}^{-}\succ w_{\alpha}^{+}$.
\end{enumerate}
\item The \texttt{firsts} and \texttt{seconds} of $I$ are monic by definition,
and $\mathcal{Q}\subseteq I\subseteq\mathcal{A}\backslash K$ since
$I$ is a proper ideal.\\
First, we show that the head terms $\left(\text{HT}\left(q\right)\right)_{q\in\mathcal{Q}}$
are distinct. Within the same ordinal $\alpha<\omega_{I}$, part \ref{enu: HT and HLL of second}
implies that $w_{\alpha}^{-}\succ w_{\alpha}^{+}$. For distinct ordinals
$\beta<\alpha<\omega_{I}$, both $w_{\alpha}^{+}$ and $w_{\alpha}^{-}$
are supported on $T_{I_{\alpha}}$ by the respective definitions of
$f_{\alpha}$ and $s_{\alpha}$. In contrast, $w_{\beta}^{+}$ and
$w_{\beta}^{-}$ are prefix-neighbors of $T_{I_{\alpha}}$ by part
\ref{enu: HT and HLL of second} and Proposition \ref{prop: Form for first using T_I},
both applied to the exposure process of $I_{\alpha}$. Hence, all
elements in the set $\left(w_{\alpha}^{\pm}\right)_{\alpha<\omega_{I}}$
are distinct.\\
Second, let $q\in\mathcal{Q}$ and let $b=\text{HTT}\left(q\right)$.
Let $\alpha<\omega_{I}$ be the ordinal such that $q\in\left\{ f_{\alpha},s_{\alpha}\right\} $.
The $I_{\alpha}$-\texttt{second} of $f_{\alpha}$ is $s_{\alpha}$,
which satisfies by the first part $\text{HTT}\left(s_{\alpha}\right)=b_{\alpha}^{-1}=\text{HTT}\left(f_{\alpha}\right)^{-1}$.
Applying Corollary \ref{cor: remainder of f modulo T_I is linear combination of Q and has smaller support}
using the combinatorially reducing system $\mathcal{Q}_{\alpha}$
for $I_{\alpha}$, we conclude that $f_{\alpha}$ is also the $I_{\alpha}$-\texttt{second}
of $s_{\alpha}$. Hence, denoting by $p$ the $I_{\alpha}$-\texttt{second}
of $q$, we have 
\[
p=\begin{cases}
s_{\alpha} & \text{if }q=f_{\alpha},\\
f_{\alpha} & \text{if }q=s_{\alpha}.
\end{cases}
\]
In both cases, $p\in\mathcal{Q}$ and ends in $b^{-1}$. Applying
Corollary \ref{cor: remainder of f modulo T_I is linear combination of Q and has smaller support}
again to $q$, we have that $qb^{-1}$ is a $K$-linear combination
of $p$ and elements $\left\{ q'\in\mathcal{Q}_{\alpha}:\text{HTT}\left(q'\right)=b^{-1}\right\} $.
Since $\mathcal{Q}_{\alpha}\subseteq\mathcal{Q}$, this is as required.
\\
Third and last, for every $q\in\mathcal{Q}$, its head term lies outside
of $T_{I}$ by part \ref{enu: HT and HLL of second} and Proposition
\ref{prop: Form for first using T_I}. Since $T_{I}$ is prefix-closed,
we deduce that $T_{I}\subseteq T_{\mathcal{Q}}^{\text{HT}}$. Let
$q\in\mathcal{Q}$ and $b=\text{HTT}\left(q\right)$. We conclude,
using Part \ref{enu: The second s_alpha is supported on T_I cup T_I * b_alpha^-1}
and Proposition \ref{prop: Form for first using T_I}, that
\[
\text{supp}\left(q\right)\subseteq T_{I}\cup T_{I}b\subseteq T_{\mathcal{Q}}^{\text{HT}}\cup T_{\mathcal{Q}}^{\text{HT}}b.
\]
\end{enumerate}
\end{proof}
\begin{rem}
\label{rem: why we define seconds this way}The support of $s_{\alpha}$
may contain other prefix-neighbors of $T_{I}$ in addition to $w_{\alpha}^{-}$.
One might therefore consider replacing $s_{\alpha}$ with a ``cleaner''
version $s_{\alpha}':=w_{\alpha}^{-}-\phi_{I}\left(w_{\alpha}^{-}\right)$.
The resulting CRS $\mathcal{Q}':=\left\{ f_{\alpha}:\alpha<\omega_{I}\right\} \cup\left\{ s_{\alpha}':\alpha<\omega_{I}\right\} $
for $I$ has the appealing property that the head term of each $q'\in\mathcal{Q}'$
is a prefix-neighbor of $T_{I}$, while the rest of its support lies
entirely within $T_{I}$. However, we define $s_{\alpha}$ as we do
because this choice remains stable throughout the exposure process:
$s_{\alpha}$ does not require fixing at later stages. At each stage
$\alpha<\omega_{I}$, precisely two new elements -- $f_{\alpha}$
and $s_{\alpha}$ -- are added to the existing Gröbner basis, maintaining
$B_{I_{\alpha}}^{\text{gr}}\subseteq B_{I}^{\text{gr}}$. 

Another approach was developed by Reinert \cite[Chapter~5]{Reinert1995},
building on earlier work with Madlener \cite{Madlener1993}. The algorithm
of \cite{Madlener1993} computes Gröbner bases for general monoid
rings, but does not specialize well to free group algebras (as noted
by Rosenmann \cite{Rosenmann1993}). Reinert's later specialization
addresses this by outputting an ``efficient'' Gröbner basis for
$I$, which we denote by $\mathcal{Q}''$, consisting of $\text{rk}I$
pairs of the form $\left(f,s\right)$ with $s=f\cdot\text{HLL}\left(f\right)^{-1}$.
The system $\mathcal{Q}''$ reflects a simpler relation between pairs,
but does not include the exposure basis $B_{I}$ (free bases are not
discussed in this work). Our construction, by contrast, ensures that
$B_{I}^{\text{gr}}$ contains $B_{I}$ throughout.
\end{rem}
The following Example illustrates Remark \ref{rem: why we define seconds this way}.
\begin{example}
\label{exa: Another option for Grobner basis}Recall Example \ref{examples of seconds}
(\ref{enu: involved example of seconds}), in which $K=\mathbb{F}_{2}$,
and $\prec$ is the shortlex order on $F=\left\langle x,y\right\rangle $
satisfying $y^{-1}\prec x^{-1}\prec x\prec y$. The exposure process
for $I=\left(y^{-2}+y+x\right)\mathcal{A}+\left(xy^{-1}+y\right)\mathcal{A}$
consists of two stages, yielding the Gröbner basis
\begin{align*}
f_{0} & =y^{-2}+y+x,\ s_{0}=y^{2}+xy+y^{-1}\\
f_{1} & =xy^{-1}+y,\ \ \ \ s_{1}=xy+x+y^{-1}.
\end{align*}
Observe that $\text{HT}\left(s_{1}\right)=xy$ also appears in $\text{supp}\left(s_{0}\right)$.
Thus, $s_{0}$ has two terms in its support outside $T_{I}$, namely
$y^{2}$ and $xy$. Since $y^{2}+x=s_{0}+s_{1}\in I$, it follows
that $\phi_{I}\left(y^{2}\right)=x$. The proposed change to $\mathcal{Q}'$
would then replace $s_{0}=y^{2}+xy+y^{-1}$ with 
\[
s_{0}'=y^{2}-\phi_{I}\left(y^{2}\right)=y^{2}+x=s_{0}+s_{1}.
\]
See Figure \ref{fig: different option for combinatorially reducing system}
for a visual description of these changes. While $s_{0}'$ is simpler
-- in that only its head term $y^{2}$ lies outside of $T_{I}$ --
this formulation complicates the inductive structure of the exposure
process: at stage $\alpha=1$, the previously exposed $s_{0}$ would
have to be retroactively fixed by adding the newly exposed $s_{1}$.\\
For comparison, Reinert's algorithm \cite{Reinert1995} outputs the
CRS $\mathcal{Q}''=\left\{ \left(f_{0}'',s_{0}''\right),\left(f_{1}'',s_{1}''\right)\right\} $,
where
\begin{align*}
f_{0}'' & =y^{-2}+y+x,\ s_{0}''=y^{2}+xy+y^{-1},\\
f_{1}'' & =xy^{-1}+y^{-2}+x,\ \ \ \ s_{1}''=xy+x+y^{-1}.
\end{align*}
Note that $s_{i}''=f_{i}''\cdot\text{HLL}\left(f_{i}''\right)^{-1}$,
illustrating the simpler relation between pairs. However, $f_{1}''$
is not an exposure basis element, and $\mathcal{Q}''$ does not contain
$B_{I}$.

\begin{figure}
\centering{}\label{fig: different option for combinatorially reducing system}\includegraphics[scale=0.5]{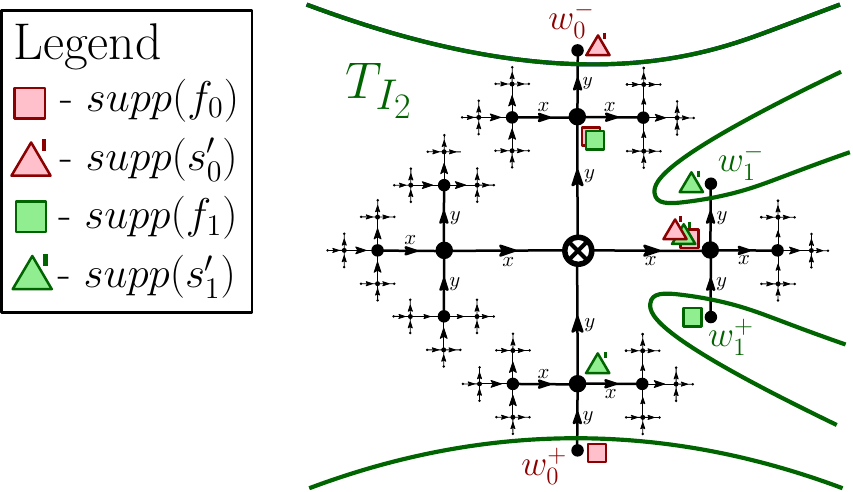}\caption{Illustration of Example \ref{exa: Another option for Grobner basis}.
The combinatorially reducing system $\mathcal{Q}'=\left\{ f_{0},s_{0}',f_{1},s_{1}'\right\} $
is shown. While $s_{1}'$ coincides with the original \texttt{second}
$s_{1}$, the element $s_{0}'$ differs from the \texttt{second} $s_{0}$
and no longer contains $\text{HT}\left(s_{1}'\right)=w_{1}^{-}$ in
its support (cf. Figure \ref{fig: involved example of seconds}).
The supports of $s_{0}'$ and $s_{1}'$ are indicated by red and green
primed triangles, respectively.}
\end{figure}
\end{example}
\begin{rem}
It is somewhat unclear to us which of the two variants of a CRS described
in Remark \ref{rem: why we define seconds this way} -- $\mathcal{Q}'$
or $\mathcal{Q}=B_{I}^{\text{gr}}$ -- was intended to serve as the
output Gröbner basis in the algorithm designed by Rosenmann \cite{Rosenmann1993}
for the shortlex order. On the one hand, the proof of correctness
indicates that each $\alpha\in GBasis$ has only its head term outside
$T_{I}$, which suggests that $\mathcal{Q}'$ may have been intended.
On the other hand, it is unclear whether the algorithm retroactively
updates earlier \texttt{seconds} as required for $\mathcal{Q}$.
\end{rem}

\subsection{\label{sec: The Exposure Basis is a Basis}The Exposure Basis is
a Basis}

In this section we finally prove Theorem \ref{thm: The exposure basis is a basis}
-- showing that the exposure basis $B_{I}=\left\{ f_{\alpha}:\alpha<\omega_{I}\right\} $
forms a basis for $I$. Recall that for every $w\in F$, the set $F_{w}$
consists of all words in $F$ beginning with $w$.
\begin{lem}
\label{lem: top coefficients vanish for exposure basis}Let $I\lneqq\mathcal{A}$
and let $\alpha<\omega_{I}$. If $f_{\alpha}g+s_{\alpha}h\in I_{\alpha}$
for some $g\in$$\mathcal{A}$ supported on $F\backslash F_{b_{\alpha}^{-1}}$
and $h\in\mathcal{A}$ supported on $F\backslash F_{b_{\alpha}}$
then $g=h=0$.
\end{lem}
\begin{proof}
By Theorem \ref{Thm: properties of second + Grobner basis is combinatorially reducing},
$B_{I_{\alpha}}^{\text{gr}}$ is a CRS for $I_{\alpha}$. Thus, we
can write canonically $f_{\alpha}g+s_{\alpha}h=\sum_{q\in B_{I_{\alpha}}^{\text{gr}}}qg_{q}$
as in Lemma \ref{lem: Canonical-Representation-using-1}. Then the
difference $f_{\alpha}g+s_{\alpha}h-\sum_{q\in B_{I_{\alpha}}^{\text{gr}}}qg_{q}$
vanishes, and is in canonical form with respect to the larger CRS
$B_{I}^{\text{gr}}$ (recall that $f_{\alpha},g_{\alpha}\notin B_{I_{\alpha}}^{\text{gr}}$).
By Lemma \ref{lem: Canonical-Representation-using-1} (\ref{enu: canonical expression using grobner - coefficients are unique-1}),
we conclude that $g=h=0$.
\end{proof}
We now prove the main Theorem of this section, Theorem \ref{thm: The exposure basis is a basis},
to finally assert that the Exposure basis $B_{I}$ is indeed a basis
for $I$.

\begin{proof}[Proof of Theorem \ref{thm: The exposure basis is a basis}]If
$I=\mathcal{A}$, the exposure basis consists of a single element
-- the unity $1\in\mathcal{A}$ -- which is trivially a basis for
$\mathcal{A}$. Assume $I\lneqq\mathcal{A}$. By the definition of
$\omega_{I}$, the set $\left\{ f_{\alpha}:\alpha<\omega_{I}\right\} $
generates $I$ . To prove the generation is free, we show that for
every $\alpha<\omega_{I}$, the intersection of $f_{\alpha}\mathcal{A}$
and $I_{\alpha}$ is trivial. Suppose $f_{\alpha}p\in I_{\alpha}$
for some $p\in\mathcal{A}$. By the definition of $s_{\alpha}$, there
exists $\mu\in K^{\times}$ such that $\mu s_{\alpha}-f_{\alpha}b_{\alpha}^{-1}\in I_{\alpha}$.
Write $p=g+h$, where $h$ is supported on $F_{b_{\alpha}^{-1}}$
and $g$ on $F\backslash F_{b_{\alpha}^{-1}}$. Then 
\[
f_{\alpha}g+s_{\alpha}\mu b_{\alpha}h=f_{\alpha}p+\left(\mu s_{\alpha}-f_{\alpha}b_{\alpha}^{-1}\right)b_{\alpha}h\in I_{\alpha}.
\]
By Lemma \ref{lem: top coefficients vanish for exposure basis}, since
$g$ is supported on $F\backslash F_{b_{\alpha}^{-1}}$ and $\mu b_{\alpha}h$
on $b_{\alpha}F_{b_{\alpha}^{-1}}=F\backslash F_{b_{\alpha}}$, both
must vanish. Hence, $p=g+h=0$.\end{proof}
\begin{cor}
\label{cor: Seconds are also a basis.}Let $I\lneqq\mathcal{A}$.
The \texttt{seconds} $\left(s_{\alpha}\right)_{\alpha<\omega_{I}}$
also form a basis for $I$.
\end{cor}
\begin{proof}
First, we prove that $\left(s_{\alpha}\right)_{\alpha<\omega_{I}}$
generate $I$. We proceed by transfinite induction. Suppose for each
$\beta<\omega_{I}$ that the \texttt{seconds} $\left(s_{\alpha}\right)_{\alpha<\beta}$
form a basis for $I_{\beta}$. Let $I'$ be the ideal generated by
$\left(s_{\alpha}\right)_{\alpha<\omega_{I}}$. Since $s_{\alpha}\in I$
for all $\alpha<\omega_{I}$, we have $I'\subseteq I$. Conversely,
let $\beta<\omega_{I}$ and let $\mu\in K^{\times}$ satisfy $f_{\beta}b_{\beta}^{-1}-\mu s_{\beta}\in I_{\beta}$.
By the induction hypothesis, $I_{\beta}$ is generated by $\left\{ s_{\alpha}:\alpha<\beta\right\} $,
so $I_{\beta}\subseteq I'$. Therefore, $f_{\beta}\in s_{\beta}\mathcal{A}+I_{\beta}\subseteq I'$.
Hence $I'$ contains the generating set $\left(f_{\beta}\right)_{\beta<\omega_{I}}$
for $I$, so it contains $I$. The generation is free by the same
argument as in Theorem \ref{thm: The exposure basis is a basis},
reversing the roles of $b_{\alpha}$ and $b_{\alpha}^{-1}$. 
\end{proof}
We finish this section by showing a result which illustrates the flexibility
of exposure bases. Let $T$ be a Schreier transversal for a right
ideal $I\leq\mathcal{A}$. Recall that Lewin \cite[Theorem~1]{Lewin1969}
used $T$ to construct a basis for $I$. We denote this basis by $B_{T}^{+}$
and call it the \emph{Lewin basis induced} by $T$. Let $\phi_{I,T}:\mathcal{A}\rightarrow\text{Sp}_{K}\left(T\right)$
be the transversal function associated to $T$ with respect to $I$.
The Lewin basis is given by $B_{T}^{+}=\left\{ v-\phi_{I,T}\left(v\right):v\in\partial_{+}T\right\} $,
where $\partial_{+}T$ is the set of prefix-neighbors of $T$ which
are either the empty word or whose last letter lies in $S$ (rather
than $S^{-1}$). We show that $B_{T}^{+}$ can be realized as an exposure
basis, with respect to a suitably chosen exposure order.
\begin{prop}
\label{prop: Lewin bases are exposure bases}Let $T$ be a Schreier
transversal for a right ideal $I\leq\mathcal{A}$. There exists an
exposure order $\prec_{T}^{+}$ on $F$ such that $B_{I,\prec_{T}^{+}}=B_{T}^{+}$.
\end{prop}
\begin{proof}
Let $T'=T\cup\partial_{+}T$. Since $T'$ is obtained by adjoining
to the prefix-closed set $T$ some of its prefix-neighbors, $T'$
is itself prefix-closed. Let $\prec$ denote a shortlex order. By
applying Claim \ref{claim: Sum of two order on tree and complement is an exposure order too}
twice, with respect to $T'$ and then $T$, we obtain an exposure
order $\prec_{T}^{+}:=\left(\prec_{T'}\right)_{T}$ which satisfies
$u\prec_{T}^{+}v\prec_{T}^{+}w$ for every $u\in T,$ $v\in\partial_{+}T$,
and $w\in F\backslash T'$. By Proposition \ref{prop: Every Schreier tranversal is minimal with respect to some exposure order},
we have $T_{I,\prec_{T}^{+}}=T$. In particular, the associated transversal
functions coincide: $\phi_{I,\prec_{T}^{+}}=\phi_{I,T}$.

Now, both $B_{T}^{+}$ and $B_{I,\prec_{T}^{+}}$ are bases of the
same ideal $I$: the former by Lewin \cite[Theorem~1]{Lewin1969},
and the latter by Theorem \ref{thm: The exposure basis is a basis}.
Hence, it suffices to prove the inclusion $B_{I,\prec_{T}^{+}}\subseteq B_{T}^{+}$.
Let $f\in B_{I,\prec_{T}^{+}}$. By Proposition \ref{prop: Form for first using T_I},
we can write 
\[
f=w^{+}-\phi_{I,\prec_{T}^{+}}\left(w^{+}\right)=w^{+}-\phi_{I,T}\left(w^{+}\right),
\]
where $w^{+}=\text{HT}_{\prec_{T}^{+}}\left(f\right)$ is a prefix-neighbor
of $T_{I,\prec_{T}^{+}}=T$. It remains to prove that $w^{+}\in\partial_{+}T$.
Suppose towards contradiction that $w^{+}$ ends in $b^{-1}$ for
some $b\in S$. Then in particular, $f\notin K$, so $I\neq\mathcal{A}$.
Let $s$ be the \texttt{second} of $f$ in the exposure process of
$I$, and set $w^{-}:=\text{HT}_{\prec_{T}^{+}}\left(s\right)$. By
Theorem \ref{Thm: properties of second + Grobner basis is combinatorially reducing},
$w^{-}$ is a prefix-neighbor of $T$ ending in $\text{HLL}_{\prec_{T}^{+}}\left(f\right)^{-1}=b$,
and moreover $w^{+}\prec_{T}^{+}w^{-}$. This contradicts the construction
of $\prec_{T}^{+}$, since $w^{+}\in F\backslash T'$ while $w^{-}\in\partial_{+}T$.
Hence $w^{+}\in\partial_{+}T$, as required.
\end{proof}

\subsection{\label{sec: Suffix Preserving Orders}Suffix-Invariant Orders}

The following section departs briefly from the main thread to examine
a narrower class of exposure orders -- those satisfying an additional
condition we call suffix-invariance. For such orders, the key objects
we previously associated with an ideal $I\leq\mathcal{A}$ -- the
minimal Schreier transversal $T_{I}$, the transversal function $\phi_{I}$
and the exposure basis $B_{I}$ -- admit simpler definitions and
improved properties (see Theorem \ref{thm: suffix preserving order main theorem}).
In particular, applying $\phi_{I}$ to any $f\notin\text{Sp}_{K}\left(T_{I}\right)$
strictly $\prec_{\text{max}}$-decreases its support, facilitating
easier control over algorithm termination (cf.\ Example \ref{exa: word smaller than its phi}).
These properties illustrate the advantages of using the shortlex order,
for example, over a general exposure order -- an approach employed
by Rosenmann (\cite{Rosenmann1993}) and others in the context of
Gröbner basis theory.
\begin{defn}
\label{def: suffix preserving}A well-order $\prec$ on $F$ is called
\emph{suffix-invariant} if whenever $u,u',v\in F$ satisfy $u\prec u'$
and there is no cancellation in the product $u'v$, then $uv\prec u'v$.
\end{defn}
To clarify, a suffix-order is required to be a well-order.
\begin{example}
The following examples illustrate Definition \ref{def: suffix preserving}.
\begin{enumerate}
\item \textbf{(Non-example)} Consider the integers $\mathbb{Z}$, viewed
as the free group on the single generator $1$. The usual order on
$\mathbb{Z}$ is not a well-order, since infinite strictly decreasing
sequences exist. However, for every $m,m',n\in\mathbb{Z}$ such that
$m<m'$ we have $m+n<m'+n$ (regardless of cancellations). This example
demonstrates that the explicit requirement of well-ordering in Definition
\ref{def: suffix preserving} is not redundant. 
\item \textbf{(Shortlex)} The shortlex order $\prec$ from Example \ref{exa: Examples of Exposure orders}
(\ref{enu: shortlex}) is suffix-invariant.
\item \textbf{(Weighted Shortlex)} Let $\prec_{\text{slex}}$ be the shortlex
order on $F=\left\langle x,y\right\rangle $ from Example \ref{exa: Examples of Exposure orders}
(\ref{enu: shortlex}). We describe a method for constructing more
elaborate suffix-invariant orders from $\prec_{\text{slex}}$ using
weights. Let $\left(G,<_{G}\right)$ be a bi-ordered group (i.e.,
such that $gh<_{G}g'h$ and $hg<_{G}hg'$ whenever $g<_{G}g'$). Denote
by $G_{+}=\left\{ g\in G:g>_{G}1_{G}\right\} $ positive cone. Suppose
that the finitely generated sub-semigroups of $G_{+}$ are well-ordered
(with respect to $<_{G}$). Fix a function $f:\left\{ x,y,x^{-1},y^{-1}\right\} \rightarrow G_{+}$,
and define the $f$-weight of a reduced word $w=w_{1}w_{2}...w_{\ell}$,
where $w_{i}\in\left\{ x,y,x^{-1},y^{-1}\right\} $, by $f\left(w\right)=f\left(w_{1}\right)f\left(w_{2}\right)...f\left(w_{\ell}\right)$.
The restriction of $<_{G}$ to $f\left(F\right)$ is a well-order
since it is finitely generated as a semigroup by the $f$-images of
$\left\{ x,y,x^{-1},y^{-1}\right\} $. Define an order $\prec_{f}$
on $F$ by 
\[
u\prec_{f}v\ \ \ \iff\ \ \ \left(f\left(u\right),u\right)<_{\text{lex}}\left(f\left(v\right),v\right),
\]
where the lexicographic order $<_{\text{lex}}$ uses $<_{G}$ on $f$-weights
and $\prec_{\text{slex}}$ to break ties. The resulting order $\prec_{f}$
is a well-order since both $\left(f\left(F\right),<_{G}\right)$ and
$\left(F,\prec_{\text{slex}}\right)$ are well-ordered. \\
To verify suffix-invariance, first note that any subword $u'$ of
a reduced word $u\in F$ satisfies $f\left(u'\right)\leq_{G}f\left(u\right)$,
with equality if and only if $u'=u$. Indeed, since $f$ takes values
in $G_{+}$, writing $u=vu'w$ without cancellation, we have 
\[
f\left(u\right)=f\left(v\right)f\left(u'\right)f\left(w\right)\geq_{G}1_{G}\cdot f\left(u'\right)\cdot1_{G}=f\left(u'\right),
\]
where equality happens if and only if $v=w=e$. Now let $u,u',v\in F$
such that $u\prec_{f}u'$ and there is no cancellation in $u'v$.
By definition of $\prec_{f}$, we have $f\left(u\right)\leq_{G}f\left(u'\right)$,
so the $f$-weights of $uv$ and $u'v$ satisfy: 
\[
f\left(uv\right)\leq_{G}f\left(u\right)f\left(v\right)\leq_{G}f\left(u'\right)f\left(v\right)=f\left(u'v\right).
\]
If these weights are unequal, it follows that $uv\prec_{f}u'v$. If,
instead, $f\left(uv\right)=f\left(u'v\right)$, then both inequalities
above must be equalities, implying $f\left(u\right)=f\left(u'\right)$.
Since $u\prec_{f}u'$, this means that $u\prec_{\text{slex}}u'$.
But $\prec_{\text{slex}}$ is suffix-invariant, and $u'v$ is without
cancellation, so $uv\prec_{\text{slex}}u'v$, which implies $uv\prec_{f}u'v$
as well. Thus, $\prec_{f}$ is suffix-invariant.\\
We give several examples of this construction, all assuming $y^{-1}\prec_{\text{slex}}x^{-1}\prec_{\text{slex}}x\prec_{\text{slex}}y$:
\begin{enumerate}
\item Let $G=\mathbb{R}$ with the usual order, and suppose that $f$ is
constant on $\left\{ x,y,x^{-1},y^{-1}\right\} $. Then the order
$\prec_{f}$ on $F$ coincides with $\prec_{\text{slex}}$, showing
that the construction generalizes the standard shortlex order.
\item Let $G=\mathbb{R}$ again, and define $f\left(x\right)=1$ and $f\left(x^{-1}\right)=f\left(y\right)=f\left(y^{-1}\right)=4$.
The smallest words in the resulting order $\prec_{f}$ are 
\[
e\prec_{f}x\prec_{f}x^{2}\prec_{f}x^{3}\prec_{f}y^{-1}\prec_{f}x^{-1}\prec_{f}y\prec_{f}x^{4}\prec_{f}y^{-1}x\prec_{f}...
\]
\item Let $G=\mathbb{Z}^{2}$, ordered lexicographically with $\left(1,0\right)>_{G}\left(0,1\right)$,
and define 
\[
f\left(x\right)=f\left(x^{-1}\right)=\left(1,0\right),\ f\left(y\right)=f\left(y^{-1}\right)=\left(0,1\right).
\]
This defines another suffix-invariant order $\prec_{f}$ on $F$,
in which every nonzero power of $y$ is $\prec_{f}$-smaller than
any nonzero power of $x$.
\end{enumerate}
\end{enumerate}
\end{example}
\begin{claim}
\label{claim: Suffix-invariant orders are exposure orders}Every suffix-invariant
order on $F$ is an exposure order.
\end{claim}
\begin{proof}
It suffices to show that for every word $v\in F$ and letter $b\in S\cup S^{-1}$,
if there is no cancellation in the product $vb$, then $v\prec vb$.
Suppose otherwise, that $v\succ vb$. Then by suffix-invariance, $vb\succ vb^{2}$.
Repeating this argument yields the infinite strictly decreasing chain
$v\succ vb\succ vb^{2}\succ...$, contradicting well-ordering.
\end{proof}
From here until the end of the section, fix a right-ideal $I\leq\mathcal{A}$.
For every order $\prec$ on $F$, denote by $T_{I,\prec}'$ the set
of words $w\in F$ which have $\prec_{\text{max}}$-minimal support
in their $I$-coset.
\begin{claim}
\label{cla: T_I_prime is a partial schreier transversal}For every
total order $\prec$ on $F$, we have $I\cap\text{Sp}_{K}\left(T_{I,\prec}'\right)=\left\{ 0\right\} $.
\end{claim}
\begin{proof}
Suppose for contradiction there exists a nonzero $f\in I\cap\text{Sp}_{K}\left(T_{I,\prec}'\right)$.
Without loss of generality, assume $f$ is monic, and let $w=\text{HT}\left(f\right)$.
Then $f':=w-f\in w+I$ , but $f'$ has support strictly $\prec_{\text{max}}$-smaller
than $w$, contradicting the assumption that $w\in T_{I,\prec}'$.
\end{proof}
\begin{claim}
\label{claim: word in the support of smallest in its I-coset is also smallest in its I-coset}Let
$\prec$ be a total order on $F$ and let $f\in\mathcal{A}$ have
$\prec_{\text{max}}$-minimal support among elements of its $I$-coset
$f+I$. Then $\text{supp}\left(f\right)\subseteq T_{I,\prec}'$.
\end{claim}
\begin{proof}
Suppose not. Then there exists $w\in\text{supp}\left(f\right)$ and
$g\in w+I$ such that $\text{supp}\left(g\right)\prec_{\text{max}}\left\{ w\right\} $.
Let $\lambda\in K^{\times}$ be the coefficient of $w$ in $f$. Then
$f-\lambda w+\lambda g\in f+I$ has smaller support than $f$, in
contradiction.
\end{proof}
Assuming the requirement on suffixes in Definition \ref{def: suffix preserving},
we have improved properties for $T_{I,\prec}'$.
\begin{claim}
\label{cla: set of smallest words in their coset is a prefix closed}Let
$\prec$ be a suffix-invariant order on $F$. The set $T_{I,\prec}'$
is prefix-closed.
\end{claim}
\begin{proof}
Let $w_{1}w_{2}\in T_{I}'$, and suppose there is no cancellation
in the product. Suppose towards contradiction that $w_{1}\notin T_{I,\prec}'$.
Then there exists $f\in w_{1}+I$ with $\text{supp}\left(f\right)\prec_{\text{max}}\left\{ w_{1}\right\} $.
Each word $u\in\text{supp}\left(f\right)$ satisfies $u\prec w_{1}$,
so by suffix-invariance $uw_{2}\prec w_{1}w_{2}$. Therefore, $\text{supp}\left(fw_{2}\right)\prec_{\max}\left\{ w_{1}w_{2}\right\} $.
But $fw_{2}\in w_{1}w_{2}+I$, contradicting $w_{1}w_{2}\in T_{I,<}'$.
\end{proof}
The following example shows that Lemma \ref{cla: set of smallest words in their coset is a prefix closed}
might fail if the order is not suffix-invariant.
\begin{example}
\label{exa: T_I' not closed under prefixs}Consider the ideal $I=\left(x-1\right)\mathcal{A}$
from Example \ref{exa: word smaller than its phi}, and suppose the
exposure order on $F$ has $1\prec x\prec xy\prec y\prec...$ as the
smallest words. Note that this order is not suffix-invariant, since
$1\prec x$ but $y\succ xy$. Since $x$ and $1$ lie in the same
$I$-coset, we have $x\notin T_{I,\prec}'$. However, $xy\in T_{I,\prec}'$
since no element with smaller support -- i.e., of the form $\alpha+\beta x$
with $\alpha,\beta\in K$ -- lies in the same coset as $xy$. To
see this, consider the ring homomorphism $\varphi:\mathcal{A}\rightarrow\mathcal{A}$
which fixes $K$ and maps $\varphi\left(x\right)=1,\varphi\left(y\right)=y$.
Since $\varphi\left(x-1\right)=0$, we have $I\subseteq\text{ker}\varphi$
, so any two elements of the same $I$-coset must have the same image
under $\varphi$. But $\varphi\left(xy\right)=y\neq\alpha+\beta=\varphi\left(\alpha+\beta x\right)$,
so $\alpha+\beta x\notin xy+I$.
\end{example}
We are now ready for the main result of this section.
\begin{thm}
\label{thm: suffix preserving order main theorem}Let $\prec$ be
a suffix-invariant order on $F$. Then:
\begin{enumerate}
\item \label{enu: T_I'=00003DT_I}The minimal Schreier transversal $T_{I}$
consists precisely of those words $u\in F$ that have minimal support
in their $I$-coset.
\item \label{enu: remainder is smallest in I-coset}For every $f\in\mathcal{A}$,
the remainder $\phi_{I}\left(f\right)$ is the unique element of $f+I$
with minimal support.
\item For every $\alpha<\omega_{I}$, the element $f_{\alpha}$ is the unique
monic element of $I\backslash I_{\alpha}$ with minimal support. 
\end{enumerate}
\end{thm}
\begin{proof}
We omit the order from the notation of $T_{I,\prec}'$ and simply
write $T_{I}'$ (as previously done with $T_{I}$) .
\begin{enumerate}
\item By Claims \ref{cla: T_I_prime is a partial schreier transversal}
and \ref{cla: set of smallest words in their coset is a prefix closed},
$T_{I}'$ is a partial Schreier transversal. By Corollary \ref{cor: T_I  is minimal among partial Schreier transversals for I},
we have $T_{I}\preceq_{\text{min}}T_{I}'$. \\
Suppose for contradiction that $T_{I}\prec_{\text{min}}T_{I}'$, and
let $w=\min\left(T_{I}\triangle T_{I}'\right)\in T_{I}\backslash T_{I}'$.
Let $f\in w+I$ be the element of minimal support in its coset. Since
$w\notin T_{I}'$, it follows that $\text{supp}\left(f\right)\prec_{\text{\ensuremath{\max}}}\left\{ w\right\} $.
By Claim \ref{claim: word in the support of smallest in its I-coset is also smallest in its I-coset},
we have $\text{supp}\left(f\right)\subseteq T_{I}'$, and by the minimality
of $w$, also $\text{supp}\left(f\right)\subseteq T_{I}$. But now
$w-f$ is a nonzero element of $I$ supported on $T_{I}$, in contradiction.
\item Let $f\in\mathcal{A}$ and let $g\in f+I$ be an element of minimal
support in its coset. By Claim \ref{claim: word in the support of smallest in its I-coset is also smallest in its I-coset},
we know $\text{supp}\left(g\right)\subseteq T_{I}'$. By part \ref{enu: T_I'=00003DT_I},
$T_{I}'=T_{I}$ so $g=\phi_{I}\left(f\right).$
\item Let $\alpha<\omega_{I}$. The ideal $I_{\alpha}$ is strictly contained
in $I$, so by Proposition \ref{prop: unique monic element supported on T with minimal support}
(applied to $T'=F$), there exists a unique monic element $g_{\alpha}\in I\backslash I_{\alpha}$
with minimal support. \\
We claim that $g_{\alpha}$ is already supported on $T_{I_{\alpha}}$.
Indeed, consider the remainder $\phi_{I_{\alpha}}\left(g_{\alpha}\right)$.
Since $g_{\alpha}\in I\backslash I_{\alpha}$, we have $\phi_{I_{\alpha}}\left(g_{\alpha}\right)\in I\backslash I_{\alpha}$
as well, and by part \ref{enu: remainder is smallest in I-coset},
it has support less than or equal to that of $g_{\alpha}$. By the
uniqueness of $g_{\alpha}$, it follows that $\phi_{I_{\alpha}}\left(g_{\alpha}\right)$
must be a scalar multiple of $g_{\alpha}$, and hence $g_{\alpha}$
is supported on $T_{I_{\alpha}}$.
\end{enumerate}
\end{proof}
Theorem \ref{thm: suffix preserving order main theorem} might fail
when working with an exposure order which is not suffix-invariant.
This was illustrated in Example \ref{exa: word smaller than its phi}. 
\begin{rem}
Theorem \ref{thm: suffix preserving order main theorem} can alternatively
be proven via an algorithmic argument. Specifically, one can show
that each reduction step in Algorithm \ref{alg:reduction_mod_TI},
which by Theorem \ref{thm: Algorithm: reduction modulo T_I given combinatorially reducing system for I}
computes $\phi_{I}\left(f\right)$, strictly $\prec_{\text{max}}$-decreases
the support. To see this, consider a reduction step of the form $r\leftarrow r-\lambda q_{u}u_{\text{suf}}$,
for some $u\in\text{supp}\left(r\right)$. Since $u=\text{HT}\left(q_{u}\right)\cdot u_{\text{suf}}$,
with no cancellation in the product, suffix-invariance implies $\text{HT}\left(\lambda q_{u}u_{\text{suf}}\right)=\text{HT}\left(q_{u}\right)u_{\text{suf}}=u.$
Therefore, this step removes $u$ from $\text{supp}\left(r\right)$,
while possibly introducing only words strictly smaller than $u$,
coming from $\text{supp}\left(q_{u}u_{\text{suf}}\right)$. Thus,
the support of $r$ strictly decreases at each step. Now consider
an element $f$ supported on $T_{I}$, and let $g\in f+I$. Applying
the algorithm to $g$ returns $\phi_{I}\left(g\right)=f.$ Hence,
$f=\phi_{I}\left(f\right)$ is the unique element with minimal support
in $f+I$.
\end{rem}
We finish by describing the resulting simplified recursive construction
of the exposure basis $B_{I}$ for suffix-invariant orders. By Theorem
\ref{thm: suffix preserving order main theorem}, the construction
takes the following simple form: Initialize $B=\emptyset$ , and as
long as the current ideal $I'$ generated by $B$ is not $I$, add
to $B$ the unique monic element of $I\backslash I'$ with minimal
support. In particular, for a shortlex order $\prec$ and a finitely
generated $I$, the resulting exposure basis $B_{I,\prec}$ coincides
with the Rosenmann basis from \cite{Rosenmann1993}, as both bases
satisfy the same uniquely-defining minimality.

\section{\label{sec: Algorithmic Framework}Algorithmic Framework}

\subsection{\label{sec: Extending an Existing Exposure Basis}Extending an Exposure
Basis}

Fix once again a general exposure order $\prec$ on $F$, which we
maintain until the end of the paper. In this section we prove Theorem
\ref{thm: Extending an exposure basis}, which provides a form of
converse for the combinatorial properties satisfied by the \texttt{firsts}
$f_{\alpha}$ and \texttt{seconds} $s_{\alpha}$. This result will
play a central role in the design of subsequent algorithms (see Section
\ref{sec: Algorithm for Computing Bases}). 

Let $I\lneqq\mathcal{A}$ be a proper ideal with exposure basis $B_{I}$,
and let $f\notin I$. Consider the larger ideal $I'$ generated by
$I$ and $f$. We wish to provide a combinatorial criterion that determines
whether the exposure process from Definition \ref{def: inductive definition for ideal},
when applied to the larger ideal $I'$, initially agrees with the
process for $I$ -- reproducing all of $B_{I}$ in order -- and
then exposes $f$, so that $B_{I'}=B_{I}\cup\left\{ f\right\} $.

We make the following observations: 
\begin{itemize}
\item For $f$ to extend the current exposure basis $B_{I}$, it must be
monic, supported on the current Schreier transversal $T_{I}$, and
cannot belong to $K$ since the larger ideal should remain proper.
\item By Proposition \ref{prop: Firsts are increasing}, $f$ must satisfy
$\text{HT}\left(f\right)\succ\text{HT}\left(f'\right)$ for all $f'\in B_{I}$. 
\item Furthermore, its $I$-\texttt{second} $s$ must satisfy $\text{HT}\left(s\right)\succ\text{HT}\left(f\right)$,
as shown in Theorem \ref{Thm: properties of second + Grobner basis is combinatorially reducing}
(\ref{enu: HT and HLL of second}).
\end{itemize}
We will show that these properties are sufficient.
\begin{defn}
\label{def: exposure extending}Let $I\lneqq\mathcal{A}$ and let
$f\in\mathcal{A}\backslash K$ be monic and supported on $T_{I}$.
Let $s$ be the $I$-\texttt{second} of $f$. We say that $f$ is
\emph{exposure-extending} \emph{for $I$} with respect to the fixed
exposure order $\prec$ if $\text{HT}\left(f\right)\succ\text{HT}\left(f'\right)$
for all $f'\in B_{I}$ and $\text{HT}\left(s\right)\succ\text{HT}\left(f\right)$. 
\end{defn}
The condition $\text{HT}\left(s\right)\succ\text{HT}\left(f\right)$
in Definition \ref{def: exposure extending} is essential to ensure
the $\prec_{\text{max}}$-minimality of the support of $f$ among
nonzero elements of $I$ supported on $T_{I}$. The following (easy)
example shows that this condition does not follow from the other requirements
on $f$, and must be explicitly stated.
\begin{example}
Consider $I=\left\{ 0\right\} $ and $f=x^{3}-x$, with respect to
any exposure order. Then $f$ is monic, supported on $T_{I}=F$ and
does not lie in $K$. However, its $I$-\texttt{second} $s$ is equal
to $x^{2}-1$, so $\text{HT}\left(s\right)=x^{2}\prec x^{3}=\text{HT}\left(f\right)$. 
\end{example}
We begin by characterizing the words in the support of the $I$\texttt{-second}
$s$ of an element $f$ that are greater than $\text{HT}\left(f\right)$.
In the following statement, note that we do not assume that $\text{HT}\left(s\right)\succ\text{HT}\left(f\right)$.
\begin{lem}
\label{lem: words in support of s larger than HT(f)}Let $I\lneqq\mathcal{A}$
and let $f\in\mathcal{A}\backslash K$ be supported on $T_{I}$ such
that $\text{HT}\left(f\right)\succ\text{HT}\left(f'\right)$ for all
$f'\in B_{I}$. Let $b=\text{HTT}\left(f\right)$, and let $s$ be
the $I$-\texttt{second} of $f$. Then every word $u\in\text{supp}\left(s\right)$
such that $u\succeq\text{HT}\left(f\right)$ ends in $b^{-1}$ and
satisfies $\text{HT}\left(f\right)\succeq ub$. In particular, $\text{HT}\left(s\right)\neq\text{HT}\left(f\right)$.
\end{lem}
\begin{proof}
We proceed by transfinite induction. Assume that for every $\alpha<\omega_{I}$,
the lemma holds for $I$$_{\alpha}$. In particular, the lemma applies
to $s_{\alpha}$ as the $I_{\alpha}$-\texttt{second} of $f_{\alpha}$
since, by Proposition \ref{prop: Firsts are increasing}, $f_{\alpha}$
has a strictly greater head term than all of $\left\{ f_{\beta}:\beta<\alpha\right\} $.

To prove the lemma for $I$, recall that by Theorem \ref{thm: introduction Grobner basis is combinatorially reducing },
$B_{I}^{\text{gr}}$ is a CRS for $I$. Hence, by Corollary \ref{cor: remainder of f modulo T_I is linear combination of Q and has smaller support},
$s$ is a $K$-linear combination of $fb^{-1}$ and elements of $\left\{ q\in B_{I}^{\text{gr}}:\text{HTT}\left(q\right)=b^{-1}\right\} $.
Let $u\in\text{supp}\left(s\right)$ such that $u\succeq\text{HT}\left(f\right)$.
We distinguish between the following cases:

\uline{Case 1:} $u\in\text{supp}\left(fb^{-1}\right)$. In this
case, $ub\in\text{supp}\left(f\right)$, so $\text{\text{HT}}\left(f\right)\succeq ub$.
By transitivity, we also obtain $u\succeq ub$. It follows that $u$
ends in $b^{-1}$.

\uline{Case 2:} $u\in\text{supp}\left(f'\right)$ for some $f'\in B_{I}$.
In this case, $u\preceq\text{\text{HT}}\left(f'\right)\prec\text{HT}\left(f\right)$,
which contradicts the assumption that $u\succeq\text{HT}\left(f\right)$.
Therefore, this case is impossible.

\uline{Case 3:} $u\in\text{supp}\left(s_{\alpha}\right)$ for some
$\alpha<\omega_{I}$ such that $b_{\alpha}=b$. Since the lemma holds
for $s_{\alpha}$ and $u\succeq\text{HT}\left(f\right)\succ\text{HT}\left(f_{\alpha}\right)$,
we conclude that $u$ ends in $b_{\alpha}^{-1}=b^{-1}$ and satisfies
$\text{HT}\left(f\right)\succ\text{HT}\left(f_{\alpha}\right)\succeq ub$,
as needed.

Finally, suppose for contradiction that $\text{HT}\left(s\right)=\text{HT}\left(f\right)$.
On the one hand, $\text{HT}\left(f\right)$ ends in $b$. On the other
hand, applying the lemma to $u=\text{HT}\left(s\right)\succeq\text{HT}\left(f\right)$,
we find that $\text{HT}\left(s\right)$ ends in $b^{-1}$, in contradiction.
Therefore, $\text{HT}\left(s\right)\neq\text{HT}\left(f\right)$.
\end{proof}
We are now ready to show the main results of this section.
\begin{lem}
\label{lem: exposure-extending f extends combinatorially reducing system}Let
$I\lneqq\mathcal{A}$. Let $f$ be an exposure-extending element for
$I$, and let $s$ be its $I$-\texttt{second}. Then $B_{I}^{\text{gr}}\cup\left\{ f,s\right\} $
is a combinatorially reducing system for the ideal generated by $I$
and $f$. 
\end{lem}
\begin{proof}
Let $\mathcal{Q}:=B_{I}^{\text{gr}}\cup\left\{ f,s\right\} $. Since
$B_{I}^{\text{gr}}$ generates $I$, the set $B_{I}^{\text{gr}}\cup\left\{ f\right\} $
generates the larger ideal generated by $I$ and $f$. Moreover, since
$s$ already lies in this ideal, $\mathcal{Q}$ generates the same
ideal. 

We now show that $\mathcal{Q}$ is a CRS. By Theorem \ref{thm: introduction Grobner basis is combinatorially reducing },
$B_{I}^{\text{gr}}$ is a CRS for $I$. Hence, $B_{I}^{\text{gr}}\subseteq\mathcal{A}\backslash K$
and its elements are monic. Moreover, $f$ is a monic element of $\mathcal{A}\backslash K$
since it is exposure-extending. Its $I$-\texttt{second} $s$ is,
by definition, monic. Since $\text{HT}\left(s\right)\succ\text{HT}\left(f\right)$,
we have $\text{HT}\left(s\right)\notin K$ (as $e$ is a minimal element
of $F$ in any exposure order). Thus, $\mathcal{Q}\subseteq\mathcal{A}\backslash K$
and its elements are monic.

We next analyze the new \texttt{second} $s$. Let $b_{f}=\text{HTT}\left(f\right)$.
Since $f$ is exposure-extending for $I$, and in particular $\text{HT}\left(s\right)\succ\text{HT}\left(f\right)$,
we can apply Lemma \ref{lem: words in support of s larger than HT(f)}
to $u=\text{HT}\left(s\right)$, which satisfies $u\succ\text{HT}\left(f\right)$.
It follows that $\text{HT}\left(s\right)$ ends in $b_{f}^{-1}$ and
satisfies $\text{HT}\left(s\right)b_{f}\preceq\text{HT}\left(f\right)$.
Moreover, since $\text{HTT}\left(s\right)=\text{HTT}\left(f\right)^{-1}$,
by Corollary \ref{cor: remainder of f modulo T_I is linear combination of Q and has smaller support}
(with respect to the CRS $B_{I}^{\text{gr}}$ for $I$), it follows
that $f$ is also the $I$-\texttt{second} of $s$. Hence, for any
$q\in\left\{ f,s\right\} $, its $I$-\texttt{second} also lies in
$\left\{ f,s\right\} $, so by the same Corollary, $q\cdot\text{HTT}\left(q\right)^{-1}$
is a $K$-linear combination of the elements $\left\{ q'\in B_{I}^{\text{gr}}:\text{HTT}\left(q'\right)=\text{HTT}\left(q\right)^{-1}\right\} $
and this $I$-\texttt{second}, all of which lie in $\mathcal{Q}$.
We now use these observations to establish the three conditions required
on $\mathcal{Q}$ to be a CRS.

First, we show that $\left(\text{HT}\left(q\right)\right)_{q\in\mathcal{Q}}$
are distinct. Since $B_{I}^{\text{gr}}$ is already a CRS, every pair
of head terms among $\left(w_{\alpha}^{\pm}\right)_{\alpha<\omega_{I}}$
are distinct. Furthermore, since $\text{HT}\left(f\right)$ and $\text{HT}\left(s\right)$
lie in $T_{I}$, they are distinct from each of the words $\left(w_{\alpha}^{\pm}\right)_{\alpha<\omega_{I}}$,
which are all prefix-neighbors of $T_{I}$ by Proposition \ref{prop: Form for first using T_I}
and Theorem \ref{Thm: properties of second + Grobner basis is combinatorially reducing}
(\ref{enu: HT and HLL of second}). Finally, $\text{HT}\left(s\right)\succ\text{HT}\left(f\right)$.

Second, let $q\in\mathcal{Q}$, and let $b_{q}=\text{HTT}\left(q\right)$.
We claim that $qb_{q}^{-1}$ is a $K$-linear combination of elements
$\left\{ q'\in\mathcal{Q}:\text{HTT}\left(q'\right)=b_{q}^{-1}\right\} $.
If $q\in\left\{ f,s\right\} ,$ this was already argued above. Otherwise,
$q\in B_{I}^{\text{gr}}$. Since $B_{I}^{\text{gr}}$ is a CRS, $qb_{q}^{-1}$
is such a combination, but with $q'\in B_{I}^{\text{gr}}$. This is
as required, since $B_{I}^{\text{gr}}\subseteq\mathcal{Q}$. 

Finally, let $q\in\mathcal{Q}$. To conclude that $\text{supp}\left(q\right)\subseteq T_{\mathcal{Q}}^{\text{HT}}\cup\partial T_{\mathcal{Q}}^{\text{HT}}$,
let $u\in\text{supp}\left(q\right)$ and $q'\in\mathcal{Q}$ such
that $\text{HT}\left(q'\right)$ is a prefix of $u$. We must show
that $u=\text{HT}\left(q'\right)$. Suppose towards contradiction
that $\text{HT}\left(q'\right)$ is a proper prefix of $u$. Since
$u\in\text{supp}\left(q\right)$, we then have $\text{HT}\left(q\right)\succeq u\succ\text{HT}\left(q'\right)$.
We split into three cases cases based on $q$, and show that each
leads to a contradiction.

\uline{case 1:} $q\in\left\{ f,s\right\} $. Then $q$ is supported
on $T_{I}=T_{B_{I}^{\text{gr}}}^{\text{HT}}$. Since $\text{HT}\left(q'\right)$
is a prefix of $u\in\text{supp}\left(q\right)$, it cannot lie in
$\partial T_{I}$, so $q'\notin B_{I}^{\text{gr}}$. Thus, $q'\in\left\{ f,s\right\} $
as well. Since $\text{HT}\left(q\right)\succ\text{HT}\left(q'\right)$,
we have $q=s$ and $q'=f$. Therefore, $u\in\text{supp}\left(s\right)$
and satisfies $u\succ\text{HT}\left(q'\right)=\text{HT}\left(f\right)$.
By Lemma \ref{lem: words in support of s larger than HT(f)}, $u$
ends in $b_{f}^{-1}$ and satisfies $\text{HT}\left(f\right)\succeq ub_{f}.$
Thus, $\text{HT}\left(f\right)=\text{HT}\left(q'\right)$ is a proper
prefix of $u$ but not a proper prefix of $ub_{f}$. It follows that
$\text{HT}\left(f\right)=ub_{f}$. This leads to a contradiction,
as then $u$ ends in consecutive $b_{f}^{-1}b_{f}$, which is impossible.

\uline{case 2:} $q=f_{\alpha}$ for some $\alpha<\omega_{I}$.
In this case, since $f$ is exposure-extending, we have $\text{HT}\left(s\right)\succ\text{HT}\left(f\right)\succ\text{HT}\left(q\right)\succ\text{HT}\left(q'\right)$.
It follows that $q'\notin\left\{ f,s\right\} $, so $q'\in B_{I}^{\text{gr}}.$
Therefore $u$, which has $\text{HT}\left(q'\right)$ as a proper
prefix of $u$, does not lie in $T_{B_{I}^{\text{gr}}}^{\text{HT}}\cup\partial T_{B_{I}^{\text{gr}}}^{\text{HT}}$.
But $B_{I}^{\text{gr}}$ is a CRS by Theorem \ref{thm: introduction Grobner basis is combinatorially reducing },
so $q\in B_{I}^{\text{gr}}$ is supported on $T_{B_{I}^{\text{gr}}}^{\text{HT}}\cup\partial T_{B_{I}^{\text{gr}}}^{\text{HT}}$,
in contradiction.

\uline{case 3:} $q=s_{\alpha}$ for some $\alpha<\omega_{I}$.
As in the previous case, if $q'\in B_{I}^{\text{gr}}$, we are done.
If, instead, $q'\in\left\{ f,s\right\} $, then $u\in\text{supp}\left(s_{\alpha}\right)$
satisfies $u\succ\text{HT}\left(q'\right)\succeq\text{HT}\left(f\right)\succ\text{HT}\left(f_{\alpha}\right)$,
so by Lemma \ref{lem: words in support of s larger than HT(f)} $u$
ends in $b_{\alpha}^{-1}$ and satisfies $\text{HT}\left(f_{\alpha}\right)\succeq ub_{\alpha}$.
But then $u\succ\text{HT}\left(q'\right)\succ ub_{\alpha}$, so $\text{HT}\left(q'\right)$
is a proper prefix of $u$ which is not a proper prefix of $ub_{\alpha}$,
in contradiction.
\end{proof}
We now prove Theorem \ref{thm: Extending an exposure basis}, the
main result of this section.

\begin{proof}[Proof of Theorem \ref{thm: Extending an exposure basis}]Let
$\mathcal{Q}:=B_{I}^{\text{gr}}$ and $\mathcal{Q}':=\mathcal{Q}\cup\left\{ f,s\right\} $.
By Lemma \ref{lem: exposure-extending f extends combinatorially reducing system},
$\mathcal{Q}'$ is a CRS for the ideal $I'$. Hence, by Theorem \ref{thm: combinatorially reducing system defines minimal schreier transversal of the ideal which it generates},
we have $T_{I'}=T_{\mathcal{Q}'}^{\text{HT}}$. In particular, $T_{I'}$
is obtained from $T_{I}=T_{\mathcal{Q}}^{\text{HT}}$ by removing
all words having $\text{HT}\left(f\right)$ or $\text{HT}\left(s\right)$
as a prefix. Since $\text{HT}\left(f\right)\prec\text{HT}\left(s\right)$,
we conclude that $\min\left(T_{I}\triangle T_{I'}\right)=\text{HT}\left(f\right)$.

Let $\left(f_{\alpha}\right)_{\alpha<\omega_{I}}$ and $\left(f_{\alpha}'\right)_{\alpha<\omega_{I'}}$
be the respective exposure basis elements of $I$ and $I$'. We first
prove that the two exposure processes coincide up to ordinal $\omega_{I}$;
that is, for every $\alpha<\omega_{I}$, the process for $I'$ does
not stop at $\alpha$, and $f_{\alpha}'=f_{\alpha}$. Proceeding by
transfinite induction, let $\alpha<\omega_{I}$ and suppose the claim
holds for every $\beta<\alpha$. Then $I_{\alpha}'=I_{\alpha}$ since
both are defined from earlier steps. Since $\alpha<\omega_{I}$, the
process for $I$ does not stop at $\alpha$, so $I_{\alpha}\subsetneqq I\subsetneqq I'$,
which implies that the process for $I'$ does not stop at $\alpha$
either.

Now consider the head terms of the elements $f_{\alpha}$ and $f_{\alpha}'$
added at ordinal $\alpha$. By Proposition \ref{prop: Form for first using T_I},
we have 
\[
\text{HT}\left(f_{\alpha}\right)=\min\left(T_{I}\triangle T_{I_{\alpha}}\right),\ \text{HT}\left(f_{\alpha}'\right)=\min\left(T_{I}'\triangle T_{I_{\alpha}}\right).
\]
 Since $f$ is exposure-extending for $I$, we know $\text{HT}\left(f\right)\succ\text{HT}\left(f_{\alpha}\right)$.
Thus the sets $T_{I'}\triangle T_{I}$ and $T_{I}\triangle T_{I_{\alpha}}$
have distinct minima. Then, by Claim \ref{claim: minimum of symmetric difference},
\begin{align*}
\text{HT}\left(f_{\alpha}'\right) & =\min\left(\left(T_{I'}\triangle T_{I}\right)\triangle\left(T_{I}\triangle T_{I_{\alpha}}\right)\right)=\text{HT}\left(f_{\alpha}\right).
\end{align*}
Since both $f_{\alpha}$ and $f_{\alpha}'$ are monic elements of
$I'$, if $f_{\alpha}\neq f_{\alpha}'$, then their difference $f_{\alpha}'-f_{\alpha}$
is a nonzero element of $I'$ with support strictly smaller than $\text{HT}\left(f_{\alpha}'\right)$,
contradicting the minimality defining $f_{\alpha}'$. Hence, $f_{\alpha}'=f_{\alpha}$.

Now consider the ordinal $\omega_{I}$. By the inductive argument
above, 
\[
\left\{ f_{\alpha}':\alpha<\omega_{I}\right\} =\left\{ f_{\alpha}:\alpha<\omega_{I}\right\} =B_{I}.
\]
The ideal $I$, generated by $B_{I}$, is strictly contained in $I'$,
since $f\in I'$ is a nonzero element supported on $T_{I}$. Therefore,
the exposure process for $I'$ does not stop at $\omega_{I}$, and
proceeds to define $f_{\omega_{I}}'$. By Proposition \ref{prop: Form for first using T_I},
this new \texttt{first} satisfies 
\[
\text{HT}\left(f_{\omega_{I}}'\right)=\min\left(T_{I}\triangle T_{I'}\right)=\text{HT}\left(f\right).
\]
 Again, both $f_{\omega_{I}}'$ and $f$ are monic elements of $I'$
supported on $T_{I}$ with the same head term. If $f_{\omega_{I}}'\neq f$,
then $f_{\omega_{I}}'-f$ would contradict the minimality of $f_{\omega_{I}}'$.
Hence, $f_{\omega_{I}}'=f$. Its associated \texttt{second} $s_{\omega_{I}}'$
is the $I$-\texttt{second} of $f$, which is $s$.

Finally, since $I'$ is generated by $I\cup\left\{ f\right\} $, the
exposure process must stop at the next ordinal: $\omega_{I'}=\omega_{I}+1$.
Therefore, the exposure basis and Gröbner bases for $I'$ are:
\begin{align*}
B_{I'} & =\left\{ f_{\alpha}':\alpha<\omega_{I'}\right\} =B_{I}\cup\left\{ f\right\} ,\\
B_{I'}^{\text{gr}} & =\left\{ f_{\alpha}':\alpha<\omega_{I'}\right\} \cup\left\{ s_{\alpha}':\alpha<\omega_{I'}\right\} =B_{I}^{\text{gr}}\cup\left\{ f,s\right\} .
\end{align*}
\end{proof}

\subsection{\label{sec: Algorithm for Computing Bases}Algorithms}

In this section we provide several algorithms for working with the
theoretical framework developed in the previous sections, in addition
to Algorithm \ref{alg:reduction_mod_TI}. Since our objectives are
now algorithmic in nature, we assume that both the field operations
in $K$ and the comparisons with respect to the fixed exposure order
$\prec$ are computable. Explicitly, for the latter, we assume that
there exists an algorithm which, given two words $u,v\in F$, decides
if $u\prec v$.

\subsubsection{Algorithm for computing $B_{I}$ and $B_{I}^{\text{gr}}$}

We begin with the main result of this section: an algorithm for computing
both the exposure basis $B_{I}$ and the Gröbner basis $B_{I}^{\text{gr}}$
of a finitely generated ideal $I$, given a finite generating set.
Recall that these bases depend on the fixed exposure order $\prec$
on $F$. In the degenerate case where $I=\mathcal{A}$, the exposure
basis $B_{I}$ consists solely of the unit $e\in\mathcal{A}$, and
we define $B_{I}^{\text{gr}}:=\left\{ e\right\} $.

The algorithm receives a finite generating set $h_{1},h_{2},...,h_{m}$
for $I$ and maintains the following lists:
\begin{itemize}
\item \emph{$B_{J}$} - A list of elements of $\mathcal{A}$ which, at each
step, forms the exposure basis of a ``current'' sub-ideal $J\subseteq I$.
This ideal $J$ is defined throughout the algorithm as the one generated
by the current contents of $B_{J}$, so it is updated whenever $B_{J}$
is updated. The supports of $B_{J}$ are distinct by Proposition \ref{prop: Firsts are increasing},
and the list is sorted in ascending order by these supports with respect
to the order $\prec_{\text{max}}$.
\item \emph{$S_{J}$} - A list of the corresponding \texttt{seconds} associated
with the elements of $B_{J}$ in the exposure process of the ideal
$J$, as in Definition \ref{def: (HLL of first, second, highest term of second)}.
The list $S_{J}$ follows the ordering of $B_{J}$, meaning that $S_{J}\left[i\right]$
is the \texttt{second} of $B_{J}\left[i\right]$.
\item \emph{Gens} - a list of yet ``unhandled'' generators. The elements
of $B_{J}$\emph{ }and\emph{ Gens} together generate $I$ throughout
the run of the algorithm. Note that elements of \emph{Gens} may have
non-distinct supports.
\end{itemize}
We give an overview of the algorithm. Initially, \emph{Gens}$=\left[h_{1},h_{2},...,h_{m}\right]$
and $B_{J}=S_{J}=[]$, so that $J=\left\{ 0\right\} $. The elements
of \emph{Gens} are then processed one by one, updating the lists $B_{J},S_{J}$
and \emph{Gens} accordingly until, upon termination, \emph{Gens} is
empty and $B_{J}$ generates all of $I$. The processing of an unhandled
generator $f$ from \emph{Gens} begins by discarding it if already
lies in the current sub-ideal $J$. Otherwise, the algorithm transforms
$f$, and may move some elements from $B_{J}$ back to \emph{Gens},
until $f$ becomes exposure-extending for $J$. When this goal is
satisfied, $f$ is appended to $B_{J}$, extending it validly by Theorem
\ref{thm: Extending an exposure basis}, and the corresponding $J$-\texttt{second}
is appended to $S_{J}$. A special case is monitored: if at any point
$\left|\text{supp}\left(f\right)\right|=1$, then $f$ already generates
all of $\mathcal{A}$, and the algorithm terminates. We now describe
the steps of how each element $f$ is processed and reduced.
\begin{enumerate}
\item Ensure that $f$ is supported on $T_{J}$: Replace $f$ with its remainder
modulo $T_{J}$, i.e., $f\leftarrow\phi_{J}\left(f\right)$, using
Algorithm \ref{alg:reduction_mod_TI}, with the Gröbner basis $B_{J}^{\text{gr}}=B_{J}\cup S_{J}$
(which, by Theorem \ref{Thm: properties of second + Grobner basis is combinatorially reducing},
forms a CRS for $J$). If this reduction yields $f=0$, then the original
$f$ already lies in $J$, and we discard it and continue to the next
element in \emph{Gens}.
\item Make $f$ monic: Replace $f$ with its monic $K^{\times}$-multiple.
From this point onward, $f$ will remain monic and supported on $T_{J}$.
\item Check if $f$ generates all of $\mathcal{A}$: If $\left|\text{supp}\left(f\right)\right|=1$,
then $f$ already generates all of $\mathcal{A}$, so the algorithm
terminates with $B_{I}=B_{I}^{\text{gr}}=\left\{ e\right\} $. This
check is slightly stronger than necessary -- to proceed with computing
the $J$-\texttt{second} later on, it suffices to verify that $f\neq e$
(which, since $f$ is monic, is equivalent to $f\notin K$).
\item Adjust $B_{J}$ for compatibility with $f$: Let $L$ be the final
segment of $B_{J}$ consisting of elements $g$ such that $\text{supp}\left(f\right)\prec_{\text{max}}\text{supp}\left(g\right)$.
We remove $L$ from the end of $B_{J}$, delete the corresponding
\texttt{seconds} from $S_{J}$, and prepend the elements of $L$ to
the list \emph{Gens} for later reprocessing. To justify this step,
let $B_{J'}$ be the initial segment of $B_{J}$ preceding $L$, and
let $J'$ be the ideal generated by $B_{J'}$. By Claim \ref{claim: Exposure basis for an ideal along the way},
$B_{J'}$ is the exposure basis of $J'$, and its \texttt{seconds}
are the corresponding initial segment of $S_{J}$. By Theorems \ref{Thm: properties of second + Grobner basis is combinatorially reducing}
and \ref{thm: combinatorially reducing system defines minimal schreier transversal of the ideal which it generates},
and since $B_{J'}^{\text{gr}}\subseteq B_{J}^{\text{gr}}$, it follows
that $T_{J'}\supseteq T_{J}$. In particular, since $f$ is supported
on $T_{J}$, it is also supported on $T_{J'}$. If $L$ is nonempty,
then $f$ has $\prec_{\text{max}}$-smaller support than the front
element of $L$. This contradicts the minimality required for the
front element of $L$ to follow $B_{J'}$ in the exposure process
of the input ideal $I$. Thus, removing $L$ and restoring its elements
to \emph{Gens} preserves correctness. After this step, $f$ is a monic
element of $\mathcal{A}\backslash K$ supported on $T_{J}$, and now
satisfies $\text{HT}\left(f\right)\succ\text{HT}\left(g\right)$ for
all $g\in B_{J}$.
\item Ensure that $f$ has $\prec_{\text{max}}$-smaller support than its
$J$-\texttt{second}: Compute the $J$-\texttt{second} $s$ of $f$
(computing again $\phi_{J}$ similarly to step $1$). If $\text{HT}\left(s\right)\succ\text{HT}\left(f\right)$,
then $f$ is exposure-extending for $J$. Thus, $f$ is appended to
$B_{J}$, its $J$-\texttt{second} $s$ appended to $S_{J}$, and
the algorithm continues to pop the next generator from \emph{Gens}.
Otherwise, by Lemma \ref{lem: words in support of s larger than HT(f)}
we have $\text{HT}\left(s\right)\neq\text{HT}\left(f\right)$, so
$\text{HT\ensuremath{\left(s\right)}}\prec\text{HT}\left(f\right)$.
In this case, we replace $f$ with $s$. Since this results in the
support of $f$ strictly $\prec_{\text{max}}$-decreasing, we return
to the third step - checking if $\left|\text{supp}\left(f\right)\right|=1$
- and repeat. Note that steps $1$ and $2$ are skipped as the $J$-\texttt{second}
$s$ is by definition monic and supported on $T_{J}$.
\end{enumerate}
\begin{thm}
\textbf{(Algorithm for Computing the Exposure and Gröbner Bases)}
Let $f_{1},f_{2},...,f_{m}$ be a finite generating set for a right
ideal $I\leq\mathcal{A}$. Then Algorithm \ref{alg:exposure_grobner}
computes both the exposure basis $B_{I}$ and the Gröbner basis $B_{I}^{\text{gr}}$
associated to $I$ with respect to the fixed exposure order $\prec$
on $F$.
\end{thm}
\begin{proof}
Throughout the execution of the algorithm, the ideal generated by
$B_{J}\cup$\emph{Gens} remains unchanged since the following operations
preserve it: (i) right-multiplication of a generator by an invertible
element of $\mathcal{A}$ (e.g., making elements monic or multiplying
by $b^{-1}$ for $b\in S\cup S^{-1}$), (ii) replacing a generator
$f$ by its remainder $\phi_{J}\left(f\right)$, where $J$ is the
ideal generated by other generators, and (iii) discarding generators
equal to zero. 

By Theorem \ref{thm: Extending an exposure basis}, any element $f$
added to $B_{J}$ satisfies the conditions required for it to extend
$B_{J}$ as an exposure basis. Therefore, $B_{J}$ always remains
the exposure basis for the ideal it generates, and $S_{J}$ contains
the corresponding \texttt{seconds}. If the algorithm terminates with
\emph{Gens} empty, then $B_{J}$ generates the ideal $I$, and hence
$B_{I}=B_{J}$ and $B_{I}^{\text{gr}}=B_{J}\cup S_{J}$.

We now show that the algorithm terminates. If at any point a generator
with support size $1$ is encountered, the algorithm halts immediately,
returning $B_{I}=B_{I}^{\text{gr}}=\left\{ 1\right\} $ as required.
Furthermore, since no new generators are created, generators can be
discarded at most $m$ times. After finitely many steps, all such
removals cease. 

To show termination of the remaining process, consider the list $Z$
formed by concatenating 3 lists: the supports of the elements in $B_{J}$,
followed by the support of the current element $f$, followed by the
supports of the elements in \emph{Gens}. The length of $Z$ is fixed
(say, $m'$) after the initial finitely many discards. Elements of
\emph{Gens} which are not supported on the current minimal Schreier
transversal $T_{J}$ are treated as having infinite support. The resulting
extended set of possible supports -- finite subsets of $F$ ordered
by $\prec_{\text{max}}$, together with the new maximal value $\infty$
-- remains well-ordered. We order such lists lists lexicographically:
$Z<Z'$ if at the first position $i$ where they differ, the support
in $Z_{i}$ is less than the support of $Z'_{i}$ with respect to
$\prec_{\text{max}}$ extended by $\infty$. Since each support is
from a well-ordered set and the lists are of fixed finite length $m'$,
the lexicographic order on such lists is itself a well-order.

We claim that each nontrivial change to the support of the current
processed generator $f$ strictly decreases $Z$ in this lexicographic
order respect to this order:
\begin{itemize}
\item If $f$ is not already supported on $T_{J}$, a step of the form $f\leftarrow\phi_{J}\left(f\right)$
reduces its support from $\infty$ to a finite value.
\item Moving a nontrivial suffix $L$ from the end of $B_{J}$ to the front
of \emph{Gens} results in placing the current $f$, which has smaller
head term than $L$, earlier in $Z$.
\item Replacing $f$ with its \texttt{second} $s$ when $\text{HT}\left(s\right)\prec\text{HT}\left(f\right)$
strictly $\prec_{\text{max}}$-decreases its support. 
\end{itemize}
Since $Z$ takes values in a well-ordered set, it can strictly decrease
only finitely many times. Eventually, no suffixes are moved and no
elements are replaced. At this point, the remaining elements of \emph{Gens}
are simply added to $B_{J}$ one by one, and the algorithm terminates.
\end{proof}
\begin{algorithm}
\caption{\textsc{ExposureAndGröbnerBasis} - Computes Exposure and Gröbner Bases of a Finitely Generated Right Ideal in $\mathcal{A}$}
\label{alg:exposure_grobner}
\begin{algorithmic}[1]

\REQUIRE A finite generating set $h_1, h_2, \ldots, h_m$ for an ideal $I$ in $\mathcal{A}$.
\ENSURE The exposure basis $B_I$ and the Gröbner basis $B_I^{\mathrm{gr}}$ for $I$, with respect to the fixed exposure order $\prec$ on $F$.

\STATE Initialize $B_J \gets []$ , $S_J \gets []$, $\mathrm{Gens} \gets [h_1, h_2, \ldots, h_m]$
\WHILE{$\mathrm{Gens} \neq \emptyset$}
    \STATE $f \gets \text{pop\_front}(\mathrm{Gens})$ \COMMENT{pop\_front(List) removes and returns the first element of List.}
    \STATE $f \gets \phi_J(f)$, where $J$ is the ideal generated by $B_J$ \COMMENT{$f$ is now supported on $T_J$.}
    \IF{$f = 0$}
        \STATE \textbf{continue} \COMMENT{this results in $f$ being removed from Gens}
    \ENDIF
	\STATE $f \gets \text{MONIC}(f)$ \COMMENT{$f$ is now monic and supported on $T_J$.}
	\STATE $f$\texttt{\_exposure\_extending\_for\_}$J$ $\gets$ \texttt{FALSE}
    \WHILE{$f$\texttt{\_exposure\_extending\_for\_}$J$ $\neq $ \texttt{TRUE}}
		\IF{$|\text{supp}(f)| = 1$}
			\STATE Set $B_I \gets \{e\}$ and $B_I^{\mathrm{gr}} \gets \{e\}$.
			\STATE \textbf{return}
	    \ENDIF
	    \STATE $L \gets$ suffix of $B_J$ consisting of elements $g$ with $\text{HT}(g) \succ \text{HT}(f)$.
	    \STATE Move $L$ from end of $B_J$ to the front of $\mathrm{Gens}$  \COMMENT{$f$ is now monic, supported on $T_J$ and $\text{HT}(f)$ is greater than all HTs of $B_J$.}
	    \STATE Remove the last $\left|L\right|$ elements from $S_J$. \COMMENT{Remove corresponding \texttt{seconds}}
	    \STATE $b \gets \mathrm{HTT}(f)$ \COMMENT{Head Term Tail of $f$}
	    \STATE $s \gets \text{MONIC}(\phi_J(fb^{-1}))$, where $J$ is the ideal generated by $B_J$ \COMMENT{$s$ is the $J$-\texttt{second} of $f$}
	    \IF{$\text{HT}(s) \prec \text{HT}(f)$}
	        \STATE $f \gets s$
	    \ELSE[By Lemma \ref{lem: words in support of s larger than HT(f)}, $HT(f)\neq HT(s)$. Thus, if we are here, then $\mathrm{HT}(s) \succ \mathrm{HT}(f)$.]
			\STATE $f$\texttt{\_exposure\_extending\_for\_}$J$ $\gets$ \texttt{TRUE}
	    \ENDIF
	\ENDWHILE
	\STATE Append $f$ to $B_J$ and $s$ to $S_J$.
\ENDWHILE
\STATE Set $B_I \gets B_J$ and $B_I^{\mathrm{gr}} \gets B_J \cup S_J$.
\STATE \textbf{return}

\end{algorithmic}
\end{algorithm}

We state a simpler special version of Algorithm \ref{alg:exposure_grobner}
-- corresponding to the process titled ``Orbit Reduction'' in Rosenmann's
original paper when using the shortlex order \cite{Rosenmann1993}
-- which applies when the ideal is generated by a single element
$h\in\mathcal{A}$.
\begin{thm}
\textbf{(Orbit-Reduction)} Let $I$ be a right ideal generated by
a single element $h\in\mathcal{A}.$ Then the exposure basis $B_{I}$
and Gröbner basis $B_{I}^{\text{gr}}$ associated to $I$ with respect
to the fixed exposure order $\prec$ on $F$ are computed by Algorithm
\ref{alg:single_generator_basis}.
\end{thm}
\begin{algorithm}
\caption{\textsc{OrbitReduction} - Computes Exposure and Gröbner Bases from a Single Generator}
\label{alg:single_generator_basis}
\begin{algorithmic}[1]
\REQUIRE An element $h \in \mathcal{A}$.
\ENSURE The exposure basis $B_I$ for the ideal $I$ generated by $h$. If  $I\lneqq \mathcal{A}$ then the Gröbner basis $B_I^{\mathrm{gr}}$ is returned as well.
\IF{$h = 0$}
    \RETURN $B_I = \emptyset $, $B_I^{\text{gr}} = \emptyset$
\ELSIF{$|\text{supp}(h)| = 1$}
    \RETURN $B_I = \{e\}$, $B_I^{\text{gr}} = \{e\}$
\ENDIF
\STATE $f \gets \text{MONIC}(h)$
\STATE $f$\texttt{\_exposure\_extending\_for\_zero\_ideal} $\gets$ \texttt{FALSE}
\WHILE{$f$\texttt{\_exposure\_extending\_for\_zero\_ideal} $\neq $ \texttt{TRUE}}
	\STATE $b \gets \mathrm{HTT}(f)$ \COMMENT{Head Term Tail of $f$}
	\STATE $s \gets$ the monic $K^{\times}$-multiple of $fb^{-1}$ \COMMENT{$s$ is the $J$-\texttt{second} of $f$, where $J=\left\{ 0\right\} $}.
	\IF{$\text{HT}(s) \prec \text{HT}(f)$}
		\STATE $f \gets s$
	\ELSE[if we are here, then $\mathrm{HT}(s)  \mathrm{HT}(f)$.]
		\STATE $f$\texttt{\_exposure\_extending\_for\_zero\_ideal} $\gets$ \texttt{TRUE}
	\ENDIF
\ENDWHILE
\STATE \textbf{return} $B_I = \{f\}$ and $B_I^{\mathrm{gr}} = \{f, s\}$.
\end{algorithmic}
\end{algorithm}
\begin{proof}
The correctness of the algorithm follows from Algorithm \ref{alg:exposure_grobner},
which it specializes. We justify the modifications by showing how
they naturally arise from the behavior of the general algorithm in
the single-generator case. 

The general algorithm begins with $B_{J}=B_{J}^{\text{gr}}=[]$, corresponding
to the zero ideal $J=\left\{ 0\right\} $, whose minimal Schreier
transversal is $T_{J}=F$ and for which the remainder function is
$\phi_{J}=id_{\mathcal{A}}$. The running element $f$ is initialized
as $h$, and is either discarded if $h=0$, or processed otherwise.
The first case is handled explicitly in this simplified algorithm:
If $h=0$, we terminate immediately, returning $B_{I}=B_{I}^{\text{gr}}=\emptyset$.
For the main processing step, note that whenever $B_{J}$ is empty,
no elements are returned to unhandled status in the list \emph{Gens},
so when $f$ is processed and added to $B_{J}$ the algorithm terminates.
This allows us to eliminate from the algorithm: (i) the loop iterating
on generators, (ii) calls to $\phi_{J}$, as these reduce to the identity
map, and (iii) any steps referring to the final segment $L$ from
$B_{J}$, which is always empty in this case. Lastly, without calls
to $\phi_{J}$, all updates to $f$ consist only of right multiplications
by units of $\mathcal{A}$. In particular, the size of the support
of $f$ remains unchanged throughout the run. This justifies checking
if $\left|\text{supp}\left(f\right)\right|=1$ only once - immediately
after verifying $f\neq0$.
\end{proof}

\subsubsection{\label{subsec: algorithm for division with remainder}Canonical Division
with Remainder with Respect to a Combinatorially Reducing System}

Let $\mathcal{Q}$ be a CRS for a finitely generated ideal $I\lneqq\mathcal{A}$.
In this subsection, we present an algorithm for division by $\mathcal{Q}$
with remainder in the form of Theorem \ref{thm: division with remainder},
with remainder in the minimal Schreier transversal $T_{I}$. The algorithm
generalizes Algorithm \ref{alg:reduction_mod_TI}, which computes
only the remainder $r=\phi_{I}\left(f\right)$ of $f$ modulo $T_{I}$
by repeatedly reducing $f$ using $\mathcal{Q}$ until it is supported
on $T_{I}$. Tracking these reduction also yields the coefficients
$\left(g_{q}\right)_{q\in\mathcal{Q}}$. Note that since $I$ is finitely
generated, $\mathcal{Q}$ is finite (see Remark \ref{rem: cardinality of combinatorially reducing system}),
which ensures that checking for possible reductions at each step involves
only finitely many options.
\begin{thm}
Let $\mathcal{Q}$ be a combinatorially reducing system for a finitely
generated ideal $I\leq\mathcal{A}$. Then Algorithm \ref{alg:division_with_remainder}
computes the unique remainder $r\in\text{Sp}_{K}\left(T_{I}\right)$
and coefficients $\left(g_{q}\right)_{q\in\mathcal{Q}}$ such that
$f=\sum_{q\in\mathcal{Q}}qg_{q}+r$ and every $g_{q}$ is supported
only on words not beginning with $\text{HTT}\left(q\right)^{-1}$.
\end{thm}
\begin{algorithm}
\caption{\textsc{DivisionWithRemainder} - Divides with Remainder by a Combinatorially Reducing System}
\label{alg:division_with_remainder}
\begin{algorithmic}[1]
\REQUIRE An element $f \in \mathcal{A}$ and a combinatorially reducing system $\mathcal{Q}$ for a finitely generated right ideal $I$.
\ENSURE \parbox[t]{0.91\linewidth}{The unique remainder $r \in \operatorname{Sp}_K(T_I)$ and coefficients $(g_q)_{q \in \mathcal{Q}}$ such that  $f = \sum_{q \in \mathcal{Q}} q g_q + r$, and each $g_q$ is supported only on words not beginning with $\operatorname{HTT}(q)^{-1}$}.
\STATE $r \gets f$
\FORALL{$q \in \mathcal{Q}$}
    \STATE $g_q \gets 0$
\ENDFOR
\WHILE{there exists $u \in \operatorname{supp}(r)$ and $q_u \in \mathcal{Q}$ such that $\mathrm{HT}(q_u)$ is a prefix of $u$}
    \STATE Let $\gamma \in K^{\times}$ be the coefficient of $u$ in $r$
    \STATE $u_{\mathrm{suf}} \gets \mathrm{HT}(q_u)^{-1} u$ \COMMENT{so $u = \mathrm{HT}(q_u) \cdot u_{\mathrm{suf}}$ without cancellation}
    \STATE $r \gets r - \gamma q_u u_{\mathrm{suf}}$ \COMMENT{cancel the term $u$ from $r$}
    \STATE $g_{q_u} \gets g_{q_u} + \gamma u_{\mathrm{suf}}$
\ENDWHILE
\RETURN $\left((g_q)_{q \in \mathcal{Q}}, r\right)$
\end{algorithmic}
\end{algorithm}
\begin{proof}
The procedure for computing $r$ and determining termination is identical
to Algorithm \ref{alg:reduction_mod_TI}, where it was shown that
the process halts after finitely many steps and correctly returns
the remainder $\phi_{I}\left(f\right)\in\text{Sp}_{K}\left(T_{I}\right)$
of $f$ modulo $T_{I}$. The same argument applies here. 

We show that both the support condition on the coefficients $g_{q}$
and the equation $f=\sum_{q\in\mathcal{Q}}qg_{q}+r$ hold at every
step of the algorithm. Initially, $r=f$ and all coefficients $g_{q}$
are zero. Each time a term $\gamma q_{u}u_{\text{suf}}$ is subtracted
from $r$, the term $\gamma u_{\text{suf}}$ is added to $g_{q_{u}}$,
so the total remains unchanged. Furthermore, since $u=\text{HT}\left(q_{u}\right)\cdot u_{\text{suf}}$
without cancellation, the suffix $u_{\text{suf}}$ cannot begin with
$\text{HTT}\left(q\right)^{-1}$, ensuring that the support condition
on $g_{q_{u}}$ is maintained when adding $\gamma u_{\text{suf}}$.
\end{proof}
\begin{rem}
When Algorithm \ref{alg:division_with_remainder} is applied to an
element $f\in I$ as input, the remainder $\phi_{I}\left(f\right)$
is zero. In this case, the algorithm computes the canonical representation
$f=\sum_{q\in\mathcal{Q}}qg_{q}$ described in Lemma \ref{lem: Canonical-Representation-using-1}
with respect to the CRS $\mathcal{Q}$ for $I$.
\end{rem}

\subsubsection{Canonically Expressing $f\in I$ with the Exposure Basis of $I$}

Let $I\leq\mathcal{A}$ be a finitely generated ideal. Denote by $m$
the rank of $I$, and let the exposure basis be $B_{I}=\left\{ f_{0},f_{1},...,f_{m-1}\right\} $,
with the elements ordered in ascending $\prec_{\text{max}}$-order
of supports, as determined in the exposure process. In this subsection,
we describe an algorithm which, given $B_{I}$ and an element $h\in I$,
computes the unique coefficients $p_{i}\in\mathcal{A}$ such that
$h=\sum_{i=0}^{m-1}f_{i}p_{i}$.

Our approach consists of three main steps: First, express the \texttt{seconds}
$s_{0},s_{1},...,s_{m-1}$ in terms of $B_{I}$. This amounts to computing
a matrix $C\in\text{Mat}_{m\times m}\left(\mathcal{A}\right)$ such
that $\left(s_{0},s_{1},...,s_{m-1}\right)=\left(f_{0},f_{1},...,f_{m-1}\right)\cdot C$.
Second, recall that by Theorem \ref{Thm: properties of second + Grobner basis is combinatorially reducing},
$B_{I}^{\text{gr}}=\left\{ f_{i}:0\leq i<m\right\} \cup\left\{ s_{i}:0\leq i<m\right\} $
is a combinatorially reducing system for $I$. We thus use Algorithm
\ref{alg:division_with_remainder} to express $h$ canonically using
$B_{I}^{\text{gr}}$. Third and last, bring this expression to the
desired form by using the matrix $C$.

We begin with the Algorithm \ref{alg:compute_matrix_C}, which computes
the change-of-basis matrix $C$ from the $B_{I}$ to $\left\{ s_{i}:0\leq i<m\right\} $.
The \texttt{seconds} $\left(s_{i}\right)_{0\leq i<m}$ are not required
as input, but computes them from $B_{I}$. The algorithm is slightly
technical, so a proof of its correctness follows.

\begin{algorithm}
\caption{\textsc{ComputeMatrixC} - Computes the Matrix Expressing the \texttt{seconds} Using the Exposure Basis}
\label{alg:compute_matrix_C}
\begin{algorithmic}[1]
\STATE \textbf{Input:} \parbox[t]{0.85\linewidth}{%
    The exposure basis $B_I = \left[ f_0, f_1, \ldots, f_{m-1}\right]$ of a proper ideal $I \lneqq \mathcal{A}$, given in ascending $\prec_{\max}$-order of supports.
}\vspace{4pt}
\STATE \textbf{Output:} \parbox[t]{0.85\linewidth}{%
    The matrix $C \in \mathrm{Mat}_{m \times m}(\mathcal{A})$ such that $(s_0, s_1, \ldots, s_{m-1}) = (f_0, f_1, \ldots, f_{m-1}) \cdot C$. Here, $s_{i}=\mathrm{MONIC}\left(\phi_{I_{i}}\left(f_{i}b_{i}^{-1}\right)\right)$ is the \texttt{seconds} corresponding to $f_i$, where $b_{i}=\mathrm{HTT}(f_i)$ and $I_i$ is the ideal generated by $f_0,f_1,...,f_{i-1}$.
}\vspace{8pt}
\STATE Initialize $C$ to be the $m\times m$ zero matrix.
\FOR{$i = 0$ \TO $m-1$}
    \STATE $b_i \gets \mathrm{HTT}(f_i)$
    \STATE $\mathcal{Q}_i \gets \{f_0, f_1, \ldots, f_{i-1}, s_0, s_1, \ldots, s_{i-1}\}$ \COMMENT{CRS for $I_i=\langle f_0,\ldots,f_{i-1}\rangle$}
    \STATE $(r_i, \{g_{f_j} : j < i\}, \{g_{s_j} : j < i\}) \gets \textsc{DivisionWithRemainder}(f_i b_i^{-1}, \mathcal{Q}_i)$ \COMMENT{$r_i=\phi_{I_{i}}\left(f_{i}b_{i}^{-1}\right)$}
	\STATE $\mu_i \gets $ coefficient of $\mathrm{HT}(r_i)$ in $r_i$.
    \STATE $s_i \gets \mu_i^{-1}r_i$ \COMMENT{$s_i$ is the monic $K^\times$-multiple of $r_i$}
	\STATE $C_{i,i} \gets C_{i,i}+\mu_i^{-1}\cdot b_i$.
    \FOR{$j = 0$ \TO $i-1$}
        \STATE $C_{j,i} \gets C_{j,i}-\mu_i^{-1}g_{f_j}$.
		\STATE $C_{\bullet,i} \gets C_{\bullet,i} - C_{\bullet,j}\cdot \mu_i^{-1}g_{s_j}$. \COMMENT{Right-multiply column $j$ by $\mu_i^{-1}g_{s_j}$ and subtract from column $i$}
    \ENDFOR
\ENDFOR
\RETURN $C$
\end{algorithmic}
\end{algorithm}
\begin{prop}
Algorithm \ref{alg:compute_matrix_C} computes the matrix $C\in\text{Mat}_{m\times m}\left(\mathcal{A}\right)$
satisfying $\left(s_{0},s_{1},...,s_{m-1}\right)=\left(f_{0},f_{1},...,f_{m-1}\right)C$.
\end{prop}
\begin{proof}
We prove by induction on $i$ that the $i$-th column of $C$ satisfies
$s_{i}=\left(f_{0},f_{1},...,f_{m-1}\right)C_{\bullet,i}$. Suppose
this holds for every $0\leq j<i$. By definition, $s_{i}=\text{MONIC}\left(\phi_{I_{i}}\left(f_{i}b_{i}^{-1}\right)\right)$,
where $b_{i}=\text{HTT}\left(f_{i}\right)$ and $I_{i}$ is the ideal
generated by $\left\{ f_{0},f_{1},...f_{i-1}\right\} $. By Theorem
\ref{Thm: properties of second + Grobner basis is combinatorially reducing},
the set $\mathcal{Q}_{i}=\left\{ f_{j}:j<i\right\} \cup\left\{ s_{j}:j<i\right\} $
forms a CRS for $I_{i}$. Dividing $f_{i}b_{i}^{-1}$ with remainder
with respect to $\mathcal{Q}_{i}$, we obtain:
\[
f_{i}b_{i}^{-1}=\sum_{j<i}f_{j}g_{f_{j}}+\sum_{j<i}s_{j}g_{s_{j}}+r_{i},
\]
where $r_{i}=\phi_{I_{i}}\left(f_{i}b_{i}^{-1}\right)$. Let $\mu_{i}\in K^{\times}$
be the coefficient of the head term $\text{HT}\left(r_{i}\right)$
in $r_{i}$. Then $s_{i}=\mu_{i}^{-1}r_{i}$.

Now let $e_{0},e_{2},...,e_{m-1}$ denote the standard basis for the
space of column vectors $\mathcal{A}^{m}$. Initially, $C_{\bullet,i}=0$.
Iteration $i$ updates $C_{\bullet,i}$ so that
\[
C_{\bullet,i}=\mu_{i}^{-1}\left(e_{i}b_{i}^{-1}-\sum_{j<i}e_{j}g_{f_{j}}-\sum_{j<i}C_{\bullet,j}g_{s_{j}}\right).
\]
Using the induction hypothesis, we conclude:
\begin{align*}
\left(f_{0},f_{1},...,f_{m-1}\right)\cdot C_{\bullet,i} & =\mu_{i}^{-1}\left(f_{i}b_{i}^{-1}-\sum_{j<i}f_{j}g_{f_{j}}-\sum_{j<i}s_{j}g_{s_{j}}\right)=\mu_{i}^{-1}r_{i}=s_{i}.
\end{align*}
\end{proof}
We now provide the algorithm for expressing an element $h\in I$ canonically
using the exposure basis $B_{I}$. In addition to $h$ and $B_{I}$,
the algorithm takes the change-of-basis matrix $C\in\text{Mat}_{m\times m}\left(\mathcal{A}\right)$
returned by Algorithm \ref{alg:compute_matrix_C}.

\begin{algorithm}
\caption{\textsc{ExpressWithExposureBasis} - Expresses an Element in an Ideal Using its Exposure Basis}
\label{alg:express_using_exposure}
\begin{algorithmic}[1]
\STATE \textbf{Input:} \parbox[t]{0.90\linewidth}{%
    An element $h\in I$, the exposure basis $B_I = \left[ f_0, f_1, \ldots, f_{m-1}\right]$ of an ideal $I \leq \mathcal{A}$, given in ascending $\prec_{\max}$-order of supports, and the matrix $C \in \mathrm{Mat}_{m \times m}(\mathcal{A})$ such that $(s_0, s_1, \ldots, s_{m-1}) = (f_0, f_1, \ldots, f_{m-1}) \cdot C$ for the \texttt{seconds} $s_0, s_1, \ldots, s_{m-1}$ of $B_I$.
}\vspace{4pt}
\STATE \textbf{Output:} \parbox[t]{0.90\linewidth}{%
    The unique coefficients $p_0,p_1,...,p_{m-1}\in \mathcal{A}$ such that $h=\sum_{i=0}^{m-1}f_{i}p_{i}$.
}\vspace{8pt}
\IF{$B_I = \left[ e \right]$}
	\RETURN $h$ \COMMENT{treat the case $I=\mathcal{A}$ separately}
\ENDIF
\STATE Set $(s_0, s_1, \ldots, s_{m-1}) \gets (f_0, f_1, \ldots, f_{m-1}) \cdot C$ \COMMENT{calculate \texttt{seconds}}
\STATE Set $\mathcal{Q} \gets \{f_0, f_1, \ldots, f_{m-1}, s_0, s_1, \ldots, s_{m-1}\}$ \COMMENT{CRS for $I$}
\STATE $(\{g_{f_j} : j < m\}\cup \{g_{s_j} : j < m\}, r) \gets \textsc{DivisionWithRemainder}(h, \mathcal{Q})$ \COMMENT{$r$ should vanish since $h\in I$}
\STATE Set $g_f \gets \left(g_{f_{0}},g_{f_{1}},...,g_{f_{m-1}}\right)^{T}$ and $g_s \gets \left(g_{s_{0}},g_{s_{1}},...,g_{s_{m-1}}\right)^{T}$ \COMMENT{$v^T$ denotes the transpose of $v$}
\STATE Set $p \gets g_f + C \cdot g_s$.
\RETURN entries of $p$.
\end{algorithmic}
\end{algorithm}
\begin{prop}
The column vector $p\in\mathcal{A}^{m}$ computed by Algorithm \ref{alg:express_using_exposure}
satisfies $h=\sum_{i=0}^{m-1}f_{i}p_{i}$.
\end{prop}
\begin{proof}
Let $f$ and $s$ denote the row vectors $\left(f_{0},f_{2},...,f_{m-1}\right)$
and $\left(s_{0},s_{1},...,s_{m-1}\right)$ respectively. The matrix
$C$ satisfies $s=fC$. Since $h\in I$, we have $\phi_{I}\left(h\right)=0$,
so the division with remainder computed by Algorithm \ref{alg:division_with_remainder}
satisfies
\[
h=\sum_{i=0}^{m-1}f_{i}g_{f_{i}}+\sum_{i=0}^{m-1}s_{i}g_{s_{i}}=fg_{f}+fCg_{s}=f\left(g_{f}+Cg_{s}\right)=fp.
\]
\end{proof}
\bibliographystyle{plain}
\phantomsection\addcontentsline{toc}{section}{\refname}\bibliography{Exposure_Orders_in_Free_Group_Algebras}

\end{document}